%% file: cs.tex
\documentclass[11pt]{article}
\input{math}
\input{common}
\input{symbols}
\usepackage[style=alphabetic,sorting=none]{biblatex}
\usepackage[margin=1in]{geometry}
\usepackage{stmaryrd}
\DeclareFieldFormat{sentencecase}{\MakeSentenceCase{#1}}

\renewbibmacro*{title}{%
  \ifthenelse{\iffieldundef{title}\AND\iffieldundef{subtitle}}
    {}
    {\ifthenelse{\ifentrytype{article}\OR\ifentrytype{inbook}%
      \OR\ifentrytype{incollection}\OR\ifentrytype{inproceedings}%
      \OR\ifentrytype{inreference}}
      {\printtext[title]{%
        \printfield[sentencecase]{title}%
        \setunit{\subtitlepunct}%
        \printfield[sentencecase]{subtitle}}}%
      {\printtext[title]{%
        \printfield[titlecase]{title}%
        \setunit{\subtitlepunct}%
        \printfield[titlecase]{subtitle}}}%
     \newunit}%
  \printfield{titleaddon}}

\usepackage{authblk}
\author[1]{Mark Fornace}
\author[1,2]{Michael Lindsey}
\affil[1]{Lawrence Berkeley National Laboratory}
\affil[2]{University of California, Berkeley}
\title{Column and row subset selection using nuclear scores: algorithms and theory for Nystr\"{o}m approximation, CUR decomposition, and graph Laplacian reduction}

\begin{document}


\maketitle

\begin{abstract}
    Column selection is an essential tool for structure-preserving low-rank approximation, with wide-ranging applications across many fields, such as data science, machine learning, and theoretical chemistry.
    In this work, we develop unified methodologies for fast, efficient, and theoretically guaranteed column selection.
    First we derive and implement a sparsity-exploiting deterministic algorithm applicable to tasks including kernel approximation and CUR decomposition.
    Next, we develop a matrix-free formalism relying on a randomization scheme satisfying guaranteed concentration bounds, 
    applying this construction both to CUR decomposition and to the approximation of matrix functions of graph Laplacians. Importantly, the randomization is only relevant for the computation of the scores that we use for column selection, not the selection itself given these scores.
    For both deterministic and matrix-free algorithms, we bound the performance favorably relative to the expected performance of determinantal point process (DPP) sampling and, in select scenarios, that of exactly optimal subset selection. 
    The general case requires new analysis of the DPP expectation.
    Finally, we demonstrate strong real-world performance of our algorithms on a diverse set of example approximation tasks.
\end{abstract}
\vspace{0.3in}


\section{Introduction}

Low-rank matrix approximation is a ubiquitous task in computational modeling, machine learning, and data analysis (e.g. \cite{Markovsky2012-hn,Halko2011-md}).
Beyond obvious applications in data compression, low-rank approximations may be used for preconditioning, noise elimination, experimental design, and model reduction \cite{Tropp2023-mi,Chaturantabut2010-wm,Benner2015-hj,Schmid2010-oy,Williams2000-ah}.
Low-rank matrix approximation is canonically performed using spectral methods, i.e., the singular value decomposition (SVD) or eigendecomposition \cite{Klema1980-sj,Abdi2010-fh}.
While exact implementations of such methods scale cubically with respect to the input matrix dimensions, recent research in iterative and randomized algorithms has expanded the scope of rigorous and accurate approximation to larger and larger systems \cite{Martinsson2020-gm}.
However, the results of such methods generally lose the mathematical structure of the input matrix, may be less interpretable by the user, and may not be suitable for downstream applications \cite{Mahoney2009-et,Chaturantabut2010-wm,Martinsson2020-gm}.

In contrast to spectral methods, column selection (or its transpose, row selection) can be used to generate a low-rank approximation which naturally preserves the mathematical structure (e.g., sparsity, positivity, etc.) of the original matrix.
Moreover, the result of a factorization based on column and/or row selection is readily interpretable and may be uniquely suitable for downstream use.
Already, column selection techniques have proved to be impactful across a variety of problems including kernel approximation \cite{Abedsoltan2023-ou}, regression \cite{Zhang2020-zl}, clustering \cite{Wang2019-xl,Tremblay2020-gk}, quadrature \cite{Epperly2023-bj}, and nonlinear model reduction \cite{Chaturantabut2010-wm}. 
On the other hand, optimal column selection is a difficult combinatorial optimization problem that is assumed to be NP-hard \cite{Civril2014-hy,Civril2013-ll}.
As a result, research in this area has revolved around \emph{approximate selection}, in which the selected column indices are suboptimal but nonetheless close to optimal in theoretical and empirical performance.

In this work, we develop methods based on \emph{nuclear maximization}, in which an objective based on the nuclear or trace norm is directly optimized via greedy augmentation of a growing set of selected indices.
While such a strategy can be more expensive than other selection methods, we demonstrate that its computational cost is actually relatively modest in a variety of scenarios.
Moreover, we are able to develop strong theoretical guarantees, and we validate both the scalability and the accuracy of this selection strategy in practice. Interestingly, such guarantees are known not to hold, e.g., for commonly used greedy augmentation strategies \cite{Higham1990-yp,Civril2009-kt}, which in fact has motivated recent work on randomized index selection \cite{Chen2022-lp} for Nystr\"{o}m approximation.

In addition, we find it to be of particular interest to design our algorithms to be \emph{matrix-free}. 
In such an algorithm, an operator $A$ is accessed only implicitly via matrix vector products (\emph{matvecs}) of the form $A x$ (for arbitrary vector $x$).
Such algorithms are ideally suited to situations in which the full $A$ is infeasible to calculate or store in memory, or in which $A$ is viewed as a matrix function of a sparse matrix and query access to individual entries of $A$ is inefficient.
An elementary motivating example is given by taking $A = B^\Mo$, where $B$ is a sparse and well-conditioned matrix.
While $A$ is often too large to be computed or stored ahead of time (being generally dense), $A x = B^\Mo x$ may be computed accurately and quickly using standard iterative methods.
In general, a matrix-free algorithm allows for flexible and efficient exploitation of operator structure across a host of diverse scenarios.
To our knowledge, our selection algorithms are the first designed explicitly for such a regime.

The deterministic algorithms that we introduce require exact calculations of the diagonals of certain intermediate matrices. 
In some contexts, these quantities can be calculated with reasonable computational cost.
In other contexts, and in a general matrix-free setting, such calculations may impose severe bottlenecks.
As a result, we use a randomized diagonal estimation scheme which can be viewed as an application of stochastic trace estimation techiques \cite{Meyer2021-pn,Vershynin2009-jx} or Johnson-Lindenstrauss subspace embedding \cite{Johnson1984-kh}. 
The estimator can be used to compute the nuclear scores that we require for column selection with constant relative error at a cost growing nearly optimally with dimension. In addition to our theoretical error bound, we find the matrix-free approach to be effective in practice across a wide-ranging set of tasks.

A problem of particular interest to us is column selection for inverse graph Laplacians.%
\footnote{Actually, we shall focus on recurrent Laplacians containing a single stationary mode and handle the ``inverse'' of this mode in an appropriate limiting sense.}
In this setting, we are given a (generally sparse) graph Laplacian system and seek to build a \Nystrom{} approximation which models its inverse well.
Even without the capacity to compute the dense inverse operator, our matrix-free approach allows us to leverage recent fast Laplacian solvers (e.g., \cite{Kyng2016-qv,Cohen2016-an,Chen2020-ya}) to choose informative columns in an efficient and theoretically controlled manner.
This problem has applications to model order reduction of Markov processes \cite{Aldous1995-bg} that we intend to explore in future work.
Existing approaches (e.g., \cite{Higham1990-yp,Chen2022-lp}) based on entrywise queries would be infeasible to perform at scale. 
Meanwhile, the inverse Laplacian setting is of special interest since we are able to prove that our chosen objective is strictly submodular on non-empty sets, allowing us to bound performance relative to optimal column subset selection.

For all of our algorithms, we present novel theoretical error bounds and proofs.
In particular, we leverage linear programming bounds of the sort developed in the submodularity literature \cite{Nemhauser1978-ci} to bound the results of our nuclear maximization scheme relative to those of determinantal point process sampling (DPP) \cite{Derezinski2021-hi} and, in turn, spectral methods. 
With small modifications, these proofs easily accommodate our randomized methods for computing the necessary diagonal scores, via the aforementioned concentration bounds controlling their relative error.
Finally, for the particular task of inverse graph Laplacian column selection, we are able to prove the strict submodularity of our nuclear objective, which bounds our algorithm relative to optimal column subset selection.
Our theoretical results are supported by strong empirical performance on a diverse set of example studies.

\subsection{Past work \label{s:past-work}}

Due to its central role in low-rank approximation, the column subset selection problem (CSSP) has been the subject of many research efforts over the past twenty years. See, e.g., \cite{Drineas2006-ak,Frieze2004-pp,Smola2000-si,Smola2000-zq} for early treatments of this problem, as well as \cite{Kireeva2024-if,Kishore_Kumar2017-ls,Martinsson2020-gm,Tropp2023-mi} for useful references.
While it is intractable to offer a comprehensive survey of the literature in this manuscript, we will highlight some of the main thrusts.
Later, in Section~\ref{s:related-methods}, we will describe in more depth the column selection strategies which are most relevant to our work and which we compare against in numerical experiments.

\paragraph{Spectral methods:}
The singular value decomposition (SVD) forms the baseline for research and comparison into low-rank approximations of matrices \cite{Klema1980-sj}.
Classical results (e.g. \cite{Gu2015-uv,Golub2013-an}, Section~\ref{s:spsd-objective}) show that the truncated $k$-SVD achieve the optimal Frobenius error over all rank-$k$ factorizations of a given rectangular matrix, and likewise, that the truncated $k$-eigendecomposition achieves the optimal nuclear norm error over all rank-$k$ factorizations of a symmetric matrix.
Thus the SVD is a ubiquitous tool in feature extraction and data analysis, often referred to as principal component analysis (PCA) within this context \cite{Abdi2010-fh}.
A great success of randomized methods has been in the development of accurate and fast algorithms for truncated SVD computation (e.g. \cite{Halko2011-md,Musco2015-fb,Kaloorazi2018-gd}).
(In this work, we use subspace iteration \cite{Gu2015-uv,Demmel1997-sp} to benchmark our results.)
Given the success of randomized SVD methods, numerous CSSP algorithms compute a truncated (often randomized) SVD as an initial preprocessing step before column selection takes place (e.g. \cite{Avron2013-od,Boutsidis2009-tt,Boutsidis2014-zn,Li2010-gj}). By contrast, our algorithms will not rely on spectral methods, except for the condition number estimation procedure that we apply to our preconditioned graph Laplacians in Section~\ref{s:sqrt-alg}.

\paragraph{Leverage score methods:}
Another major category of column selection approaches is based on independently sampling columns according to leverage scores \cite{Drineas2008-tp,Mahoney2009-et, Drineas2011-nz, Papailiopoulos2014-wq, Alaoui2015-cn, Calandriello2017-zy, Musco2017-fb, Rudi2017-ik, Rudi2018-yp}. The naive cost of computing leverage scores is bottlenecked by the cost of computing a truncated SVD, though we point out that recent work \cite{Cohen2017-st} has explored alternative approaches based on ridge leverage scores, for example, achieving input-sparsity complexity. Leverage scores have also been explored as a tool for tensor factorization \cite{Woodruff2022-wt, Bharadwaj2023-jn}. While leverage score sampling permits transparent theoretical guarantees, requiring $O(k \log k)$ columns for a target rank of $k$, the major practical drawback is that the selection procedure is non-adaptive in the sense that columns are selected independently of one another, yielding subsets that are typically far less compact than alternative, possibly heuristic, approaches. It is worth pointing out that leverage score sampling can be viewed as a preprocessing step for a more refined selection procedure \cite{Derezinski2021-hi}.

\paragraph{Pivoted decompositions:}
Many column selection strategies may be viewed as \emph{pivoting} routines.
In particular, the QR decomposition with column pivoting (QRCP) \cite{Businger1971-gy} is a workhorse tool, and more powerful, albeit more expensive, \emph{strong} rank-revealing QR (RRQR) factorizations have also been introduced \cite{Gu1996-wz}. Contemporary applications range from quantum chemistry \cite{Damle2015-lv} to randomized projection preprocessing \cite{Voronin2017-ct}. 
Meanwhile, the LU decomposition with column pivoting forms the basis for another group of algorithms
\cite{Chen2020-de,Ekenta2022-qb,Pan2000-rv,Boutsidis2009-tt}, and 
the discrete empirical interpolation method (DEIM) \cite{Chaturantabut2010-wm} is based on a selection strategy that is related to the pivoted LU decomposition \cite{Sorensen2016-is}. 
Finally, pivoted Cholesky factorizations \cite{Gu2004-eg,Chen2022-lp} can be used to implement symmetric positive semidefinite (SPSD) \Nystrom{} approximation.
In fact, there is an established connection \cite{Higham1990-yp} between the pivoted QR decomposition and the pivoted Cholesky factorization of the induced kernel matrix in the sense of Section~\ref{s:cur-objective}. In particular, this relates standard QRCP with the method we refer to as diagonal maximization in Section~\ref{s:related-methods} and beyond.

\paragraph{Additional sampling approaches:}
Additional column selection methods include techniques based on on Markov chain Monte Carlo (MCMC) sampling (e.g. \cite{Belabbas2009-il}), relaxation to a continuous optimization problem \cite{Mathur2023-az}, and local determinant-maximizing swaps \cite{Fanuel2021-oc,Anari2016-ai,Damle2024-sn}.
Some of algorithms introduce randomized pre-sketches of input matrices \cite{Musco2017-nf,Martinsson2016-sy}.
One may consult \cite{Kumar2012-xk,De_Hoog2011-jq,Cohen2015-rq,Cortinovis2020-me,Kumar2009-fa} for discussions of non-adaptive and adaptive scoring metrics.
Finally, a few recent works have targeted the extension of column selection to the streaming regime \cite{Jiang2021-wf,Woodruff2023-qp,Tropp2017-kn}.

\paragraph{Experimental design:}
Outside of the field of numerical linear algebra, column selection has a rather long history in the field of experimental design (e.g. \cite{Pukelsheim2006-vh,Krause2008-iu}).
In this context, given a covariance kernel $K$ (e.g., from a Gaussian process or specified information matrix), in which each row/column index corresponds the output of a particular experiment or sensor, one seeks to determine the most informative subsequent experiments to conduct or positions in which to place sensors.
A number of heuristics for this task have been proposed \cite{Boyd2004-yv,Pukelsheim2006-vh}, including
(1) $A$-optimality, maximizing the trace of $K_{\I, \I}$,
(2) $D$-optimality, maximizing the determinant of $K_{\I, \I}$,
(3) $E$-optimality, maximizing the minimum eigenvalue of $K_{\I, \I}$, and 
(4) $I$-optimality, maximizing the average conditional variance $\sum_{\j} K_{\j,\j} - K_{\j,\I} (K_{\I,\I})^\Mo K_{\I,\j}$.
From this point of view, our work can be viewed as targeting $I$-optimality.
Additional objective formulations include multi-objective optimization via Pareto optimality \cite{Qian2020-xf} and ``fair'' objectives \cite{Matakos2023-fb}.

\paragraph{Submodularity:}
A main thread in the experimental design and sensor placement literature concerns the optimization of \emph{submodular} functions (Section~\ref{s:submodularity}).
Roughly speaking, if an objective can be proven to be submodular, then a greedy approach is well-bounded relative to optimal subset selection \cite{Nemhauser1978-ci}.
This observation has led to research on submodular objectives, including quite prominently \cite{Krause2008-iu}.
Unfortunately, many relevant objectives are not in general submodular.
As a result, numerous researchers have focused on the notion of \emph{weak submodularity}, using a theoretical or empirical bound on the degree to which submodularity is violated in order to bound greedy approximation quality \cite{Bian2017-se,Bogunovic2018-pm,Das2011-ub,Elenberg2018-rf,Harshaw2019-gj,Khanna2017-dl}.
Relevant works include those of \cite{Harshaw2019-gj,Hashemi2018-fw,Sun2006-xj}, where specific formulations for covariance matrices including diagonal noise terms are bounded using notions of weak submodularity.

\subsection{Alternative and related methods \label{s:related-methods}}
Here we describe in greater depth the methods which are most relevant to our work.
Except for uniform sampling, we consider only adaptive approaches here.%
\footnote{We refer the reader for basics to the recent review \cite{Kireeva2024-if}, although not all methods we discuss are mentioned therein.}

\subsubsection{Uniform sampling}
A useful benchmark is \emph{uniform sampling}, in which each index subset is generated (without replacement) with equal probability (e.g. \cite{Williams2000-ah}).
While extremely simple and cheap to implement, the quality of approximation achieved is poorly bounded in general, and results on practical problems are often unsatisfactory.
A large body of literature, however, studies this selection strategy as related to the concept of matrix \emph{coherence} \cite{Candes2012-yq}.
For matrices with low coherence, uniform sampling has strongly bounded error and can be effective for sketching, preconditioning, and approximation applications.
This strategy is popular for many kernel approximation tasks \cite{Wang2019-xl,Gittens2011-fd,Gittens2013-kl,Kumar2009-wd} but has also been studied for CUR decomposition \cite{Cai2021-km,Hamm2020-xf,Chiu2013-wc} and as a preprocessing step for other selection methods \cite{Khanna2017-dl}.

\subsubsection{Diagonal maximization}
For \Nystrom{} approximation of a SPSD matrix $K$, the popular \emph{diagonal maximization} algorithm selects the index corresponding to the largest diagonal element, followed adaptively by the indices corresponding to the largest diagonal elements in a sequence of updated matrices formed by eliminating prior choices \cite{Higham1990-yp}.
This approach may be viewed as a greedy maximization of the submatrix determinant, since for a fixed index subset $\I$ and a new disjoint index $\j$ to be selected, we have that
$
    \logdet(K_{\Ij,\Ij}) = \logdet(K_{\I,\I}) + \log(K_{\j,\j} - K_{\j,\I} (K_{\I,\I})^\Mo K_{\I, \j})
    \req{logdet}
$ 
\cite{Abadir2005-wf}.
This strategy is common in practice and has been studied and extended by numerous authors 
\cite{Anari2016-ai,Belabbas2009-il,Bhaskara2023-mq,Bach2005-lx,Civril2009-kt,Deshpande2006-do,Deshpande2006-uc,Dong2023-gg,Gu2004-eg}. It is often simply called the ``greedy'' approach (cf. \cite{Chen2022-lp}), but we shall not refer to it as such in order to avoid confusion with our methods.

\subsubsection{Diagonal sampling}
\cite{Chen2022-lp} have recently introduced the randomly pivoted Cholesky (RPCholesky) method for kernel approximation, which is algorithmically the same as diagonal maximization, except that each subsequent index is chosen randomly according to weights given by the appropriate diagonal element magnitudes.
We shall refer to this approach as \emph{diagonal sampling}. Note that this method is induced by the randomly pivoted QR factorization \cite{Deshpande2006-do,Deshpande2006-uc} via the aforementioned connection between pivoted QR and pivoted Cholesky factorizations \cite{Higham1990-yp}. 
In some sense, this is an intermediate approach between uniform sampling (all indices of equal probability to be chosen) and diagonal maximization (only the index corresponding the largest diagonal element to be chosen).
In \cite{Chen2022-lp} it is shown that this approach is superior to diagonal maximization in its theoretical worst case and for some example problems; followup work has discussed different weighting metrics \cite{Steinerberger2024-pp} and applications \cite{Diaz2023-kl}.

\subsubsection{Determinantal point process sampling}
In \emph{determinantal point process} (DPP) sampling (e.g., \cite{Kulesza2012-pp,Derezinski2021-hi,Gillenwater2014-to}), or volume sampling, a subset of columns from a symmetric positive semidefinite (SPSD) square matrix $K$ is generated according to a probability density in which the probability of a given subset $\I$ is proportional to the principal minor $\det(K_{\I, \I})$.
The theory behind such sampling is quite rich \cite{Guruswami2012-rr,Derezinski2020-fb}.
In particular, the expectation of DPP sampling is independent of the input basis and can be derived purely from the spectrum of input $K$ \cite{Guruswami2012-rr}.
Building on this connection, \textcite{Guruswami2012-rr,Derezinski2020-fb,Belabbas2009-il} showed that DPP sampling is close to worst-case optimal and achieves strong theoretical guarantees in expectation.
On the other hand, while DPP sampling enjoys attractive theoretical properties, it is relatively complicated to perform algorithmically in practice.
Exact sampling may be performed using eigendecomposition of the input matrix $K$, incurring $O(n^3)$ cost \cite{Derezinski2021-hi}.
To treat larger systems, approximate sampling methods based on Markov chain Monte Carlo (MCMC) and approximate scores have been developed \cite{Calandriello2020-kz,Derezinski2019-zj,Gautier2019-ml,Poulson2020-zc,Deshpande2010-uy}.
In our view, DPP represents a gold standard for theoretical bounds but remains relatively difficult to scale to large problems.
In this work, we use DPP sampling exclusively as a theoretical tool rather than as an algorithm, meaning that its practical implementation is not of immediate concern for us.

\subsubsection{Nuclear maximization}
We study, extend, and invent greedy algorithms for \emph{nuclear maximization} in the current work.
From our reading, this paradigm has received less attention in the literature, perhaps due to a perceived increase in the computational complexity. This approach has been implicitly considered in works pertaining to both \Nystrom{} approximation and CUR factorization, but via the reduction of CUR factorization to \Nystrom{} approximation reviewed below in Section~\ref{s:cur-objective}, we will review them in a unified context.
The theory of this approach is studied most in depth in \cite{Altschuler2016-gd}, in which a notion of weak submodularity is leveraged to bound approximation accuracy. The content and proof of the bounds within this reference are very different than ours, involving an \emph{a posteriori} bound on the minimal singular value of the optimal $k$-subset which we believe to be less useful. Moreover, the algorithmic details of this work, which pertains to CUR decomposition, are quite different; in particular, sketching is used as a preprocessing step to maintain computational efficiency. We believe that our approaches to CUR factorization are more practical, and in particular they can naturally exploit sparsity.
Aside from that reference, this general scheme has been studied or mentioned (in varying levels of detail) in \cite{Farahat2011-zm,Farahat2015-zl,Derezinski2020-fb,Ordozgoiti2018-wn,Ordozgoiti2019-xz}.
In particular, we believe that \cite{Farahat2011-zm} and \cite{Ordozgoiti2018-wn} give methods which match our Algorithm~\ref{alg:exact-cholesky} in time and space complexities for dense matrices (though not for sparse matrices).
In the sensor placement literature, a related algorithm involves uniformly sampling subsets of candidate sensor positions followed by nuclear maximization within each drawn subset \cite{Hashemi2018-fw,Hibbard2023-fg}.

\subsubsection{Nomenclature}
To avoid confusion due to under-specification of the terms ``greedy'' and ``random,'' we adopt in this paper the following nomenclature for the selection of $k$ column indices $\I$ out of a SPSD matrix $K$:
\begin{enumerate}
    \item \emph{Uniform sampling} refers to sampling columns from $K$ with equal probabilities (without replacement).
    \item \emph{Diagonal maximization} refers to the (deterministic) greedy maximization of $\det(K_{\I,\I})$. In practice this means picking the next element $\j$ to maximize $K_{\j,\j} - K_{\j,\I} (K_{\I,\I})^\Mo K_{\I, \j}$.
    \item \emph{Diagonal sampling} refers to modifying diagonal maximization such that each additional element is chosen instead with probability proportional to $K_{\j,\j} - K_{\j,\I} (K_{\I,\I})^\Mo K_{\I, \j}$.
    \item \emph{Nuclear maximization} refers to the (deterministic) greedy maximization of the nuclear norm $\Tr[(K^2)_{\I,\I} (K_{\I,\I})^\Mo]$. This is the basic method considered in our work.
    \item \emph{DPP sampling} refers to the sampling of subsets $\I$ with probability proportional to $\det(K_{\I,\I})$. 
    (This is sometimes otherwise labeled as volume sampling.)
    We denote DPP sampling over all subsets of size $k$ as $k$-DPP sampling.
\end{enumerate}
In the context of CUR decomposition of a matrix $A$, the above conventions apply to $K=A A^\t$ or $K=A^\t A$, in the cases of row or column selection, respectively.
Note that we do not foresee any advantage to performing randomized index selection based on the nuclear scores, so this possibility is omitted.
In this work, randomization is used solely to approximate the nuclear scores for subsequent deterministic maximization.

\subsection{Outline of problems and contributions}
While we develop methods which are applicable to a variety of general scenarios, in this work we are particularly focused on three.
In Section~\ref{s:objectives}, we describe our mathematical formulation for each of these problems and how, in a basic sense, they may be reduced to the problem of column selection for a symmetric positive semidefinite (SPSD) matrix $K$: 
\begin{enumerate}
    \item \emph{Kernel approximation}: 
    In this use-case, we are presented \emph{a priori} with a (possibly sparse) SPSD matrix $K$ which supports $O(1)$ element query access.
    Our task is to form a low-rank \Nystrom{} approximation of $K$ using a subset of its columns \cite{Drineas2005-dh}.
    The output column selection may be used for data compression, preconditioning, or experimental design.
    \item \emph{CUR decomposition}:
    In this application, we are presented an arbitrary rectangular matrix $A$, and we would like to compute an approximate factorization of it using subsets of its rows and columns \cite{Goreinov1997-yi,Mahoney2009-et}.
    The resultant factorization may be used for data analysis, interpretation, and other downstream use.
    \item \emph{Inverse graph Laplacian rank reduction}:
    Finally, in this case we are presented with a (typically sparse) graph Laplacian corresponding to a recurrent and reversible random walk \cite{Aldous1995-bg,Merris1994-aw}.
    The goal is to minimize the trace of the Schur complement upon removing select nodes from the graph.
    We intend to explore further applications of this computation in upcoming work.
\end{enumerate}

In Section~\ref{s:exact-alg}, we motivate and describe deterministic algorithms for column selection.
We first detail the implementation and computational complexity for a given, possibly sparse, SPSD matrix $K$. In this case, our algorithm matches a past treatment \cite{Farahat2011-zm} for dense matrices, though our exploitation of sparsity is novel to the best of our knowledge.
Then we describe the novel extension of this algorithm to inverse Laplacian column selection, which requires us to accommodate the null space due to the stationary state.
These algorithms require computations of row norms and diagonal entries and are not fully matrix-free.

In Section~\ref{s:matrix-free}, we introduce our fully matrix-free approaches based on randomized diagonal estimation.
We offer efficient specialized approaches to this estimation in the special cases of sparse CUR decomposition and inverse graph Laplacian compression. 
More generally, the approach can be applied to any operator permitting certain efficient matvec implementations.
First, we describe a matrix-free algorithm for the usual objective given an arbitrary positive semidefinite operator.
Then we extend this treatment to the inverse Laplacian case where, again, we modify the algorithm to accommodate the stationary state. 
Finally, we perform a theoretical analysis to determine and bound the accuracy of our approaches for estimating the diagonals that we require.
While we focus on a particular column selection algorithm, namely nuclear maximization, our approach actually allows for matrix-free extensions of existing approaches such as diagonal maximization and diagonal sampling, which are novel to the best of our knowledge and which we implement for comparison in our numerical experiments. 

In Section~\ref{s:error-analysis}, we perform a thorough error analysis of both our deterministic and matrix-free algorithms, treating both cases simultaneously. 
We develop a generalization of the linear programs used in submodular function maximization (Theorem~\ref{th:submodular}) to show that (1) the objective value achieved by our approach is relatively tightly bounded relative to the expectation of DPP sampling (Theorem~\ref{th:submodular-dpp} and Corollary~\ref{cor:dpp-greedy-bound}), (2) that we can bound our algorithms' error relative to the spectrum of the input operator $K$ (Theorem~\ref{th:greedy-general-bound}), and (3) that we can prove the submodularity of our inverse Laplacian objective (Theorem~\ref{th:laplacian-is-submodular}) to bound the quality of our selected columns relative to that of the exactly optimal set of columns (Theorem~\ref{th:laplacian-trace-bound}).
To the best of our knowledge, all of this analysis is novel except where explicitly noted. In particular, our proofs are generally unrelated to past approaches based on weak submodularity \cite{Altschuler2016-gd,Bian2017-se}. Such approaches bound the degree of violation of submodularity, necessitating an \emph{a posteriori} bound on the minimum singular value of the selected submatrix, cf. \cite{Altschuler2016-gd}, for example.

Finally, in Section~\ref{s:benchmarks}, we present the results of computational experiments in kernel approximation, CUR decomposition, and inverse Laplacian rank reduction.
We first construct illustrative examples to show when alternative approaches may fail to produce accurate approximations.
Then we run our algorithms on a diverse set of more realistic problems, culminating in our evaluation of the discussed column selection methods for inverse Laplacians with $n>10^6$.
Overall, our empirical studies support the quality of our devised methodology. 
In no case that we consider does an investigated alternative appear superior. 
While we observe that diagonal maximization often performs well on many of our real-world examples, producing qualitatively similar performance to nuclear maximization, we can easily pose examples where diagonal maximization is ineffective. 
Indeed, our method possesses substantially stronger theoretical guarantees in general.
Thus we conclude that the empirical results support the accuracy, reliability, and trustworthiness of our algorithms.

For concision, we defer pseudocode, proofs, and elaborations to the following appendices:
\begin{itemize}
    \item Appendix~\ref{a:alg-details} includes pseudocode and details for our algorithms.
    \item Appendix~\ref{a:objective-formulation} contains elementary proofs and further explanations of our objectives.
    \item Appendix~\ref{a:dpp-proofs} contains our proofs related to DPP sampling, including what we believe to be a novel proof of its $k$-concavity; these proofs might be of independent theoretical interest.
    \item Appendix~\ref{a:error-bounds} contains proofs and extended results for our theoretical error bounds.
\end{itemize}

Overall, we advocate our methods on the grounds of their compelling theoretical guarantees, computational cost, and empirical performance, while highlighting that these competing desiderata are otherwise difficult to maintain simultaneously.
Table~\ref{tab:complexities} summarizes the computational complexities of our algorithms.
Efficient and parallelized implementations of our algorithms are available \href{https://github.com/mfornace/nuclear-score-maximization}{online}.

\begin{table}
    \centering
    \begin{tabular}{lll}
    \toprule
    Algorithm setting          & Deterministic scoring & Matrix-free scoring \\
    \midrule
    Kernel approximation                 &   $O(n k^2 + k \nnz(K))$  &  $O(n k^2 z + k z \nnz(K))$  \\
    CUR decomposition  &   $O((m + n) k^2 + k \nnz(A) + \chi)$ & $O((m + n) k^2 z + k z \nnz(A))$        \\
    Laplacian reduction        &   $O(n^2 k + \vartheta )$       &       $O(n k^2 z) + \tilde{O} (k z \nnz(L))$        \\
    \bottomrule
    \end{tabular}
    \caption{
        Computational complexities of proposed algorithms for \Nystrom{} approximation of an $n \times n$ kernel matrix $K$, CUR decomposition of an $m \times n$ matrix $A$, and rank reduction of a rescaled Laplacian $L$.
        Here $k$ is the rank of the approximation or number of selected columns/rows, $\chi$ is the cost of computing both $A A^\t$ and $A^\t A$,  $\vartheta$ is the cost of computing $L^+$, and $z$ is the number of random vectors used in randomized diagonal approximation. 
        To maintain a constant relative scoring error $\epsilon$ with high probability, $z$ should be chosen to scale logarithmically with $\max(m,n)$, cf. Section~\ref{s:concentration}.
        It is assumed that a sparse approximate factorization of $L$ is provided with $\tilde{O}(\nnz(L))$ nonzeros yielding an $O(1)$ preconditioned condition number. Here $\tilde{O}$ indicates the omission of log factors. 
        See Section~\ref{s:complexity-details} for an explanation of each entry.
        \label{tab:complexities}
    }
\end{table}

\section{Formulation of objectives \label{s:objectives}}

In this section, we formulate the nuclear norm maximization which motivates all of our developed algorithms.
We then describe the mathematical formulations of three application problems.
In general, we shall find that each problem is reducible (with mild caveats) to the same problem of column selection from a symmetric positive semidefinite operator:
\begin{enumerate}
    \item In the context of \emph{kernel approximation}, a SPSD matrix $K$ may be supplied directly.
    For instance, for a set of points $\SetDef{x_i}{i=1, \dots, n}$ in $\RR^d$, a typical Gaussian process treatment might explicitly construct a kernel matrix $K$ with $K_{i,j} \propto \exp(-\alpha \norm{x_i - x_j}^2_2)$. 
    \item In the context of \emph{CUR decomposition}, an arbitrary rectangular matrix $A$ is supplied.
    Column selection on $K=A A^\t$ and $K = A^\t A$ may be used to yield a rigorously bounded factorization.
    \item In the context of \emph{inverse Laplacian reduction}, we are given a SPSD rescaled Laplacian $L$ and associated null eigenvector $\h$ corresponding to a recurrent and reversible Markov chain. 
    In this scenario, we effectively perform column selection on the operator $(L + \Shift h h^\t)^\Mo$ in the $\Shift \rightarrow 0$ limit. Detailed motivation for this case is discussed in Appendix~\ref{s:markov}.
\end{enumerate}
While these are the only applications we discuss in this work, we expect that there may be additional applications of our methodologies beyond those listed above.

\subsection{Column selection from a symmetric positive semidefinite matrix \label{s:spsd-objective}}

We begin with the \Nystrom{} approximation of symmetric positive semidefinite (SPSD) matrix $K$ formed from selecting a subset of (disjoint) column indices $\I$:
\begin{eqn}
    K \approx K_{:,\I} (K_{\I,\I})^\Mo K_{\I,:} 
\end{eqn}
In this work, we will generally assume for simplicity that all explicitly considered submatrices $K_{\I,\I}$ are non-singular, but the \Nystrom{} approximation generalizes naturally by considering pseudoinverses. 
Straightforwardly, $K_{:,\I} (K_{\I,\I})^\Mo K_{\I,:} \leq K$ in the Loewner order, and so the Schur complement
\begin{eqn}
\req{Ktilde0}
    \tK(\I) \Eq K - K_{:,\I} (K_{\I,\I})^\Mo K_{\I,:}
\end{eqn}
is also symmetric and positive semidefinite.%
\footnote{In a Gaussian process interpretation, $\tK(\I)$ is the conditional covariance matrix upon measurement of outputs at indices $\I$.}
A natural way then of quantifying the approximation error of subset $\I$ is via the trace of the remainder $\tK(\I)$:%
\begin{eqn}
    \mathcal{E}_K(\I) &\Eq \Tr[\tK(\I)] = \Tr[K] - \Tr[K_{:,\I} (K_{\I, \I})^\Mo K_{\I, :}] 
    \req{error-def}
\end{eqn}
Since the remainder is SPSD, this trace coincides with its \emph{nuclear norm}, which trivially upper bounds the spectral norm.
Therefore, we may define the following non-negative trace objective, which we seek to maximize:
\begin{eqn}
    \Lk_K(\I) \Eq \Tr[K_{:,\I} (K_{\I, \I})^\Mo K_{\I, :}] = 
    \Tr[(K^2)_{\I, \I} (K_{\I, \I})^\Mo].
    \req{objective-def}
\end{eqn}
Note that we have used the cyclic invariance of the matrix trace on the right-hand side.

Let $\Psi_k^n$
denote the set of all $k$-subsets of indices $\{ 1, \dots, n \}$.
Then we may naturally define an \emph{optimal column subset} $\Opt_k$ via: 
\begin{eqn}
    \Lk_K(\Opt_k) \Eq \max_{\I \subset \Psi^n_k} \Tr[(K^2)_{\I, \I} (K_{\I, \I})^\Mo].
    \req{optimal-def}
\end{eqn}
Although each of the matrices $(K^2)_{\I, \I}$ and $K_{\I,\I}$ are SPSD matrices of only size $k \times k$, making evaluation cheap for small $k$, exact optimization of (\ref{eq:optimal-def}) requires enumeration of $\Psi_k^n$ and is infeasible for all but the most trivial problems (e.g. \cite{Civril2014-hy}).
To upper bound the optimal approximation, it is useful to observe that
\begin{eqn}
    \Lk_K(\Opt_k) &= \max_{\I \in \Psi_k^n} \Tr[(\Id^\t_{:,\I} K^2 \Id_{:,\I}) (\Id_{:,\I}^\t K \Id_{:,\I})^\Mo],
    \req{id-def}
\end{eqn}
where $\Id$ denotes the identity matrix of size $n \times n$.
By comparison, the sum of the largest $k$ eigenvalues of $K$ may be posed as the unconstrained maximization problem\footnote{See Proposition~\ref{prop:ky-fan} in Appendix~\ref{a:eig-proof} for a concise proof.}:
\vspace{-1em}
\begin{eqn}
    \Tp{k}{K} &\Eq \max_{F \in \mathbb{R}^{n \times k}} \Tr[(F^\t K^2 F) (F^\t K F)^\Mo] = \sum_{i=1}^k \lambda_i(K),
    \req{eig-def}
\end{eqn}
with $\lambda_i(K)$ the $i$-th eigenvalue of $K$ in descending order. 
Since \eq{eig-def} is a relaxation of \eq{id-def}, we immediately obtain the upper bound on the approximation quality of any column subset selection that $\Lk_K(\Opt_k) \leq \Tp{k}{K}$.

Having introduced our general problem formulation, we next explore three applications of the low-rank approximation problem introduced above in Sections~\ref{s:app-sensor}-\ref{s:app-laplacian}.
In each, by assuming a kernel matrix $K$ of a particular form, we reduce the approximation task in each application to the same maximization of $\Lk_K(\I)$ \eq{objective-def} over index subsets $\I$. 

\subsection{Sensor placement and experimental design \label{s:app-sensor}}
Consider a kernel or covariance matrix $K$ induced by a Gaussian process on a set of points, e.g., possible sensor placements.
Assuming that a function is drawn from the Gaussian process prior, the conditional covariance of function outputs on all points after observation of outputs at points $\I$ \cite{Williams2006-yk} is the Schur complement $K - K_{:, \I} (K_{\I, \I})^\Mo K_{\I, :}$.
Therefore maximization of $\Lk_K$, our chosen objective, is equivalent to minimization of the average variance over all outputs conditioned on the observations of output on $\I$.

This correspondence suggests applicability to experimental design and sensor placement.
Given a function which is expensive to evaluate (for example, a computational or physical experiment), the goal of this application is to learn as much information about the function with as few evaluations as possible.
Although exact determination of the optimal evaluation subset is infeasible in such settings, the cost of any appproximate approach with reasonable computational scaling may be relatively insignificant compared to the cost of the evaluations themselves. 

An emphasis in the recent literature on large-scale Gaussian processes has been to develop \emph{sparse kernel representations} of the covariance kernel $K$ (e.g. \cite{Liu2020-ed,Noack2023-ar,Noack2021-nb}). 
Therefore, it is important to develop column selection methods that can seamlessly exploit kernel sparsity. Indeed, sparsity of the kernel matrix can be exploited naturally by matrix-free approaches since the cost of matrix-vector multiplication is $O(\mathrm{nnz}(K))$.

\subsection{CUR decomposition via row and column selection \label{s:cur-objective}}

Next we review how the problem of computing the CUR decomposition of a general $m \times n$ matrix $A$ \cite{Goreinov1997-yi} can be reduced to column selection for SPSD matrices.

First consider the problem of computing a CX or interpolative decomposition, i.e., $A \approx C X$, where $C$ consists of columns of $A$ and $X$ is arbitrary. Given such a matrix $C = A_{:, \J}$, the choice of $X$ that minimizes the Frobenius norm error $\Vert A - CX \Vert_{F}$ is $X = C^+ A$, where where $^+$ indicates the pseudoinverse. Then it is natural to attempt to choose $\J$ to minimize the combinatorial objective:
\begin{eqn}
    \epsilon_A(\J) \Eq \min_{B \in \RR^{k,n}} \norm{A_{:, \J} B - A}_F^2 = \Vert A_{:,\J} A_{:,\J}^+ A - A \Vert_F^2
    \req{cur_err}
\end{eqn}
In fact, minimization of this objective is equivalent to the optimization problem \eq{optimal-def} introduced above for the choice $K = A^\t A$ (see proof in Appendix~\ref{a:objective-formulation}):
\begin{restatable}{lemma}{TraceCUR}
$\epsilon_A(\J) = \mathcal{E}_K (\J)$, where $K = A^\t A$.
\label{lem:cur-reduction}
\end{restatable}

Thus, we can select rows $R = A_{\I, :}$ by computing a CX factorization for $A^\t$, which reduces to column selection for $K = AA^\t$. Given columns and rows $C$ and $R$, the choice of $U$ minimizing $\Vert A - CUR\Vert_F$ is given by $U = C^+ A R^+$.
Although rows and columns are selected independently from one another, in fact the Frobenius norm error of our CUR decomposition is controlled by our nuclear norm objectives for the column and row selections (see proof in Appendix~\ref{a:objective-formulation}):
\begin{restatable}{lemma}{IndependentCUR}
Given an $m \times n$ matrix $A$, let $C = A_{:,\J}$ and $R = A_{\I, :}$, and choose $U = C^+ A R^+$. Then
$\Vert A - CUR \Vert_F \leq \mathcal{E}_{A A^\t}^{1/2} (\I) + \mathcal{E}_{A^\t A}^{1/2} (\J)$.
\end{restatable}

In summary, we approach CUR decomposition with the following procedure:
\begin{enumerate}
    \item Select columns $\I$ of the SPSD matrix $A A^\t$ to maximize $\Lk_{A A^\t}(\I)$.
    \item Select columns $\J$ of the SPSD matrix $A^\t A$ to maximize $\Lk_{A^\t A}(\I)$.
    \item Return $C=A_{:, \J}$, $R=A_{\I, :}$, and $U = C^+ A R^+$.
\end{enumerate}
The key difficulty arises from steps 1 and 2, i.e., column selection for suitable SPSD matrices $K = A^\t A, AA^\t$. Note that it may be desirable to avoid ever forming $K$ as such a matrix-matrix product, and in fact our matrix-free approach (Section~\ref{s:matrix-free}) will achieve this.

\subsection{Column subset selection for graph Laplacians \label{s:app-laplacian}}

Finally, we consider the rank reduction of operators based on graph Laplacians.
Such operators are of core interest in graph theory and Markov chain analysis with applications to networks, flow problems, and chemical systems \cite{Merris1994-aw,Aldous1995-bg,Chung2000-gd,Van_Mieghem2014-dl}.
A significant body of research (e.g., \cite{Spielman2004-ex,Cohen2018-ji,Kyng2016-qv,Cohen2016-an}) has demonstrated that sparse Laplacian linear systems can be solved in nearly linear time, a feature that we will exploit in our matrix-free algorithms (Section~\ref{s:matrix-free}).
We shall restrict ourselves to connected and symmetric graph Laplacians, defined as follows:

\begin{definition}[Graph Laplacians and rescaled Laplacians] 
A matrix $\bar L$ is a \emph{graph Laplacian} if it can be written in terms of a symmetric, entrywise nonnegative matrix $(w_{i,j})$ as 
\begin{eqn}
    \bar L_{i,j} = \begin{cases}
        \sum_{k \neq i} w_{i,k} & i = j \\
        -w_{i,j} & i \neq j
    \end{cases}
\end{eqn}    
Assuming that the nonzero pattern of $(w_{i,j})$ defines a connected graph, $\bar L$ is SPSD with one-dimensional null space spanned by $\One$.
More generally, we say that a matrix $L$ is a \emph{rescaled Laplacian} if it may be written $L = \Diag(\h^\Mo) \bar L \Diag(\h^\Mo)$ for some graph Laplacian $\bar L$ and an entrywise positive vector $\h$ with $\norm{h}_2 = 1$. 
Such a matrix $L$ is SPSD with a one-dimensional null space spanned by $\h$ (such that $L \h = \Zero$).
\label{def:laplacian}
\end{definition}

In this work, we consider $L$ and $\h$ to both be pre-specified inputs. 
For a graph Laplacian, $\h \propto \One$ is trivially known ahead of time.
On the other hand, we are particularly interested in applications of our computations to reversible  Markov chains \cite{Aldous1995-bg}.
In this context, given pre-specified rate matrix $R$ and stationary distribution $\pi$, we may consider that $\h = \sqrt{\pi}$ and $L = \Minus \Diag(\h) R \Diag(\h^{\Mo})$. Note that indeed $\norm{\h}_2^2 = \One^\t \pi = 1$.
This interpretation is described in depth in Appendix~\ref{s:markov}. 
More generally, it is possible to consider strictly positive definite matrices $L$ corresponding to killed Markov chains (e.g. \cite{Sharpe1988-pi}). For concreteness we will not discuss this case, which is in fact simpler due to the lack of a null space.

Regardless of interpretation, in this work we shall specifically consider the following minimization:
\begin{eqn}
	\min_{\Abs{\I}=k} \Tr[(L_{\Comp \I, \Comp \I})^\Mo]
    \req{laplacian-min}
\end{eqn}
where $\Comp \I$ is the set complement of $\I$, i.e. $\Comp \I \Eq \{1, \dots, n\} \setminus \I$.
While we might expect this minimization to be equivalent (via Schur complements) to a maximization of $\Lk_{K}$ \eq{objective-def} with $K = L^\Mo$, this inverse does not exist due to the stationary mode $\h$.
Fortunately, we may still re-pose \eq{laplacian-min} by considering it in the limit $\lim_{\Shift \rightarrow 0} (L+ \Shift h h^\t)^\Mo$.

Specifically, to do so we will take $K \Eq L^+$, the Moore-Penrose pseudoinverse of $L$, known in the Markov chain literature as its discrete Green's function \cite{Chung2000-gd} or fundamental matrix \cite{Aldous1995-bg}, modulo diagonal scalings.
With this definition, $K$ is a finite SPSD matrix satisfying $K \h = \Zero$ and is efficiently applicable via existing Laplacian solvers.
While
$ 
	\Tr[(L_{\Comp \I, \Comp \I})^\Mo] \neq \Tr[K] - \Tr[(K^2)_{\I, \I} (K_{\I,\I})^\Mo]
$, 
in general, due to the pathology of the graph Laplacian's stationary state, 
we can still reformulate \eq{laplacian-min} as an analogous maximization problem via the following lemma (proof in Appendix~\ref{a:laplacian-schur}):

\begin{restatable}[Complementary trace formulation for inverse Laplacians]{lemma}{LaplacianComplement}
Let $L$ be a rescaled Laplacian (Definition~\ref{def:laplacian}), let $h$ satisfy $L h = \Zero$ and $\norm{h}_2 = 1$, let $K = L^+$, and let $\I$ be a non-empty subset of $\lrb{1, \dots, n}$. Then:
\begin{eqn}
    \Tr[(L_{\Comp \I, \Comp \I})^\Mo] = \Tr[K] - \Tr[(K^2)_{\I, \I} (K_{\I,\I})^\Mo] + \frac{1+h_{\I}^\t (K_{\I,\I})^\Mo (K^2)_{\I,\I} (K_{\I,\I})^\Mo h_{\I}}{h_{\I}^\t (K_{\I,\I})^\Mo h_{\I}}
    \req{laplacian-complement}
\end{eqn}
\end{restatable}
In essence, our goal is to form a low-rank approximation to the inverse (roughly speaking) of $L$, rather than $L$ itself. Our resulting column selection converges in the $\Shift \rightarrow 0$ limit, even if the inverse $(L + \Shift h h^\t)^\Mo$ does not.
In short, consideration of the stationary state in our formulation adds a single term to our objective maximization. Moreover, this extra term is relatively easy to compute exactly, as it reflects only a rank-one component.
For a graph Laplacian $L$ and its generalized inverse $K$, we therefore consider column selection to maximize the following objective: 
\begin{eqn}
    \Lap_L^h(\I) \Eq \Tr[(K^2)_{\I, \I} (K_{\I,\I})^\Mo] - \frac{1+h_{\I}^\t (K_{\I,\I})^\Mo (K^2)_{\I,\I} (K_{\I,\I})^\Mo h_{\I}}{h_{\I}^\t (K_{\I,\I})^\Mo h_{\I}}
    \req{laplacian-objective}
\end{eqn}
The advantages of formulating our optimization in terms of $K$ include the facts (1) that it is a finite operator efficiently applicable via modern graph Laplacian solvers and (2) that the subsets $\I$ are much smaller than $\Comp \I$ when $k \Eq \Abs{\I} \ll n$.

\section{Deterministic algorithms \label{s:exact-alg}}

In this section, we describe a method for maximizing the nuclear norm objective $\Lk_K(\I)$ posed in \eq{objective-def}.
The algorithm is based on greedy maximization of this objective.
We develop an efficient Cholesky-based implementation, showing that its computational complexity is tractable for many problems of interest.
Algorithm~\ref{alg:dumbest} conveys our deterministic approach at a high level. Ultimately our practical implementation will efficiently produce the same mathematical output.
\begin{algorithm}[h] \begin{algorithmic}
    \setstretch{1.2}
    \caption{Conceptual approach to nuclear maximization}
	\Procedure{ConceptualNuclearMaximization}{$K,k$}
    \State $\I \gets \emptyset$
	\For{$t = 1, \dots, k$}
		\State $\I \gets \I \cup \argmax_{\j \notin \I} \Tr[(K^2)_{\Ij, \Ij} (K_{\Ij,\Ij})^\Mo]$
	\EndFor
	\EndProcedure
    \label{alg:dumbest}
\end{algorithmic} \end{algorithm}
In each iteration of Algorithm~\ref{alg:dumbest}, the gain of adding an additional column is calculated, and a column corresponding to the maximum gain is chosen.
Full enumeration of $k$-subsets is avoided in favor of a greedy approach scaling only linearly with $k$.
While approximate, this selection strategy is adaptive, with column scores dependent on past column choices.
We demonstrate in Sections~\ref{s:error-analysis} and \ref{s:benchmarks} (respectively) that this approach produces approximations which are of theoretically guaranteed and empirically validated accuracy.

\subsection{Reduction to single column selection \label{s:cd-approach}}

In this section, we begin the work of turning Algorithm~\ref{alg:dumbest} into a computationally viable method.
We do so first by working out a more efficient formulation for the marginal increase in $\Lk$ for a chosen column at any given iteration.
To introduce our approach, it is convenient to examine the initial iteration.
In this case, the greedy column choice solves the maximization problem:
\begin{eqn}
    \max_\j \Lk_K(\{\j\}) = \max_\j \frac{(K^2)_{\j,\j}}{K_{\j,\j}} = \max_\j \frac{\norm{K_{:,\j}}^2_2}{K_{\j,\j}}.
    \req{first-obj}
\end{eqn}
Thus both the diagonal entries and the column norms of $K$ are involved in the selection procedure.

Next, consider the augmentation of a column subset $\I$ with one additional index $\j$.
Then from standard use of block matrix identities, the objective $\Lk_K$ cleanly decomposes:
\begin{eqn}
    \Lk_K(\Ij) 
    = \Lk_K(\I) + \frac{(\tK^2(\I))_{\j,\j}}{(\tK(\I))_{\j,\j}}.
    \req{one-off}
\end{eqn}
See Appendix~\ref{a:trace-split-column} for an elementary proof. As a result, we see that the greedy augmentation may be performed by the following maximization:
\begin{eqn}
    \max_{\j} \frac{(\tK^2(\I))_{\j,\j}}{(\tK(\I))_{\j,\j}} 
    = \max_{\j} \left[ \frac{\Diag(\tK^2(\I))}{\Diag(\tK(\I))} \right]_{\j}
\end{eqn}
where vector division is defined elementwise.%
\footnote{Comparing \eq{one-off} to \eq{logdet} makes it clear that greedy nuclear selection maximizes the ratio $(\tK)^2_{\j,\j} / \tK_{\j,\j}$, while greedy diagonal selection maximizes only the denominator term $\tK_{\j,\j}$.}
In other words, the greedy augmentation may be easily performed as long as this ratio of diagonal factors can be estimated with high accuracy and low cost.
This motivates Algorithm~\ref{alg:simple-score-version}, which returns the same column selection as Algorithm~\ref{alg:dumbest}.

\begin{algorithm}[H] \begin{algorithmic}
    \setstretch{1.2}
    \caption{Nuclear maximization via modification of $K$ \label{alg:naive-greedy}}
	\Procedure{IterativeNuclearMaximization}{$K,k$}
    \State $\I \gets \emptyset$
	\For{$t = 1, \dots, k$}
        \State $\tK \gets K - K_{:,\I} (K_{\I,\I})^\Mo K_{\I,:}$
	\State $\I \gets \I \cup \argmax_{\j \notin \I} \dfrac{(\tK^2)_{\j,\j}}{\tK_{\j,\j}}$
	\EndFor
	\EndProcedure
    \label{alg:simple-score-version}
\end{algorithmic} \end{algorithm}

\subsection{Efficient exact implementation \label{s:exact-implementation}}

While much more efficient than Algorithm~\ref{alg:dumbest}, Algorithm~\ref{alg:simple-score-version} involves far more computation than is actually necessary.
We thus present an efficient Cholesky-derived implementation in Algorithm~\ref{alg:exact-cholesky} (Appendix~\ref{s:alg-exact-cholesky}).%
\footnote{As mentioned before, \cite{Farahat2011-zm} presents a similar algorithm. We believe that the complexities for dense $K$ are the same, although that paper does not address potential sparsity or structure in $K$.}
In particular, Algorithm~\ref{alg:exact-cholesky} maintains exact representations of $\Diag(\tK^2)$ and denominator $\Diag(\tK)$ terms which are iteratively updated. For general (possibly sparse) $K$ supporting $O(1)$ element access, the time complexity is $O(n k^2 + k \nnz(K))$.
This conclusion should be compared with the $O(nk^2)$ complexity of standard pivoted Cholesky, which may be used to perform either diagonal maximization or diagonal sampling. The difference in complexity arises from the matvecs needed to update the numerator $\Diag(\tK^2)$ term.

In general, Algorithm~\ref{alg:exact-cholesky} is efficient when $K$ admits fast computation of $\Diag(K)$, $\Diag(K^2)$, and matvecs.
In addition to the case of sparse $K$, we might also consider $K$ to be of Kronecker product form: $K = K^{(1)} \otimes \cdots \otimes K^{(d)}$. 
Suppose in this case that each $K^{(i)}$ is dense and of equal size $(m,m)$ yielding $K$ of size $n = m^d$.
Then, $\Diag(K) = \mvec(\Diag(K^{(1)}) \otimes \cdots \otimes \Diag(K^{(d)}))$ is computable in $O(n d)$ time,
$\Diag(K^2) = \mvec(\Diag([K^{(1)}]^2) \otimes \cdots \otimes \Diag([K^{(d)}]^2))$ is computable in $O(d m^2 + n d)$ time,
and matvecs are 
computable in $O(d m n)$ time.

\subsection{Reduction of graph Laplacians \label{s:deterministic-laplacian}}

We may also extend our deterministic treatment to one targeting the maximization of $\Lap_L^h$ \eq{laplacian-objective} for graph Laplacian $L$ and stationary eigenvector $\h$.
(Recall that $L h = \Zero$, $\norm{h}_2=1$, and we define $K = L^+$, such that $K h = \Zero$.)
We shall see that this objective is computable in similar fashion, so long as the stationary state is properly treated.
Namely, recall from \eq{laplacian-objective} that we seek to maximize:
\begin{eqn}
    \Lap_L^h(\I) = \Tr[(K^2)_{\I, \I} (K_{\I,\I})^\Mo] - \frac{1+h_{\I}^\t (K_{\I,\I})^\Mo (K^2)_{\I,\I} (K_{\I,\I})^\Mo h_{\I}}{h_{\I}^\t (K_{\I,\I})^\Mo h_{\I}}
\end{eqn}
We can compute this objective for a singleton $\I = \lrb{\j}$:
\begin{eqn}
	\Lap^h_L(\lrb{\j}) = \Minus \frac{K_{\j, \j}}{h_\j^2}
    \req{first-column-laplacian}
\end{eqn}
Next, let $\I$ be a \emph{non-empty} set of indices from $\lrb{1, \dots, n}$ and define
\begin{eqn}
	\hat{K}(\I) &\Eq \tilde{K}(\I) + \frac{ (h - K_{:,\I} (K_{\I, \I})^\Mo h_{\I})(h - K_{:,\I} (K_{\I, \I})^\Mo h_{\I})^\t }{h_{\I}^\t (K_{\I, \I})^\Mo h_{\I}},
    \req{khat}
\end{eqn}
Then, as we verify in Appendix~\ref{a:trace-split-column}, the gain in the nuclear objective from adding column $\j$ into non-empty set $\I$ can be computed as:
\begin{eqn}
	\Lap_L^h(\Ij) - \Lap_L^h(\I) = \frac{(\hat{K}^2(\I))_{\j,\j}} {(\hat{K}(\I))_{\j,\j}}
    \req{latter-columns-laplacian}
\end{eqn}

Identities \eq{first-column-laplacian} and \eq{latter-columns-laplacian} permit a similar algorithm to the one considered above, in which the selection is based on entrywise ratios of suitable matrix diagonals.
In Appendix~\ref{s:alg-exact-laplacian}, we formally present Algorithm~\ref{alg:exact-laplacian}, which extends the Cholesky-based approach in Algorithm~\ref{alg:exact-cholesky} to this setting while maintaining the same computational complexity.
This algorithm is not intended for use with large graph Laplacians, as it assumes precomputation of $K$, which might incur $O(n^3)$ cost in general.
For a more scalable approach, we design a matrix-free algorithm in the next section.

\section{Matrix-free algorithms \label{s:matrix-free}}

Existing methods for \Nystrom{} approximation generally assume that $K$ is stored in memory or supports efficient queries to arbitrary entries.
In this work, however, we highlight and address scenarios in which $K$ does not support efficient element access.
Specifically, we consider:
\begin{enumerate}
    \item CUR decomposition of a rectangular, possibly sparse matrix $A$, in which we must consider $K = A^\t A$ or $K=AA^\t$ (cf. Section~\ref{s:cur-objective}). Whether $A$ is sparse or dense, it is desirable to avoid explicit computation of these matrix-matrix products.
    \item Inverse Laplacian reduction, in which we take $K = L^+$ for a graph Laplacian $L$ (cf. Section~\ref{s:app-laplacian}).
    Often $L$ is large, sparse, and poorly conditioned.
    While $K$ is generally infeasible to compute, matvecs by $K$ may be achieved using modern Laplacian solvers.
\end{enumerate}
Although we discuss only these scenarios in depth, we expect that our matrix-free formalism could be extensible to a variety of considered operators $K$ encountered in various other contexts.
The central requirements of our method are (1) that $K$ supports efficient matvecs and (2) that $K$ admits a factorization $K = C C^\t$ such that matvecs by $C$ are efficient. Note that our algorithms do not require exact calculation or storage of either $K$ or $C$.

To achieve this, in particular we use a randomized estimator for the diagonal of $\tilde{K}$. We emphasize that the use of this estimator is not tied exclusively to our selection method based on nuclear maximization. In fact, it may be used just as well to perform column selection based on diagonal maximization or sampling.
On the other hand, we observe that in the matrix-free contexts that we consider, the cost of our approach differs only minorly from these alternative strategies. In Section~\ref{s:benchmarks}, we compare the performance of each approach on a diverse set of examples.

\subsection{Randomized diagonal estimation}
\label{s:randdiag}
In this section, we develop a matrix-free approach for greedy optimization of $\Lk_K$, where we assume only matvec access to operators $K$ and $C$ satisfying $K = C C^\t$.
We will modify this construction for the graph Laplacian case in Section~\ref{s:sqrt-alg}.
Our approach is motivated by the following identity expressing the diagonals we need to compute as expectations: 
\begin{eqn}
    \Diag(K^2) = \Ex_{x \sim \mathcal{N}(\Zero, \Id)} [(K x)^\otwo], \quad 
    \Diag(K) = \Ex_{x \sim \mathcal{N}(\Zero, \Id)} [(C x)^\otwo].
    \req{random-asymptote}
\end{eqn}
Here $^\otwo$ denotes elementwise squaring.

We will discuss how to approximate \eq{random-asymptote} with empirical averages, but first we must note that while \eq{random-asymptote} can be used to compute the scores that we need for the selection of the first column, we must also handle computation of $\Diag(\tK^2(\I))$ and $\Diag(\tK(\I))$ during later iterations.
Our estimator can be extended accordingly by considering $\tK(\I) = K-K_{:,\I}(K_{\I,\I})^{-1} K_{\I,:}$ defined in \eq{Ktilde0} as a projection of $K$ onto the subspace defined by indices $\Comp{\I}$:
\begin{eqn}
    \tK(\I) = (\Id - K_{:,\I} (K_{\I,\I})^\Mo \Id_{\I,:}) K (\Id - \Id_{:,\I} (K_{\I,\I})^\Mo K_{\I, :}).
\req{Ktilde}
\end{eqn}
It follows that $\tK(\I) = \tilde{C}(\I) \, \tilde{C}(\I)^\t$, where we define
\begin{eqn}
\tilde{C}(\I) := C - K_{:,\I} (K_{\I,\I})^\Mo C_{\I,:}.
\req{Ctilde}
\end{eqn}
Note that both $\tK = \tK (\I)$ and $\tilde{C} = \tilde{C}(\I)$ allow efficient matvecs via \eq{Ktilde} and \eq{Ctilde}, respectively, given that we can perform matvecs by $K$ and $C$, and we have, similar to \eq{random-asymptote}, that
\begin{eqn}
    \Diag(\tK^2) = \Ex_{x \sim \mathcal{N}(\Zero_n, \Id_n)} [(\tK x)^\otwo], \quad 
    \Diag(\tK) = \Ex_{x \sim \mathcal{N}(\Zero_n, \Id_n)} [(\tilde{C} x)^\otwo].
    \req{random-asymptote2}
\end{eqn}
We can use empirical averages to approximate \eq{random-asymptote2}, yielding
\begin{eqn}
    \frac{\Diag(\tK^2)}{ \Diag(\tK)}  \approx \frac{ (\tK Z^{(1)})^\otwo \One} 
        { (\tilde{C} Z^{(2)})^\otwo \One }
        \req{random-finite},
\end{eqn}
where $Z^{(1)}$ and $Z^{(2)}$ are matrices of size $(n,z)$ with elements independently drawn from $\mathcal{N}(0, 1)$. Alternatively, we can write this approximation entrywise as
\begin{eqn}
    \frac{ [\tK^2]_{\j,\j} } { \tK_{\j,\j}}  \approx \frac{ \norm{ [\tK Z^{(1)}]_{\j,:} }_2^2 } { \norm{ [\tilde{C} Z^{(2)}]_{\j,:} }_2^2 }
        \req{random-finite-elementwise}.
\end{eqn}
As the numerator and denominator terms are positive almost surely, this approach yields positive score estimates.
In Section~\ref{s:concentration}, we will prove that the estimator achieves controlled relative error with high probability for a number of matvecs $z = O(\log n)$.

Algorithm~\ref{alg:randomized-cholesky} (Appendix~\ref{s:alg-randomized-cholesky}) presents our resulting algorithm for matrix-free column selection, in which $z$ denotes the number of of independent Gaussian vectors used for matvecs in each iteration.
With this notation, the time complexity of Algorithm~\ref{alg:randomized-cholesky} is $O(n k^2 z + k z \psi)$, where $\psi$ is the greater cost of performing a matvec by $K$ or $C$. 
As mentioned above, in the next section we will show that we need only take $z = O(\log n)$ to achieve fixed relative error for these diagonals with high probability.
The randomization introduced thus eliminates the explicit computations of $\Diag(K)$ and $\Diag(K^2)$ at the cost of additional $O(k z) \sim O(k \log n)$ matvecs. 

\subsection{Concentration bounds \label{s:concentration}}

The relative error of our randomized score approximation may be rigorously controlled by extending stochastic trace estimation results \cite{Rudelson2013-dc,Meyer2021-pn} to the problem of diagonal estimation.
(Full proofs are deferred to Appendix~\ref{a:concentration-proofs}.)
In particular, the following formulation is similar to \cite{Lindsey2023-fn}:

\begin{restatable}[Stochastic diagonal approximation]{theorem}{DiagonalApprox}
Let $X$ be an $n\times n$ SPSD matrix, and suppose that $X = Y Y^\t$ for some matrix $Y$ which need not be square. 
Let $Z$ be an $n \times z$ matrix whose elements are independently drawn from $\mathcal{N}(0,1)$.
Then for any $\epsilon > 0$, 
\begin{eqn}
    (1-\epsilon) \Diag(X) \leq \tfrac{1}{z} \Diag(Y Z Z^\t Y^\t) \leq (1+\epsilon) \Diag(X)
\end{eqn}
with probability $1-\delta$ if $z \geq c \log(n/\delta) / \epsilon^2$, where $c$ is a universal constant and $\delta \in (0,1/2]$. Note that $\Diag(C Z Z^\t C^\t) = (CZ)^\otwo \One$.
\label{th:trace-concentration}
\end{restatable}

We can apply Theorem~\ref{th:trace-concentration} in (1) the case where $X = \tK^2$, $Y = \tK$ and (2) the case where $X = \tK$ and $Y = \tilde{C}$, to get relative error bounds on the numerator and denominator, respectively, of our estimator \eq{random-finite}-\eq{random-finite-elementwise} for the diagonal scores. 
In turn, these relative error bounds yield a relative error bound for the quotient itself, which can be used to control the error of our approximate nuclear maximization relative to exact nuclear maximization:
\begin{restatable}[Stochastic ratio maximization]{corollary}{StochasticRatio}
    Let $\tK = \tK (\I)$ and $\tilde{C} = \tilde{C} (\I)$, and let $Z^{(1)}$ and $Z^{(2)}$ be matrices of size $n \times z$, with all elements drawn independently from $\mathcal{N}(0,1)$.
    Then:
    \begin{eqn}
        \max_{\j} \frac{ \norm{ [\tK Z^{(1)}]_{\j,:} }_2^2 } { \norm{ [\tilde{C} Z^{(1)}]_{\j,:} }_2^2 }
        \geq \Orand^\Mo \max_{\j} \frac{({\tK}^2)_{\j,\j}}{ \tK_{\j,\j}}
    \end{eqn}
    with probability at least $1-4\delta$ if $z \geq c' \log(n/\delta) / \Rand^2$ where $c'$ is a universal constant, $\delta \in (0,1/4]$, and $\Rand \in (0, 1]$.
\end{restatable}
Thus, for fixed relative approximation error and failure probabilities, it suffices to choose $z$ scaling only logarithmically with $n$. 
For our theoretical error bounds, we shall assume that we have guaranteed approximation error $\Rand$ across all $k$ iterations of our selection algorithm.
Algorithmically, we will always choose independent random vectors $Z^{(1)}$ and $Z^{(2)}$ in each iteration.
In this case, taking a union bound over all iterations trivially yields that we may choose $z \geq c' \log(k n / \delta) / \Rand^2$ to achieve approximation error $\Rand$ over all iterations with probability $1-4\delta$.

\subsection{Extension to approximate factorizations \label{s:sqrt-alg}}

In the context of CUR decomposition, we are provided an \emph{exact} symmetric factorization $K = CC^\t$ by construction. To wit, $C = A^\t$ for column selection and $C=A$ for row selection.
However, we are also interested in settings where only \emph{approximate} factorization is possible.
We specifically examine our graph Laplacian application, although the idea should be useful when generally for $K$ that can be viewed as a matrix function of a sparse matrix and suitable preconditioning is possible.

In the graph Laplacian setting, for preconditioning we draw on recent developments in the fast Laplacian solver literature. 
In particular, we shall exploit the concept of approximate sparse Cholesky factorization, i.e., the computation of a sparse triangular matrix $R$ such that $a L \leq R R^\t \leq b L$ for some constants $0 < a \leq b$. 
Given an approximate factorization of a graph Laplacian $\bar L \approx \bar R \bar R^\t$, we may choose $R = \Diag(h^{\Mo}) \bar R$ to yield $L \approx R R^\t$ to obtain the approximate factorization we require.

Most relevantly, \textcite{Kyng2016-qv} invented a particularly simple nearly linear-time algorithm yielding $R$ in $O(\nnz(\bar L) \log^3 n)$ expected time with $\nnz(R) = O(\nnz(\bar L) \log^3 n)$, which is a sparse Cholesky factor for $a=1/2$ and $b=3/2$ with probability $1 - 1/\mathrm{poly}(n)$.
While this algorithm is expensive in practice, practical alternatives with weaker theoretical guarantees have also been recently developed \cite{Gao2023-vh, Chen2020-ya}.
We use the randomized Cholesky (\texttt{rchol}) algorithm and software of \cite{Chen2020-ya} in our examples, diagonally rescaling the output factor as described above.

To implement our approach, suppose we have a sparse Cholesky factor $R$ for $L$ and some $a,b$. Note that necessarily $R^\t$ and $L$ share the same null space, while $R$ and $L$ share the same column space.
Also, as $R$ is triangular, matvecs by $R^+$ or $R^{+\t}$ are simple to perform in $O(\nnz(R))$ time.

Given our convention that $K = L^+$, then:
\begin{eqn}
K = R^{+} (R^{+} L R^{+\t} )^+ R^{+ \t} = C C^\t,
\end{eqn}
where, with $^{+/2}$ indicating the square root of the pseudoinverse:
\begin{eqn}
     C := R^{+} (R^{+} L R^{+\t} )^{+/2}.
\end{eqn}
Considering our objective \eq{laplacian-objective} for one column therefore yields, via \eq{first-column-laplacian}:
\begin{eqn}
	\Lap^h_L(\lrb{\j}) = \Minus \frac{K_{\j,\j}}{h_\j^2} = \Minus \frac{\Ex_{x \sim \mathcal{N}(\Zero,\Id)}[[C x]_\j^2]}{h_\j^2}.
\end{eqn}
For subsequent columns ($\I \neq \emptyset$), it is helpful to 
define intermediate quantities
\begin{eqn}
    \tau(\I) \Eq h_{\I}^\t (K_{\I,\I})^\Mo h_{\I} \qquad
    y(\I) \Eq h - K_{:,\I} (K_{\I,\I})^\Mo h_{\I}
    \req{laplacian-complicated-defs}
\end{eqn}
such that, omitting $(\I)$ where implied, using $\tilde K$ from \eq{Ktilde0} and $\tilde C$ from \eq{Ctilde}, and applying \eq{latter-columns-laplacian}:
\begin{eqn}
	\Lap_L^h(\Ij) - \Lap_L^h(\I) 
    = \frac{\Ex_{x \sim \mathcal{N}(\Zero,\Id)}[[(\tilde K + \tau^\Mo y (y-h)^\t) x]_\j^2] + y_{\j}^2 / \tau^2}
    {\Ex_{x \sim \mathcal{N}(\Zero,\Id)}[[\tilde C x]_\j^2] + y_{\j}^2 / \tau}.
    \req{laplacian-complicated-estimator}
\end{eqn}
(See Appendix~\ref{a:laplacian-factor} for details.)
As the right-hand terms of the numerator and denominator are non-negative and exactly computable, controlled relative error is guaranteed with high probability using the same diagonal estimators from Section~\ref{s:concentration}.
Thus nuclear score estimation can be performed at all iterations given the ability to perform matvecs by $K$, $R^+$, and the pseudoinverse square root of $R^+ L R^{+\t}$. Matvecs by $K$ can be performed efficiently using the preconditioned conjugate gradient method, in which $R$ determines the preconditioner.

To approximately compute matvecs by $(R^{+} L R^{+\t})^{+/2}$, we will leverage our capacity to perform matvecs by $R^{+} L R^{+\t}$. We assume that the condition number $\kappa \Eq b/a$ of $R^{+} L R^{+\t}$ on its image is relatively small. Specifically, our complexity analysis assumes that $\kappa = O(1)$, independent of dimension. 
Given this assumption, a variety of approximate techniques for matrix function computation (e.g., \cite{Hale2008-ae}) could suffice.
We choose to invoke the Chebyshev approximant of $f(x)= ( [a+b+(b-a) x]/2)^{\Minus 1/2}$ on $x \in [-1,1]$ \cite{Trefethen2020-cm}.
An analysis of this approximation is given in Lemma~\ref{lem:chebyshev} (Appendix~\ref{s:chebyshev-analysis}), showing that the number $n_\mathrm{cheb}$ of matvecs by $R^{+} L R^{+\t}$ needed to achieve relative approximation error $\ve$ in the infinity norm is asympotically bounded as $n_\mathrm{cheb} \lesssim \frac{1}{2} \sqrt{\kappa } \log (\kappa \sqrt{n} / \ve)$. 
We provide our full matrix-free algorithm for graph Laplacian rank reduction in Algorithm~\ref{alg:randomized-laplacian} (Appendix~\ref{s:alg-randomized-laplacian}).

\section{Error analysis \label{s:error-analysis}}

\begin{figure}
    \centering
    \includegraphics[width=5.5in]{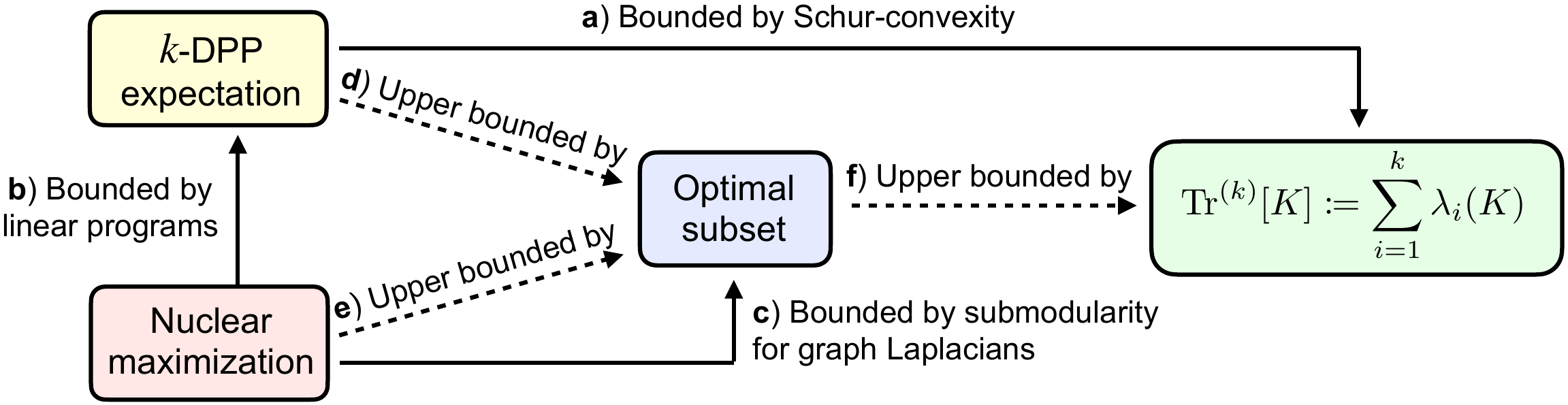}
    \caption{Overview of theoretical error bounds.
    Dashed arrows refer to inequalities in the form of upper (best-case) bounds.
    Solid arrows refer to inequalities in the form of lower (worst-case) bounds.
    Inequalities (d), (e), and (f) are elementary and indicate the best performance one could achieve.
    Inequality (a) is from the theory of $k$-DPP sampling and elementary symmetric polynomials, while
    inequalities (b) and (c) are new contributions.
    \label{fig:overview}
    }
\end{figure}

We now turn to the task of establishing performance bounds of our given deterministic and matrix-free algorithms.
We mostly defer proofs and extended results to Appendix~\ref{a:dpp-proofs} (for proofs about DPP sampling) and Appendix \ref{a:error-bounds} (for all other proofs and results).
Our bounds are summarized visually in Figure~\ref{fig:overview}.

First, in Section~\ref{s:submodularity}, we establish a generic set of linear programming bounds (Theorem~\ref{th:submodular}) generalized from the theory of submodular functions.
In particular, these bounds are able to handle approximate greedy selection even though our objectives are not submodular, and they can encompass simple stopping conditions (Corollaries~\ref{cor:accumulated-constraint} and \ref{cor:initial-constraint}).

Next, in Section~\ref{s:kernel-error}, we demonstrate that for an arbitrary SPSD matrix $K$, this set of linear programs yields immediate bounds on the discrepancy in $\Lk_K$ between the outputs of nuclear maximization (Algorithms~\ref{alg:exact-cholesky} and \ref{alg:randomized-cholesky}) and the expectation value of DPP sampling (Theorem~\ref{th:submodular-dpp}).
This in turn motivates us to prove certain inequalities to bound the error of DPP sampling itself (Appendix~\ref{a:dpp-proofs}), yielding novel proofs of possible independent interest (e.g., Theorem~\ref{a:dpp-concave}).
As a result of our efforts, we are able to pose multiplicative (Theorem~\ref{th:greedy-general-bound}) and additive (Corollary~\ref{cor:additive-re-bound}) error bounds for our algorithms as applied to nearly low-rank matrices.

Finally, in Section~\ref{s:laplacian-analysis}, we extend our results to rank reduction of a rescaled Laplacian $L$.
Here, we are able to use the random walk interpretation of graph Laplacians to establish the submodularity of our objective function $\Lap_L^h$ on non-empty sets (Theorem~\ref{th:laplacian-is-submodular}).
Next, by posing a linear program modified to suitably handle the graph Laplacian's stationary state,
we are able to control the discrepancy between the outputs of nuclear maximization (Algorithms~\ref{alg:exact-laplacian} and \ref{alg:randomized-laplacian}) and the exactly optimal set of $k$ columns (Theorem~\ref{th:laplacian-trace-bound});

We note explicitly that our analysis neglects floating point error. Furthermore, for our iterative approach for graph Laplacians, we neglect the errors of Chebyshev and conjugate gradient computations, assuming that these techniques can be converged within numerical precision without too much expense.
These assumptions appear well-founded judging from the empirical results in Section~\ref{s:benchmarks}.

To bound the quality of columns selected by nuclear maximization, we shall repeatedly reference the following definition, which encompasses the outputs of our deterministic and matrix-free algorithms with pre-specified relative error tolerance $\Rand$.
Taking $\Rand=0$ accommodates the deterministic methods (Algorithms~\ref{alg:exact-cholesky} and \ref{alg:exact-laplacian}), while more generally we must take $\Rand>0$ for our matrix-free methods (Algorithms~\ref{alg:randomized-cholesky} and \ref{alg:randomized-laplacian}, subject to the failure probability described in Section~\ref{s:concentration}).
It is a strength of our proof technique that these cases may be handled straightforwardly in a unified manner.

\begin{definition}[Greedy, approximately greedy, and optimal subsets]
    Given objective function $\subm$ and relative error $\Rand \geq 0$, let $\Ge_{t} \subseteq \lrb{1, \dots, n}$ be an increasing sequence of sets such that $\Ge_{0}\Eq\emptyset$ and, for all  $t=1,\dots,k$:
    \begin{eqn}
        \Abs{\Ge_{t}} &=t \\
        \Ge_{t-1} &\subset \Ge_{t} \\
        \subm(\Ge_{t}) - \subm(\Ge_{t-1}) & \geq \Orand^\Mo \max_{\j} \left\{ \subm(\Ge_{t-1} \cup \lrb{\j}) - \subm(\Ge_{t-1}) \right\}
        \req{approx-greedy-def}
    \end{eqn}
    In particular, $\Gg_k$ is then a subset of size $k$ chosen to greedily maximize $\subm$ with no approximation error.
    Additionally, as shorthand, let $\G_k \Eq \subm(\Ge_k) - \subm(\Ge_{k-1})$ be the gain in objective $\subm$ upon addition of the $k$-th element, greedily selected up to maximum relative error $\Rand$.
    Finally, let $\Opt_k$ be a $k$-subset which exactly maximizes the objective function, i.e. $\Opt_k = \argmax_{\Abs{\I}=k} \subm(\I)$.
    \label{def:greedy}
\end{definition}

\subsection{Linear programming relations \label{s:submodularity}}

In order to bound the performance of our algorithms, we will first derive a set of linear programming relations naturally applicable to greedy selection procedures.
Our results draw on the theory of submodular functions, the basics of which have been laid out by others, cf. \cite{Nemhauser1978-ci}.
Submodularity (Definition~\ref{def:submodular}) is often described intuitively as a notion of ``diminishing returns''; if $A \subseteq B$, the gain in $\subm$ achieved by adding a given element to set $B$ is bounded above by the value of adding the same element to set $A$.
It has most famously been used to prove the near-optimality of greedy algorithms in common cases involving discrete set selection problems which might be NP-hard to exactly solve.
Once it is proved that a set function is submodular, its maximization via greedy augmentation can be bounded relative to the optimal result.

In this section (and Appendix~\ref{a:lp-relations}), we will prove some moderate extensions of the linear program considered in existing theory while giving a self-contained matrix-based formulation of the proofs. 
Because we will apply these programs outside of the scope of submodular functions, we shall keep our presentation in this section general:


\begin{restatable}[Generalized linear programming bound]{theorem}{GeneralSubmodularityConstraint}
Consider the set of inequalities 
\begin{eqn}
    \Lhs &\leq \sum_{i=1}^{t-1} \Y_i + f_t \Y_t, \qquad t = 1, \dots, k
    \req{inequalities}
\end{eqn}
for fixed parameters $\Lhs \geq 0$ and $f_t > 1$ for all $t \in \lrb{1,\dots,k}$.
Then:
\begin{eqn}
    1 - \frac{ \sum_{i=1}^k y_i }{\Lhs} < \exp \Big( \Minus \sum_{i=1}^k f_i^\Mo \Big)
    \req{submodular-general-bound}
\end{eqn}
and in particular, if $\bar f = f_1 = f_2 = \dots = f_k$ then:
\begin{eqn}
    1 - \frac{  \sum_{i=1}^k y_i  }{\Lhs} < e^{\Minus k / \bar f}
    \req{simple-bound}
\end{eqn}
\label{th:submodular}
\end{restatable}


One of the strengths of this bounding technique is that it can encompass additional \emph{ad hoc} constraints.
For instance, the following corollaries consider stopping criteria for the greedy augmentation procedure: (1) stopping if the $(k+1)$-th gain is smaller than $\alpha$ times the cumulative prior gain, and (2) stopping if the $(k+1)$-th gain is smaller than $\beta$ times the initial gain: 

\begin{restatable}[Constraints based on accumulated gain]{corollary}{AccumulatedConstraint}
    Suppose that in addition to the constraints \eq{inequalities}, one adds the constraint that 
    $
        \Y_{k+1} \leq \alpha \sum_{i=1}^{k} \Y_i
    $ 
    for real parameter $0 < \alpha < 1$.
    Then, in addition to obeying \eq{submodular-general-bound}, one has the bound:
    \begin{eqn}
        1 - \frac{\sum_{i=1}^{k} \Y_i}{\Lhs} \leq \frac{\alpha f_{k+1}}{1 + \alpha f_{k+1}} 
        \req{submodular-sum-bound}
    \end{eqn}
    \label{cor:accumulated-constraint}
\end{restatable}
    
\begin{restatable}[Constraints based on initial gain]{corollary}{InitialConstraint}
    Suppose instead that in addition to the constraints \eq{inequalities}, one adds the constraint that
    $
        \Y_{k+1} \leq \beta \Y_1
    $
    for real parameter $0 < \beta < 1$.
    Then, in addition to obeying \eq{submodular-general-bound}, one has the bound:
    \begin{eqn}
        1 - \frac{\sum_{i=1}^{k} \Y_i}{\Lhs} \leq \frac{1}{\prod_{j=2}^k (1-f_j^\Mo)^\Mo+\lrp{\beta  f_{k+1}}^\Mo}
        \req{submodular-first-bound}
    \end{eqn}
    \label{cor:initial-constraint}
\end{restatable}

Appendix~\ref{a:lp-relations} shows how, by applying Theorem~\ref{th:submodular}, Corollary~\ref{cor:accumulated-constraint}, and Corollary~\ref{cor:initial-constraint} to a submodular function $\subm$, one may substitute $\Y=\G$ and $\Lhs = \subm(\Opt)_s$ to upper bound the quantity $1 - \sum_{t=1}^k \G_t / \subm(\Opt_s)$, i.e., the maximum relative performance error between greedy selection of $k$ elements and the exactly optimal set of $s$ elements.
Doing so extends the results of \cite{Nemhauser1978-ci}.
In the following sections, by contrast, we show how the above linear programming results also bound nuclear maximization for kernel approximation (for which the objective is not submodular), and inverse Laplacian reduction (which we prove to be submodular, but only on non-empty sets).

\subsection{Error bounds for kernel approximation \label{s:kernel-error}}

In this section, we turn to our main theoretical question: how well does greedy nuclear maximization (performed exactly or with bounded error) work to optimize the objective \eq{objective-def}?
Most essentially, we will give theoretical lower bounds for the following quantity:
\begin{eqn}
    \Lk_K(\Ge_k) = \Tr[(K^2)_{\I,\I} (K_{\I, \I})^\Mo] \; \Big|_{\I = \Ge_k}
\end{eqn}
where $\Ge_k$ is the set of $k$ columns chosen greedily to maximize objective $\Lk_K$ up to maximum relative error $\Rand$, cf. Definition~\ref{def:greedy}. We mostly defer our proofs in this section to Appendix~\ref{a:kernel-approximation-error}.

While $\Lk_K$ is not a submodular function for general SPSD $K$,  we demonstrate here that  the same  linear programs (Theorem~\ref{th:submodular}) bound the discrepancy between nuclear maximization and DPP sampling.
Essentially, while we cannot substitute $\Lhs  = \Lk_K(\Opt_s)$ in the prerequisite inequalities of Theorem~\ref{th:submodular}, we can provably substitute $\Lhs =\Ex_{\I \sim s\text{-DPP}(K)}[\Lk_K(\I)]$.
To do so, we first define the expectation of our objective under DPP sampling:

\begin{definition}[DPP sampling]
    For a given $n \times n$ SPSD matrix $K$, the corresponding \emph{$s$-DPP} is a probability density over the set $\Psi_s^n \Eq {\{1, \dots, n\} \choose s}$ of all $s$-subsets of $\{1, \ldots ,n \}$, defined by the probability mass function: 
    \begin{eqn}
        p_{s\text{-DPP}(K)}(\I) \propto \det(K_{\I, \I})
    \end{eqn}
    The expectation $\dpp_s(K) \Eq \Ex_{\I \sim s\text{-DPP}(K)}[\Lk_K(\I)]$ is then:
    \begin{eqn}
        \dpp_s(K) &\Eq \frac{\sum_{\I \in \Psi_s^n} \Lk_K(\I) \det(K_{\I,\I})}{\sum_{\I \in \Psi_s^n} \det(K_{\I,\I})}
        \req{dpp-def}
    \end{eqn}
    Because $\dpp_s(K)$ is a (non-negative) weighted average of $\Lk_K(\I)$ over all possible column subsets $\Psi_s^n$, it is immediate that $\dpp_s(K) \leq \Lk_K(\Opt_s)$.
    \label{def:dpp}
\end{definition}

One of the primary attractions of DPP sampling is that $\dpp_s(K)$ can be expressed in closed form in terms of the eigenvalues of $K$ alone.
First, for a vector $x$ of $n$ numbers, let $e_i(x)$ be the elementary symmetric polynomial in $x$ of order $i$, adopting the conventions that $e_0(x) \Eq 1$ and $e_{n+1}(x) \Eq 0$ (see Appendix~\ref{a:dpp-proofs} for details).
Then $\dpp_s(K)$ may be expressed elegantly as follows \cite{Guruswami2012-rr}:

\begin{proposition}[Expectation of $s$-DPP sampling \cite{Guruswami2012-rr}]
Let $K$ be an $n \times n$ SPSD matrix with eigenvalues $\lambda = (\lambda_1 , \ldots ,\lambda_ n)$.
Then for any integer $0 \leq s \leq n$:
\begin{eqn}
    \dpp_s(K) = \dpp_s(\Diag(\lambda)) = e_1(\lambda) - (s+1) \frac{e_{s+1}(\lambda)}{e_{s}(\lambda)}
    \req{dpp-sym}
\end{eqn}
In particular, $e_1(\lambda) = \One^\t \lambda = \Tr[K]$, and so we have that $\dpp_0(K) = 0$ and $\dpp_n(K) = \Tr[K]$.
\label{prop:dpp-sym}
\end{proposition}


For this work, we shall rely on the following salient features of DPP sampling (see Appendix~\ref{a:dpp-proofs}): 
\begin{enumerate}
    \item The $s$-DPP expectation increases with $s$, i.e. $\dpp_{s+1} \geq \dpp_{s}$  for $0 \leq s \leq n$.
    We provide a proof of this inequality in Theorem~\ref{th:dpp-monotonic}, which is a straightforward consequence of Newton's inequalities \cite{Maclaurin1729-mh}.
    \item The $s$-DPP expectation is concave in $s$, in the sense that $\dpp_{s+1} - \dpp_{s} \leq \dpp_{s} - \dpp_{s-1}$ for $1 \leq s \leq n$.
    We provide a novel proof of this property via induction in Theorem~\ref{th:sym-convex}.
    (We are unaware of any prior proof in the literature.)
    \item The $s$-DPP expectation is Schur-convex in its arguments. 
    That is, for any unitary $U$, $\dpp_s(U \Diag(\lambda) U^\t)$ is Schur-convex in $\lambda$.
    We give a self-contained proof of this property in Theorem~\ref{th:schur-convexity}, relying on only Newton's inequalities rather than the more complicated machinery of \cite{Guruswami2012-rr}. 
    \item The $s$-DPP expectation is subadditive, i.e. $\dpp_s(A+B) \leq \dpp_s(A) + \dpp_s(B)$ for symmetric positive semidefinite matrices $A$ and $B$; see Corollary~\ref{cor:dpp-subadditive}, following straightforwardly from \cite{Marcus1957-cl}.
\end{enumerate}
In the next theorem, we show how these properties can be used to replace the optimality metric in conventional submodular analysis with the result derived from DPP sampling. 
We include its proof in the main text to highlight its structure.

\begin{theorem}[Linear program inequalities for nuclear maximization]
    Given SPSD matrix $K$ and objective function $\Lk_K$ \eq{objective-def}, let $k \geq s \geq 1$, $\Ge_t$ and $\G_t$ be defined per Definition~\ref{def:greedy}, and $\dpp_s(K)$ be defined per Definition~\ref{def:dpp}.
Then:
\begin{eqn}
    \sum_{i=1}^{t-1} \G_i(K) + \Orand s \G_t(K) \geq \dpp_s(K) \qquad t = 1, \dots, k
    \req{submodular-dpp}
\end{eqn}
\label{th:submodular-dpp}
\end{theorem}
\begin{proof}
    We will rely on the following facts for general SPSD matrices $A$ and $B$.
    First, the 1-DPP estimate is upper bounded by the exact optimum, which is in turn bounded to the approximately maximizing subset of size 1:
    \begin{eqn}
        \Orand \G_1(A) \geq \dpp_1(A)
        \req{ing-1}
    \end{eqn}
    Second, by Corollary~\ref{cor:dpp-concave}, $\dpp_s$ is concave in $s$,
    from which one clearly has that:
    \begin{eqn}
        s \dpp_1(A) &\geq \dpp_s(A)
        \req{ing-2}
    \end{eqn}
    Third, following from the definition of $\Lk_K$ \eq{objective-def}, $\tK$ \eq{Ktilde0}, and \eq{dpp-sym}, we have that:
    \begin{eqn}
        \Lk_K(\I) = \Tr[K - \tilde{K}(\I)] \geq \dpp_s(K - \tilde{K}(\I))
        \req{ing-3}
    \end{eqn}
    Fourth, by Corollary~\ref{cor:dpp-subadditive}, $\dpp_s$ is subadditive, such that 
    $
    \dpp_s(A) + \dpp_s(B) \geq \dpp_s(A+B)
    $. 
    Thus:
    \begin{eqn}
        \dpp_s(K - \tilde{K}(\I)) + \dpp_s(\tilde{K}(\I)) \geq \dpp_s(K)
        \req{ing-4}
    \end{eqn}
    Finally, by combining these ingredients:
    \begin{eqn}
        \sum_{i=1}^{t-1} \G_i + \Orand s \G_t 
        &\geq \Lk_K(\Ge_{t-1}) + s \dpp_1(\tilde{K}(\Ge_{t-1})) & \text{from \eq{ing-1}} \\
        &\geq \Lk_K(\Ge_{t-1}) + \dpp_s(\tilde{K}(\Ge_{t-1})) & \text{from \eq{ing-2}} \\
        &\geq \dpp_s(K - \tilde{K}(\Ge_{t-1})) + \dpp_s(\tilde{K}(\Ge_{t-1})) & \text{from \eq{ing-3}} \\
        &\geq \dpp_s(K) & \text{from \eq{ing-4}}
    \end{eqn}
    for any $t \geq 1$. Consideration of $t = 1, \dots, k$ yields \eq{submodular-dpp}.
\end{proof}

Having achieved the system of inequalities in Theorem~\ref{th:submodular} with $f_1 = \cdots = f_k = \Orand s$ and $\Lhs = \dpp_s(K)$, we can immediately bound the expected performance of $s$-DPP relative to that of nuclear maximization with $k$ columns:
\begin{restatable}[Discrepancy between nuclear maximization and DPP sampling]{corollary}{GreedyDPP}
    For $\Ge_k$, a subset of $k \geq s$ columns greedily selected with approximation error $\Rand$ (Definition~\ref{def:greedy}):
    \begin{eqn}
        1 - \frac{\Lk_K(\Ge_k)}{\dpp_s(K)}  < e^{\Minus k / (\Orand s)}
    \end{eqn}
    \label{cor:dpp-greedy-bound}
\end{restatable}


While Corollary~\ref{cor:dpp-greedy-bound} lets us bound nuclear maximization to DPP sampling, the spectral guarantees of DPP sampling itself let us go one step further, relating relating $\Lk_K(\Ge_k)$ (the objective achieved by nuclear maximization) to $\Tp{r}{K} = \sum_{i=1}^r \lambda_i(K)$, i.e., the sum of the top $r$ eigenvalues of $K$.
This formulation is most useful in considering matrices $K$ which are nearly rank-$r$, with a small remainder of the trace falling outside of the first $r$ eigenvalues.
In this regime, it is a core goal for column selection methods (e.g., \cite{Derezinski2021-hi,Chen2022-lp}) to achieve an additive $r,\ve$ bound, defined as follows:

\begin{definition}[Additive error control via an $r,\ve$ bound]
    Let $r$ be a positive integer and $\ve > 0$. A column selection algorithm $\I(K)$ satisfies an $r,\ve$ bound if, for any SPSD $K$:
    \begin{eqn}
        \Lk_K(\I(K) ) \leq (1 + \ve) \big( \Tr[K] - \Tp{r}{K} \big)
        \req{re-bound}
    \end{eqn}
\end{definition}

A standard literature result (e.g., \cite{Derezinski2021-hi}) establishes an $r,\ve$ bound for DPP sampling by leveraging the following statement deriving from the Schur-convexity of $\dpp_s$ (Theorem~\ref{th:schur-convexity}):

\begin{lemma}[Low-rank approximation via DPP sampling \cite{Belabbas2009-il,Guruswami2012-rr,Chen2022-lp}]
    For SPSD matrix $K$, the expectation of $s$-DPP sampling $\dpp_s(K)$ \eq{dpp-def} satisfies:
    \begin{eqn}
        \Tr[K] - \dpp_s(K) \leq \lrp{1 + \frac{r}{s-r+1}} \big( \Tr[K] - \Tp{r}{K} \big)
        \req{dpp-rp-bound}
    \end{eqn}
    \label{lem:dpp-error}
\end{lemma}

By rearrangement of Lemma~\ref{lem:dpp-error}, we can say that $s$-DPP sampling satisfies an $r,\ve$ bound, in expectation, when $s  \geq \frac{r}{\ve} + r - 1$.
In fact, optimal selection satisfies an $r,\ve$ bound only with $k \geq r/\ve$ columns \cite{Chen2022-lp,Guruswami2012-rr}, so this result for DPP sampling is quite close to optimal, justifying its status as a relative gold standard for column selection.
By combining Lemma~\ref{lem:dpp-error} with Corollary~\ref{cor:dpp-greedy-bound} (proof in Appendix~\ref{a:kernel-approximation-error}), we can most straightforwardly achieve the following \emph{multiplicative} error bound for nuclear maximization:

\begin{restatable}[Relative approximation error of nuclear maximization]{theorem}{NuclearRelativeError}
    Consider SPSD matrix $K$, $r \geq 1$, and let $\nu \Eq \frac{\Tr[K] - \Tp{r}{K}}{\Tp{r}{K}}$.
    Then in the $k \rightarrow \infty$ limit,
\begin{eqn}
    1 - \frac{\Lk(\Ge_k)}{\Tr^{(r)}[K]} \lesssim \frac{\nu r}{\frac{k}{\log(\nu r)} - r + 1}
\end{eqn}
In the limit $\Rel \rightarrow 0$ from above, to guarantee that $1 - \frac{\Lk(\Ge_k)}{\Tp{r}{K}} \leq \Rel$ it suffices to choose a number of columns $k \geq k^*$ where
\begin{eqn}
    k^* \lesssim \Orand \left(r \nu/\Rel +r-1\right) \log \left(\frac{r-1}{r \nu}+\frac{1}{\Rel }\right).
    \req{general-relative-bound}
\end{eqn}
\label{th:greedy-general-bound}
\end{restatable}

By rearrangement (Appendix~\ref{a:kernel-approximation-error}), we can also establish an (additive) $r,\ve$ bound for nuclear maximization, though one should keep in mind the dependence of the bound upon parameter $\nu$ associated with input $K$: 
\begin{restatable}[Additive $r,\ve$ error bound for nuclear maximization]{corollary}{NuclearAdditiveError}
Consider SPSD matrix $K$, $r \geq 1$, and let $\nu \Eq \frac{\Tr[K] - \Tp{r}{K}}{\Tp{r}{K}}$.
In the $\ve \rightarrow 0$ limit, there exists $k^*$ satisfying
\begin{eqn}
    k^* \lesssim \Orand \left(\frac{r}{\ve }+r-1\right) \left(\log \left(\nu^\Mo\right)+\log \left(\ve^\Mo-r^\Mo+1\right)\right)
    \req{general-re-bound}
\end{eqn}
such that selected columns $\Ge_k$ satisfy an $r,\ve$ bound if $k \geq k^*$, with $\Ge_k$ satisfying Definition~\ref{def:greedy} for objective $\Lk_K$ and reflecting the $k$-column selection by nuclear maximization with error $\Rand$. 
\label{cor:additive-re-bound}
\end{restatable}

We may compare \eq{general-re-bound} to the tabulated results in \cite{Chen2022-lp} covering the expected performance of randomized alternatives.
As mentioned, DPP sampling achieves $r,\ve$ additive approximation in $k^* \lesssim \frac{r}{\ve} + r - 1$ columns.
Defining $\eta \Eq 1 - \Tp{r}{K}/\Tr[K] = \nu \Tp{r}{K}/\Tr[K]$, 
ridge leverage score sampling achieves the same in $k^* \lesssim 136 \frac{r}{\ve} \log \lrp{\frac{3 r}{\ve \eta}}$ columns. Diagonal sampling (i.e., RPCholesky) achieves the same in $k^* \lesssim \frac{r}{\ve} + r + \log_+ \lrp{\frac{1}{\ve \eta}}$ columns.
In comparison, \eq{general-re-bound} gives a 
bound
with a comparable if slightly weaker logarithmic multiplier.

\subsection{Error bounds for graph Laplacians \label{s:laplacian-analysis}}

In this section, we analyze the error of our column selection approach for graph Laplacian reduction.
Specifically, we establish approximation bounds for Algorithm~\ref{alg:exact-laplacian} (exact nuclear maximization) and Algorithm~\ref{alg:randomized-laplacian} (the matrix-free extension).
As before, we seek to bound the following objective \eq{laplacian-objective} for an $n \times n$ rescaled Laplacian $L$ given greedily selected subset $\Ge_k$ (Definition~\ref{def:greedy}):
\begin{eqn}
    \Lap_L^h(\Ge_k) = \Tr[(K^2)_{\I, \I} (K_{\I,\I})^\Mo] - \frac{1+h_{\I}^\t (K_{\I,\I})^\Mo (K^2)_{\I,\I} (K_{\I,\I})^\Mo h_{\I}}{h_{\I}^\t (K_{\I,\I})^\Mo h_{\I}} \; \Big|_{\I = \Ge_k}
\end{eqn}
where $k \geq 1$, $\Rand \geq 0$, $L h = 0$, $\norm{h}_2 = 1$, and $K=L^+$.

In one way, our task turns out to be easier, because in Theorem~\ref{th:laplacian-is-submodular}, we show that our objective $\Lap_L^h$ is submodular on the space of non-empty sets.
However, because our core interest is in graph Laplacians possessing a unique stationary state, we shall again have to take particular care to handle this feature.
If the graph Laplacian's random walk were killed (i.e., $L$ non-singular), in fact, we could immediately conclude by using the approach of Section~\ref{s:submodularity}, substituting (1) $L^\Mo$ in for $K$ and (2) the optimal $s$-subset's objective $\Lk_{L^\Mo}(\Opt_s)$ in for the DPP expectation $\dpp_s(K)$.
As $\Lk_K(\Opt_s) \geq \dpp_s(K)$ for all positive semidefinite $K$, this would only strengthen the results of Section~\ref{s:kernel-error} for this subcase.
Note, on the other hand, that the alternative objective $\Lk_{L^+}$, which is \emph{not} the same as our target objective from $\Lap_L^h$, is not submodular in general (as we have verified empirically). Thus if we were to consider $\Lk_{L^+}$, the results of Section~\ref{s:kernel-error} would be the best we could do. We will not consider this objective any further.

Given the stationary state, the structure of our optimization demands that we break apart our error analysis into one component, deriving from the first column selected, and a second, deriving from all of the other columns selected.
Procedurally, we shall first establish the submodularity of $\Lap_L^h$ on all non-empty subsets in Theorem~\ref{th:laplacian-is-submodular}.
Roughly speaking, this will control the deviation from optimality of all but the first column chosen.
Thereafter, we will variationally bound the value of the first column in Lemma~\ref{lem:laplacian-initial-bound}, yielding a complete bound in Theorem~\ref{th:laplacian-trace-bound} relating the greedily chosen $k$-subset to the optimal $k$-subset and the trace of $L$'s pseudoinverse $\Tr[L^+]$.

For brevity, we defer all proofs to Appendix~\ref{a:laplacian-error} but outline our results here.
First, we may rigorously define submodularity as follows:

\begin{definition}[Submodular function]
    A function $\subm$ is submodular over a space $\Psi$ if, for all sets satisfying
    $A \subseteq \Psi$, $B \subseteq \Psi$, and $C \subseteq \Psi$, it holds that 
\begin{eqn}
    \subm(A \cup B \cup C) - \subm(A \cup B) \leq \subm(A \cup C) - \subm(A).
    \req{submodular-def}
\end{eqn}
Equivalently, we can demand that for any set $A \subset \Psi$ and elements $i,j \in \Psi$, it holds that
\begin{eqn}
    \subm(A \cup \lrb{i,j}) - \subm(A \cup \lrb{i}) \leq \subm(A \cup \lrb{j}) - \subm(A).
\end{eqn}
\label{def:submodular}
\end{definition}

Given the next theorem, the objective function $\Lap_L^h$ is submodular on the space of all nonempty sets.
Appendix~\ref{a:laplacian-error} presents its proof, which relies on the continuous time Markov chain interpretation of a graph Laplacian:
\begin{restatable}[Submodularity of $\Lap_L^h$ for graph Laplacians]{theorem}{LaplacianSubmodularity}
    For $n \times n$ graph Laplacian $L$, $\Lap_L^h$ is submodular for all nonempty sets of indices $\{\I \, | \,  \I \subset \{1, \dots, n\}, \,  \I \neq \emptyset \}$.
    That is, for all sets $A,B,C \subset \{1, \dots, n\}$ where $A$ is nonempty:
    \begin{eqn}
        \TrL{\Comp{A}} - \TrL{\Comp{A\cup B}} - \TrL{\Comp{A \cup C}}  + \TrL{\Comp{A\cup B\cup C}} \geq 0
        \req{submodular-trace}
    \end{eqn}
\label{th:laplacian-is-submodular}
\end{restatable}


While Theorem~\ref{th:laplacian-is-submodular} proves the submodularity of $\Lap_L^h$ over the space of all non-empty subsets, the lack of proof for empty subsets prevents the trivial application of submodularity results.
We do not expect that any such proof could be given, in fact, in light of the difference in formulas for choosing the first column versus the latter columns, i.e., \eq{first-column-laplacian} versus \eq{latter-columns-laplacian}.

Thankfully, we may still establish a strong bound for nuclear maximization using the same generic linear program of Theorem~\ref{th:submodular}.
To do so, we essentially work with (1) trace-based bounding of the error upon selection of the first column and (2) submodularity to bound the near-optimality of the rest of the chosen columns.
To start with, we simply upper bound the error upon selection of the first column as follows:

\begin{restatable}[Bounding of objective after one iteration]{lemma}{LaplacianFirstBound}
Consider a rescaled Laplacian $L$ with stationary eigenvector $h$ as defined in Definition~\ref{def:laplacian}.
Suppose that $\frac{(L^+)_{\j,\j}}{h^2_\j} \leq \Orand \min_i \left\{ \frac{(L^+)_{i,i}}{h^2_i} \right\}$ for some index $\j$. Then 
\begin{eqn}
    \Lap_L^h(\lrb{\j}) \geq \Minus \Orand \Tr[L^+]
    \req{one-iter-3}
\end{eqn}
\label{lem:laplacian-initial-bound}
\end{restatable}

Next, by applying submodularity relations to the subsequent column selection and utilizing Theorem~\ref{th:submodular}, we obtain the following set of bounds:

\begin{restatable}[Guaranteed error bounds for graph Laplacians]{theorem}{LaplacianTraceBound}
    Consider the set $\Ge_k$ of size $k$ derived from greedy maximization of $\Lap_L^h$ up to maximum relative error $\Rand$ during the estimation of scores at each iteration.
    Letting $\Opt_s$ denote an $s$-subset which exactly maximizes $\Lap_L^h$, $1 \leq s \leq k$, it holds that:
\begin{eqn}
    \Lap_L^h(\Opt_s) - \Lap_L^h(\Ge_k) 
    \leq (2 + \Rand) \Tr[L^+] e^{-(k-1)/(s \Orand)}
    \req{laplacian-trace-bound}
\end{eqn}
\label{th:laplacian-trace-bound}
\end{restatable}

Thus, for our target problem of graph Laplacian reduction, nuclear maximization yields a subset of columns with objective that is at least exponentially close to that of optimal subset selection in its worst case.

\section{Computational experiments \label{s:benchmarks}}

In this section, we present the results of computational experiments for our proposed algorithms.
In Section~\ref{s:pathological}, we examine a few adversarially constructed toy problems for which diagonal maximization, uniform sampling, and diagonal sampling perform poorly relative to our approaches.
These cases allow us to draw some intuitive conclusions about when we expect the various algorithms to perform well or not.
Approaching somewhat more realistic scenarios, in Section~\ref{s:benchmark-kernel}, we present the results of kernel approximation for a few squared exponential Gaussian process kernels. 
Next, we present results for CUR decompositions of a set of large sparse matrices in Section~\ref{s:benchmark-cur}.
Finally, we close with an examination of our algorithms for rescaled Laplacians induced by DNA secondary structure kinetic landscapes, including results for sparse Laplacians from about $n=10^3$ to $n=10^6$.

Overall, our investigations support the use of our algorithms, and we find that the matrix-free implementations using randomized diagonal approximation generally match the deterministic results with exact diagonal evaluation.
Comparing selection strategies, we find that:
\begin{enumerate}
    \item Nuclear maximization generally performs as well or better than the other methods.
    \item Diagonal maximization performs poorly on adversarially designed examples but often performs well on our real-world problems.
    \item Uniform sampling is generally ineffective, especially on larger real-world problems.
    \item While diagonal sampling outperforms diagonal maximization for some designed examples, it performs worse than diagonal maximization in our real-world examples. 
    (However, see \cite{Chen2022-lp,Diaz2023-kl} for real-world problems where the opposite conclusion holds.) 
    From what we have observed, the statistical fluctuations of the results of diagonal sampling do not have a major impact on its practical performance.
    \item On the same designed examples mentioned above, the performance of nuclear maximization is still the best in our experiments, suggesting that a greedy approach can in fact work best across all regimes.
\end{enumerate}

When we apply diagonal maximization and sampling to the Laplacian reduction problem, we use the scores $h_{\j}^2$ for the first column selection and $(\hat{K}(\I))_{\j,\j}$ otherwise, cf. \eq{khat}. These choices are based on the limits of conventional diagonal maximization and sampling strategies as applied to $(L + \gamma h h^\t)^{-1}$ when $\gamma \rightarrow 0$.

\begin{figure}[t]
    \centering
    \includegraphics[width=6.5in]{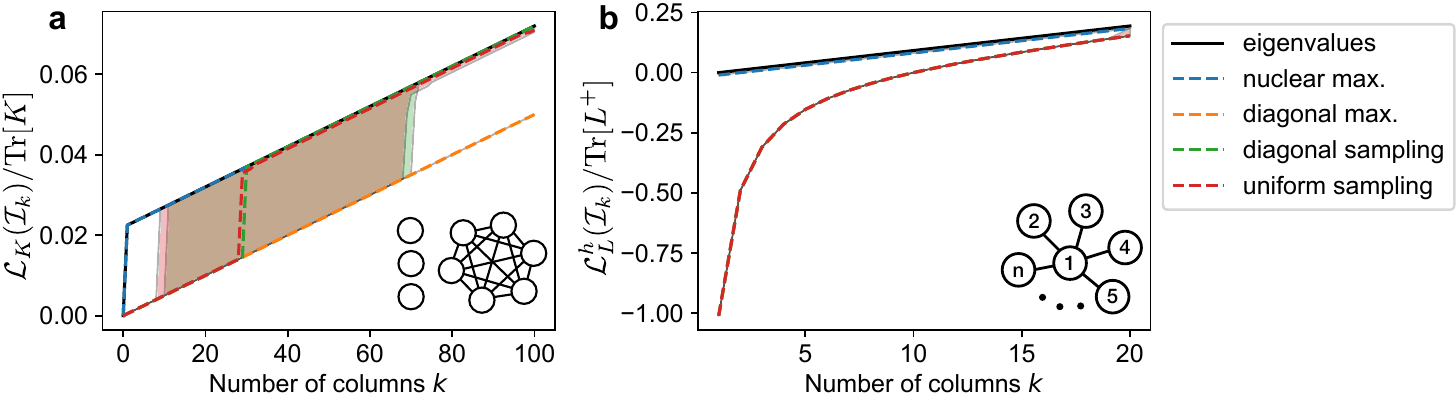}
    \caption{
        Adversarial examples for kernel and graph Laplacian reduction. 
        (a) Relative objective values $\Lk_K(\I) / \Tr[K]$ for the example \eq{pathological1} with $n=2000,n_c=45,\alpha=1.00001$.
        Nuclear maximization matches the eigenvalue bound, while each other method performs substantially worse.
        Diagonal maximization is worst of all, picking outlier nodes at each iteration, but uniform and diagonal sampling perform worse as well (and are visually coincident).
        (b) Relative objective values $\Lap_L^h(\I) / \Tr[L^+]$ obtained for the star graph Laplacian construction \eq{pathological-laplacian} with $n=100,\beta=0.9999$.
        Nuclear maximization matches the eigenvalue bound, while the other methods are qualitatively inferior (and are visually coincident).
        Convergence in the latter methods is substantially slower since the central node is picked only after a number of iterations.
        In each panel, for randomized column selection approaches, medians are plotted in 20\% and 80\% quantile envelopes over 1000 replicates.
        \label{fig:pathological-examples}
    }
\end{figure}

\subsection{Explicit adversarial examples \label{s:pathological}}

In this section, we only consider exactly performed (deterministic) score estimation, though recall that uniform and diagonal sampling have inherent randomization.
We construct two illustrative examples which, while pathological, exemplify when different selection strategies might fail.
Beyond the elementary theoretical analysis that we perform, we show numerical results in Figure~\ref{fig:pathological-examples} comparing the different selection strategies.

\subsubsection{Kernel and CUR approximation}

Motivated by a kernel matrix for a set of points consisting of a (1) a subset of $n_d$ disconnected outlier  points and and (2) a subset of $n_c$ tightly packed points, consider 
\begin{eqn}
	K = \begin{pmatrix}
		\alpha &  &  &  &  \\
		 & \ddots &  &  & \\
		 & & \alpha &  &  \\
		 &  & & 1 & \cdots & 1 \\
		 &  & & \vdots & \ddots & \vdots \\
		 &  & & 1 & \cdots & 1
	\end{pmatrix} = A A^\t 
    \ \text{ where } \ 
    A = \begin{pmatrix}
		\sqrt{\alpha} &  &    \\
		 & \ddots &   \\
		 & & \sqrt{\alpha} &  &  \\
		 &  & & 1 &  \\
		 &  & & \vdots  \\
		 &  & & 1 & 
	\end{pmatrix}
    \req{pathological1}
\end{eqn}
for $\alpha = 1 + \delta$ where $\delta > 0$ is viewed as small. We have taken $\alpha > 1$, only slightly, to break ties in diagonal maximization.
Clearly, the eigenvalues of $K$ are $n_c, \alpha, \dots, \alpha$.
For this problem, it is elementary to work out the average achieved objective $\Lk_K$ of each strategy for the first column selection.
In particular, nuclear maximization properly selects one of the connected nodes.
On the other hand, diagonal maximization selects one of the disconnected nodes, meaning that:
\begin{eqn}
    \frac{\Lk_{K}( \{ \j_\text{nuclear max} \} )}{\Lk_{K}( \{ \j_\text{diagonal max} \} )} = \frac{n_c}{\alpha} \approx n_c
\end{eqn}
Note that this ratio can be made arbitrarily large by increasing $n_c$.

Next, consider diagonal sampling, and more specifically take $n_c \sim \sqrt{n}$.
Then in the large $n$ limit, for $\alpha \approx 1$:
\begin{eqn}
    \frac{\Lk_{K}(\{ \j_\text{nuclear max} \} )}{\Ex [\Lk_{K}( \{ \j_\text{diagonal sample}) \} ]} \sim 
    \frac{n \left(\alpha-2 \left(\sqrt{n+1-\alpha}+1\right)\right)}{\alpha^2-4 n} \sim 
    \frac{n}{2 \sqrt{n}-1}
    \sim \tfrac{1}{2} \sqrt{n}
\end{eqn}
As we would expect, diagonal sampling protects partially against the pathology of the test problem, but it is still far from optimal.
Numerical results for multi-column selection are shown in Figure~\ref{fig:pathological-examples}a, validating the superiority of nuclear maximization for this test problem.

\subsubsection{Inverse Laplacian rank reduction}

Consider the following graph Laplacian $\bar L$ yielded by a connected star graph on $n$ nodes:
\begin{eqn}
    \bar L &\Eq \begin{pmatrix}
        n-1 & \Mo & \cdots & \cdots & \Mo  \\
        \Mo &  1 & &  &  \\
        \vdots &   & \ddots & &  \\
        \vdots &   &  & \ddots & \\
        \Mo &   &  &  & 1
    \end{pmatrix} 
\end{eqn}
The eigenvalues of $\bar L$ are $n, 1, \dots, 1, 0$.
To break ties in score selection, consider the slightly rescaled Laplacian $L$ determined by: 
\begin{eqn}
    h &= \begin{pmatrix} \beta & 1 & \cdots & 1 \end{pmatrix}^\t / \sqrt{n - 1 + \beta^2} \\
    L &= \Diag(h^\Mo)\bar L \Diag(h^\Mo) 
    \req{pathological-laplacian}
    \end{eqn}
Again we view $\beta = 1 - \delta$ where $\delta > 0$ is small.

It is elementary to verify that nuclear maximization correctly selects the central node for its first column.
On the other hand, by \eq{first-column-laplacian}, the initial probability density for diagonal sampling is $h^\otwo$ for this problem, because the diagonal of $(L + \Shift h h^\t)^\Mo$ is asymptotically proportional to $h^\otwo$ in the $\Shift \rightarrow 0$ limit. 
Thus diagonal maximization selects a column from the latter ($i>1$) set of nodes first, while uniform and diagonal sampling are identical for $\beta \approx 1$.
Therefore, recalling that our objective $\Lap_L^h$ is negative for a single column \eq{first-column-laplacian}, we may calculate that:
\begin{eqn}
    \frac{\Lap_L^h(\lrb{\j_\text{diagonal max}})}{\Lap_L^h(\lrb{\j_\text{nuclear max}})} &= \frac{\beta^4+n^2+2 \beta ^2 (n-2)-3 n+2}{n-1} \sim n \\
    \frac{\Ex[\Lap_L^h(\lrb{\j_\text{diagonal sample}})]}{\Lap_L^h(\lrb{\j_\text{nuclear max}})} &\sim\frac{\Ex[\Lap_L^h(\lrb{\j_\text{uniform sample}})]}{\Lap_L^h(\lrb{\j_\text{nuclear max}})} = \frac{\beta^4+n^2+2 \beta ^2 (n-2)-3 n+3}{n} \sim n
\end{eqn}
Note that the numerators and denominators on the left-hand sides are all negative, so these asymptotics show that nuclear maximization is dramatically superior for the first column selection in the large-$n$ limit.

Numerical results for multi-column selection are shown in Figure~\ref{fig:pathological-examples}b.
In this context, all three alternatives perform worse than nuclear maximization, which closely follows the eigenvalue bound.

\subsection{Experiments for kernel approximation \label{s:benchmark-kernel}}

\begin{figure}[h]
    \centering
    \includegraphics[width=\textwidth]{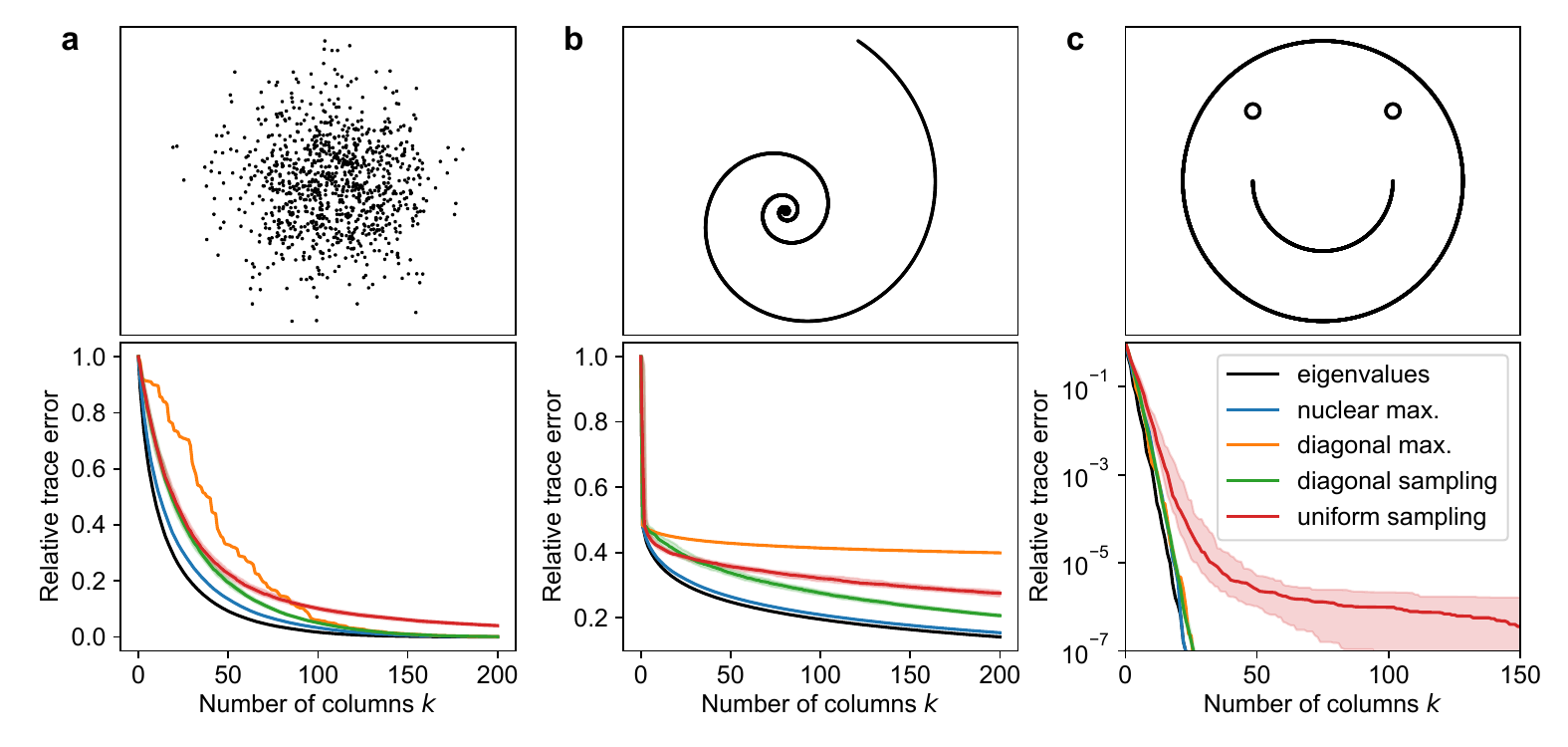}
    \caption{
        Illustrative examples of relative trace error ($1 - \Lk_K(\I)/\Tr[K]$) derived from kernel approximation.
        (a) A squared exponential kernel ($\sigma=0.4$) on a cloud of 1000 randomly distributed points from $\mathcal{N}(\Zero_2, \Id_2)$.
        (b) From \cite{Chen2022-lp}, a squared exponential kernel ($\sigma=10^3$) over $10^4$ points arranged in a spiral $(e^{t/5} \cos(t), e^{t/5} \sin(t))$ for evenly spaced $t$ in $[0,64]$).
        (c) From \cite{Chen2022-lp}, a squared exponential kernel ($\sigma=2$) over 10,000 points distributed uniformly in a 2D smiley face (50 per eye, 1,980 for the smile, and 7,920 for the outline).
        In each subfigure, uniform and diagonal sampling results are plotted according to their medians over 100 replicates, with error envelopes from 20\% and 80\% quantiles depicted as well.
        \label{fig:kernel-examples}
    }
\end{figure}

In this section, we present column selection results for three kernels, depicted in Figure~\ref{fig:kernel-examples}. The latter two examples are taken from \cite{Chen2022-lp}.
Each kernel matrix is derived from a squared exponential Gaussian process in two spatial dimensions.
These examples demonstrate the clear pitfalls of diagonal maximization and uniform sampling.
For these problems, nuclear maximization performs the best of all methods that we consider.

\paragraph{Randomly distributed point cloud:}
The points defining the kernel matrix in the first example are drawn from the two-dimensional standard normal distribution (Figure~\ref{fig:kernel-examples}a).
Points are more closely packed near the origin, with sparser coverage away from the origin.
The furthest points from the origin are ``outliers.''
We observe that nuclear maximization performs best of all methods considered, closely matching the lower bound provided by the eigenvalues of the kernel matrix.
Diagonal maximization performs quite poorly, selecting outlier points until none remain before selecting points near the origin.
Buy contrast, uniform sampling predominantly selects points from the central region, achieving satisfactory initial convergence but ultimately oversampling this region, resulting in slow convergence for large $k$.
Diagonal sampling nicely compromises between the uniform sampling and diagonal maximization strategies but performs worse than nuclear maximization.

\paragraph{Spiral:}
For our second example, we consider a manually prescribed set of points arranged in a spiral, following the construction in \cite{Chen2022-lp}.
Points towards the outer limits of the spiral are essentially complete outliers, with negligible covariance even to neighboring points.
In such a construction, diagonal maximization performs particularly poorly, selecting only outliers and yielding slow convergence. 
Uniform and diagonal sampling are closer to optimality, but nuclear maximization performs best, again closely approximating the top eigenvalues of the kernel matrix.

\paragraph{Smiley face:}
For our final example, we consider the smiley face construction of \cite{Chen2022-lp}.
With the layout described in that reference, the covariance pattern is highly heterogeneous, resulting in poor performance for uniform sampling. 
Each of the other discussed methods performs similarly, resulting in very quick convergence.

\subsection{Experiments for CUR decomposition \label{s:benchmark-cur}}

\begin{table}
    \centering
    \begin{tabular}{clrrrcl}
    \toprule
    System & SuiteSparse ID &  $m$ & $n$ & nnz & Symmetry & Origin \\
    \midrule
    a         & \texttt{bayer01}     & 57,735 & 57,735 & 275,094    & no &    \href{https://sparse.tamu.edu/Grund/bayer01}{chemical simulation}     \\
    b         & \texttt{bcsstk36}    & 23,052 & 23,052 & 1,143,140  & yes   &   \href{https://sparse.tamu.edu/Boeing/bcsstk36}{structural }     \\
    c         & \texttt{c-67b}       & 57,975 & 57,975 & 530,583    & yes &    \href{https://sparse.tamu.edu/Schenk_IBMNA/c-67b}{nonlinear optimization}     \\
    d         & \texttt{c-69}        & 67,458 & 67,458 & 623,914    & yes &    \href{https://sparse.tamu.edu/GHS_indef/c-69}{nonlinear optimization}      \\
    e         & \texttt{cbuckle}     & 13,681 & 13,681 & 676,515    & yes &    \href{https://sparse.tamu.edu/TKK/cbuckle}{structural}      \\
    f         & \texttt{crankseg\_2} & 63,838 & 63,838 & 14,148,858 & yes &     \href{https://sparse.tamu.edu/GHS_psdef/crankseg_2}{structural}      \\
    g         & \texttt{ct20stif}    & 52,329 & 52,329 & 2,600,295  & yes &    \href{https://sparse.tamu.edu/Boeing/ct20stif}{structural}      \\
    h         & \texttt{g7jac200sc}  & 59,310 & 59,310 & 717,620    & no &    \href{https://sparse.tamu.edu/Hollinger/g7jac200sc}{economics}      \\
    i         & \texttt{venkat01}    & 62,424 & 62,424 & 1,717,792  & no   &   \href{https://sparse.tamu.edu/Simon/venkat01}{fluid dynamics}     \\
    \bottomrule
    \end{tabular}
    \caption{Details of $m \times n$ sparse matrices used to evaluate selection methods for CUR decomposition.
    We mostly consider the same examples as \cite{Ekenta2022-qb} from the SuiteSparse matrix collection \cite{Kolodziej2019-fd}.}
    \label{tab:suitesparse}
\end{table}

\begin{figure}
    \centering
    \includegraphics[width=\textwidth]{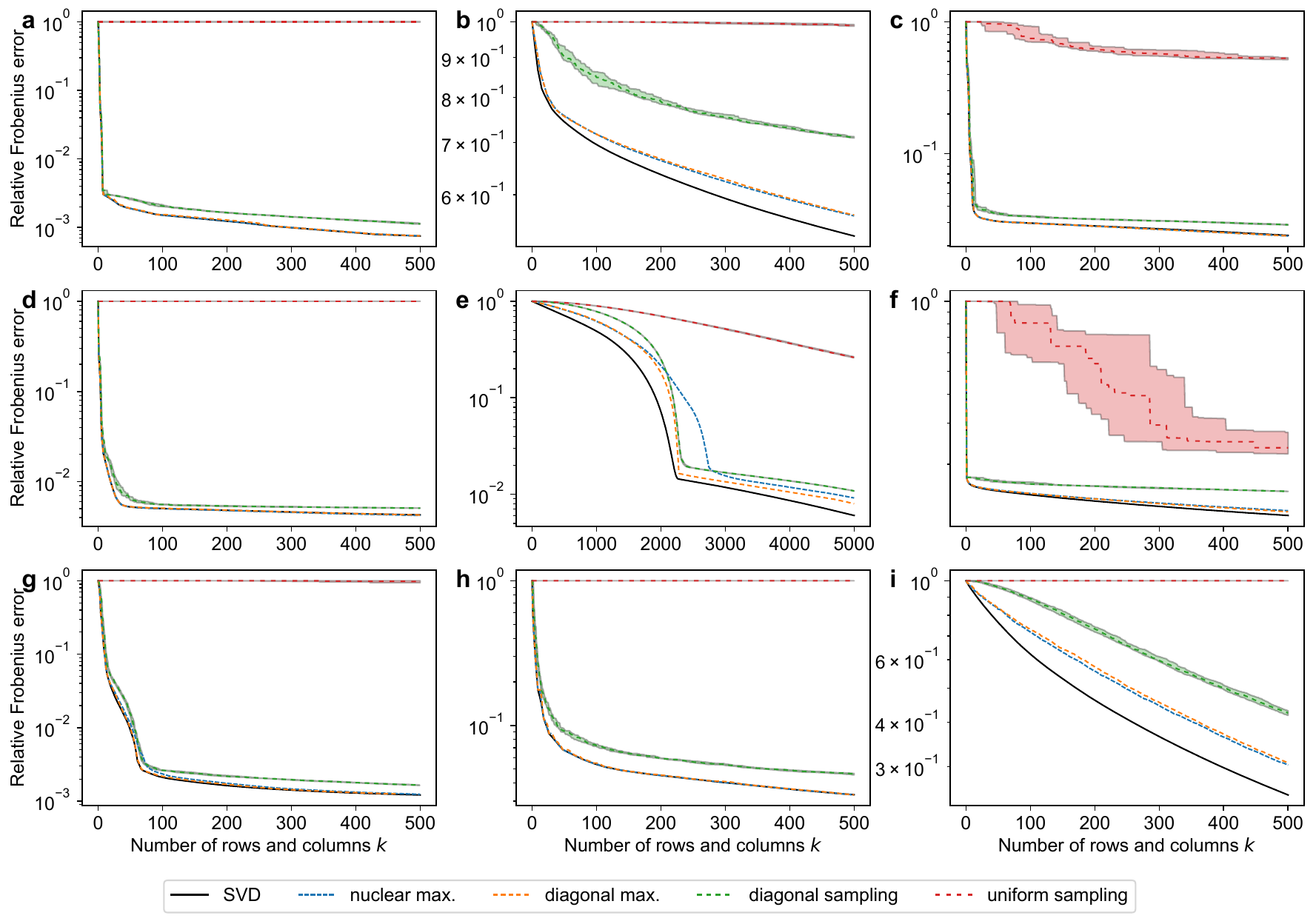}
    \caption{
        Relative Frobenius error $\norm{CUR - A}_F / \norm{A}_F$ for collected examples of CUR decomposition of sparse matrices from Table~\ref{tab:suitesparse} using deterministic scoring (Algorithm~\ref{alg:exact-cholesky}).
        For uniform and diagonal sampling, medians are plotted with error envelopes from 20\% and 80\% quantiles, based on results over 10 trials.
        Top singular values were computed using subspace iteration. 
        \label{fig:cur-suitesparse-deterministic}
    }
\end{figure}

\begin{figure}
    \centering
    \includegraphics[width=\textwidth]{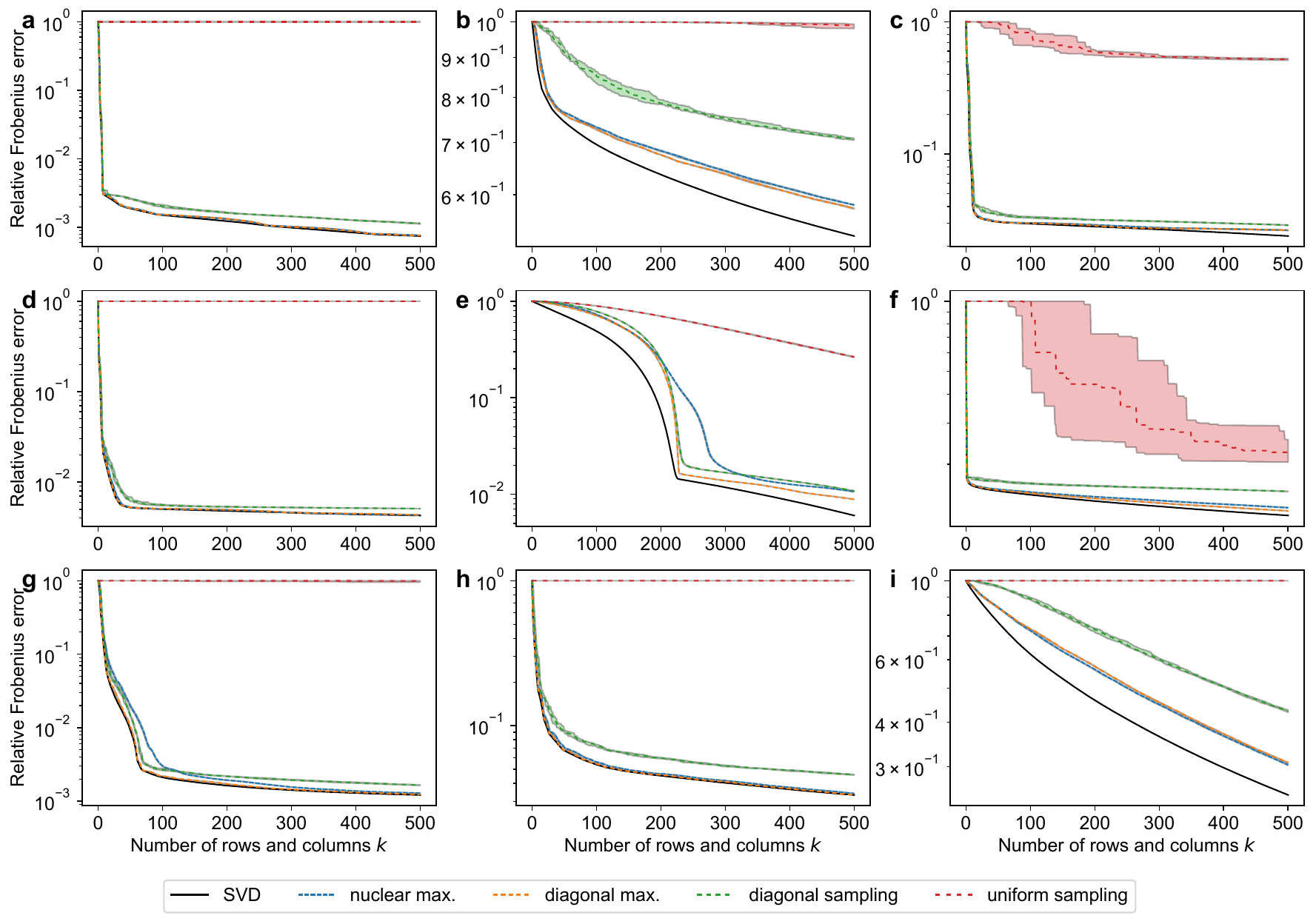}
    \caption{
        Matrix-free equivalent to Figure~\ref{fig:cur-suitesparse-deterministic}, using Algorithm~\ref{alg:randomized-cholesky} and $z=200$ for each randomized diagonal approximation. 
        Medians are plotted in a 20\% and 80\% quantile envelope for all non-SVD results over 10 trials.
        \label{fig:cur-suitesparse-random}
    }
\end{figure}

For CUR decompositions, we follow the procedure outlined in Section~\ref{s:cur-objective} and always measure error in the Frobenius norm. 
For the column and row selection subroutines, we consider both the deterministic approach of Section~\ref{s:exact-implementation} and the matrix-free approach of Section~\ref{s:randdiag}. 
For the deterministic approach, we precompute the $\Diag(K^2)$ and $\Diag(K)$ for $K = A^\t A$ and $K=A A^\t$ to perform column and row selection, respectively. 
For the matrix-free approach, this step is unnecessary.

We study nine sparse matrices $A$ drawn from the SuiteSparse (\cite{Kolodziej2019-fd}) matrix library; see Table~\ref{tab:suitesparse} for details. 
We chose these examples to follow the same benchmark as in \cite{Ekenta2022-qb} to the extent that they were well-specified enough to replicate, i.e., we did not choose them by hand.
These sparse matrices display a range of compressibility; some matrices (e.g., \texttt{bayer01}) are approximately low-rank, while others (e.g., \texttt{cbuckle}) have singular values which decay slowly.
We show results from deterministic selection in Figure~\ref{fig:cur-suitesparse-deterministic} and matrix-free selection in Figure~\ref{fig:cur-suitesparse-random}. Note that our matrix-free algorithm is actually novel for the diagonal maximization and diagonal sampling strategies, not just nuclear maximization. 
Overall, the matrix-free results match the deterministic results quite well, as anticipated by theory, and the matrix-free approach is to be preferred for its superior computational scaling.

In every case, uniform sampling is highly suboptimal, to the point of offering almost no effective error decrease.
Nuclear and diagonal maximization perform similarly, and the two methods are visually coincident in many cases.
We have not attempted to find SuiteSparse examples for which diagonal maximization performs poorly, and in many real-world examples, particularly those without strong outlier contributions, we expect that diagonal maximization \emph{may} perform quite well.
On the other hand, diagonal maximization has inferior worst-case performance guarantees and cannot be expected to succeed in all circumstances.
Finally, diagonal sampling yields a compromise between uniform sampling and diagonal maximization which often results in substantially worse practical performance than the latter.

\subsection{Experiments for graph Laplacians \label{s:benchmark-laplacian}}

Finally, we consider matrix-free methods for the optimization of $\Lap_L^h$ for physically derived systems corresponding to the normalized graph Laplacians associated with continuous-time Markov chains modeling discrete DNA secondary structure chemical kinetics.
With each graph edge corresponding to the addition or deletion of a single base pair, the graph Laplacians of such systems are relatively sparse.
Because detailed balance is assumed, the normalized graph Laplacian is symmetric with stationary eigenvector $h = \pi^{\odot 1/2}$, where $\pi$ is the normalized Boltzmann distribution corresponding to the equilibrium probability of each secondary structure in the space.

In such cases, while the generator $L$ is sparse, its pseudoinverse $K$ is completely dense,
making its computation and storage particularly demanding.
As such, for the largest systems that we consider, we investigate only matrix-free methods.
In contrast with the case for dense kernel matrices, there is no complexity difference among diagonal maximization, diagonal sampling, and nuclear maximization. Note again that our matrix-free algorithm is actually novel for the diagonal maximization and diagonal sampling strategies, not just nuclear maximization. 
Nuclear maximization necessitates only $z$ additional $K$ matvecs per column selection, but this cost is small compared to the $z$ matvecs by $C$ performed via Chebyshev expansion.
Of course, uniform sampling is much cheaper than diagonal maximization, diagonal sampling, and nuclear maximization, but the the results are so poor that we mostly omit them from our results.

\begin{table}[t]
    \centering
    \begin{tabular}{crrrccc}
    \toprule
    System  & $n$ & nnz$(L)$ & nnz$(R)$ & $\kappa_0$ & $\kappa$ & $n_\mathrm{cheb}$ \\
    \midrule
    1          & 702 &  4,750 & 5,466 & 2.97e3 & 73.25 & 99 \\
    2          & 26,045 &  252,107& 329,925 & 2.57e4 & 49.36 & 81 \\
    3          & 1,087,264  &  12,519,798 & 17,989,447  & 3.37e5 & 217.1 & 176  \\
    \bottomrule
    \end{tabular}
    \caption{Details of graph Laplacian systems.
        Each Laplacian $L$ is an $n \times n$ matrix with $\mathrm{nnz}(L)$ nonzero entries and condition number $\kappa_0$.
        Given the approximate factorization  $L \approx R R^{\t}$, $\kappa$ is the condition number of the preconditioned matrix $R^{\pinv} L R^{\pinv \t}$,
        and $n_\mathrm{cheb}$ is the number of Chebyshev nodes used to approximate matvecs by $(C^\pinv L C^{\pinv\t})^{\pinv/2}$ with guaranteed relative error $10^{\Minus 8}$. Condition numbers were estimated using power and inverse power iteration.
        See Appendix~\ref{a:computational-details} for details.
    }
    \label{tab:laplacian-inputs}
\end{table}

\begin{figure}[b!]
    \centering
    \includegraphics[width=\textwidth]{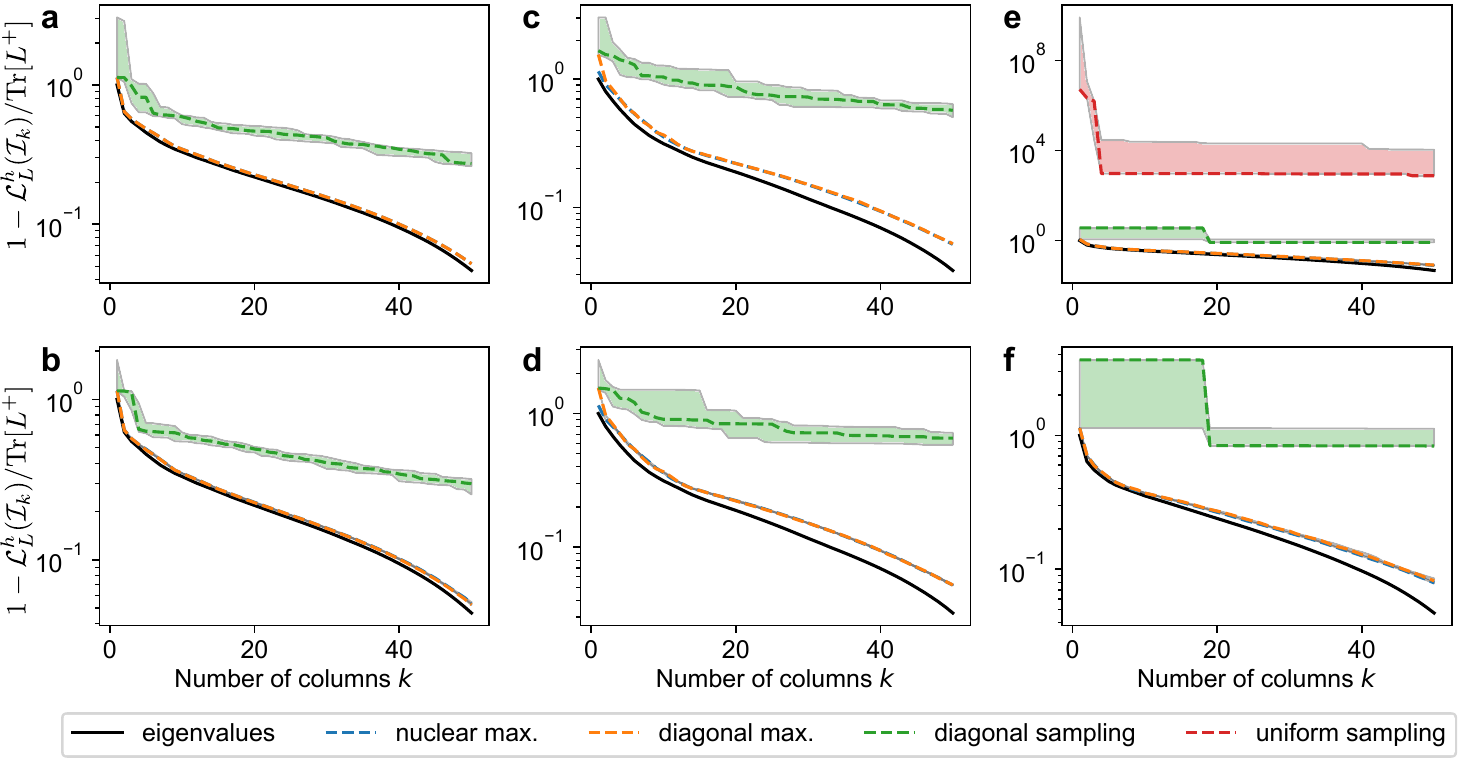}
    \caption{
        Results for example Laplacians described in Table~\ref{tab:laplacian-inputs}, including
        System 1, using exact inversion (a) or matrix-free approximation (b),
        System 2, using exact inversion (c) or matrix-free approximation (d),
        and System 3, using matrix-free approximation including (e) or excluding (f) uniform sampling results.
        20\%, 50\%, and 80\% quantiles over 30 replicates are depicted in (a-d);
        minima, medians, and maxima over 3 replicates are shown in (e-f).
        The normalized trace error is plotted in each panel. 
        When omitted, errors for uniform sampling were substantially outside the plotting window. We provide panel (e) as an illustrative example of its performance.
        \label{fig:laplacian-results}
    }
\end{figure}

We consider three example problems, presented in Table~\ref{tab:laplacian-inputs}, each induced by a randomly generated DNA strand of some length (sequences in Table~\ref{tab:laplacian-sequences}).
For each sequence, an associated continuous time Markov chain was constructed over all compatible unpseudoknotted secondary structures \cite{Flamm2000-dq}.
We utilized the \texttt{dna04} parameters in NUPACK 4 \cite{Fornace2022-ar} at default (37\textdegree C, 1 M Na$^+$) conditions in the \texttt{nostacking} ensemble; with these parameters, only Watson-Crick (\texttt{A}$\cdot$\texttt{T} and \texttt{C}$\cdot$\texttt{G}) base pairs are considered \cite{Fornace2020-gs}.
By construction, the stationary state probability of structure $s$ is proportional to its Boltzmann factor $\pi(s) \propto \exp(\Minus \Delta G(s) / k_B T)$, with $\Delta G(s)$ defined by the NUPACK 4 free energy model and \texttt{dna04} parameters.
Transition rates were generated according to the Kawasaki rate function, i.e. with rate $r(s_i \rightarrow s_j) \propto \sqrt{\pi(s_j) / \pi(s_i)} = \h_j / \h_i$ for neighboring secondary structures $s_i$ and $s_j$ \cite{Schaeffer2015-cz}.

In order to apply our approach, we used the \texttt{rchol} software package and algorithm \cite{Chen2020-ya}, conceptually upon that of \cite{Kyng2016-qv}, to provide an approximate Cholesky factorization. 
Table~\ref{tab:laplacian-inputs} summarizes the results of this preconditioning. 
For the systems that we studied, the approximate factorization yielded a preconditioned condition number $\kappa$ ranging from 50 to 200, cf. Table~\ref{tab:laplacian-inputs}.

For the cases benchmarked in Figure~\ref{fig:laplacian-results}, nuclear and diagonal maximization similarly and yield closed approximations to the eigenvalue-derived lower bound.
Diagonal sampling produces slower convergence than either of the maximization methods. 
In Figure~\ref{fig:laplacian-results}ab, the smaller system sizes permit us to perform each algorithm using exact, deterministic calculation of column scores.
These results confirm that (for our choice $z=200$), the matrix-free algorithms closely match their deterministic analogs.

\section*{Acknowledgments}
This work was partially supported by the Applied Mathematics Program of the US Department of Energy (DOE) Office of Advanced Scientific Computing Research under contract number DE-AC02-05CH11231 (M.F. and M.L.).

\clearpage
\newrefcontext[sorting=anyt]
\printbibliography[heading=bibintoc, title={Bibliography}]

\clearpage
\appendix

\section{Detailed algorithms \label{a:alg-details}}

\subsection{Deterministic kernel column selection \label{s:alg-exact-cholesky}}

In this section, we present our optimized deterministic algorithm to perform column selection via nuclear maximization.
As discussed in Section~\ref{s:cur-objective}, one may perform CUR decomposition of matrix $A$ using the same algorithm with $K=A A^\t$ and $K=A^\t A$.
We use $\odot$ to denote elementwise multiplication.


\input{algorithms/exact-kernel}

Throughout Algorithm~\ref{alg:exact-cholesky}, as $\I$ is augmented, $w$ is the diagonal of $\tK^2(\I)$, while $d$ is the diagonal of $\tK(\I)$.
Upon completion, $U$ satisfies $U_{:,\I} U_{:,\I}^\t = (K_{\I,\I})^\Mo$, while $S$ satisfies $S=K U^\t$ and $S S^\t = K_{:,\I} (K_{\I,\I})^\Mo K_{\I,:}$. 
The computational complexity of Algorithm~\ref{alg:exact-cholesky} is $O(n k)$ in space and $O(n k^2 + k \psi + \phi)$ in time, where $\psi$ is the cost of performing a matvec by $K$ and $\phi$ is the cost of computing $\Diag(K^2)$ and $\Diag(K)$. 
For general (possibly sparse) $K$ supporting $O(1)$ element access, the time complexity is therefore $O(n k^2 + k \nnz(K))$.
Compared to standard Cholesky factorization, the extra cost of Algorithm~\ref{alg:exact-cholesky} is dominated by the single matvec per iteration $K S_{:,t}$ in the final line.

\subsection{Deterministic inverse Laplacian rank reduction \label{s:alg-exact-laplacian}}

In Algorithm~\ref{alg:exact-laplacian}, we present a method for deterministic column selection using nuclear maximization with respect to an input rescaled Laplacian $L$ with stationary eigenvector $h$.
This algorithm is not intended to scale well to large graph Laplacians, as it makes use of explicit computations of the diagonals of $K$ and $K^2$. 
(Recall that in this setting, $K \Eq L^+$.)
On the other hand, it may be useful for modest system sizes or other contexts.

\input{algorithms/exact-laplacian}

Supposing that we compute $K \Eq L^+$ ahead of time (which could involve $O(n^3)$ computation), we will assume that $L^+$ is generally dense (incurring $O(n^2)$ storage by itself).
Beyond these requirements, the computational complexity of Algorithm~\ref{alg:exact-laplacian} is $O(k n^2)$, and it requires $O(k n)$ storage.

\subsection{Matrix-free kernel column selection \label{s:alg-randomized-cholesky}}

In Algorithm~\ref{alg:randomized-cholesky}, we present an algorithm for kernel approximation via matrix-free column selection using nuclear maximization.
This algorithm takes as inputs SPSD $K$ and $C$ satisfying $C C^\t = K$.
However, these operators may be accessed only via matvecs, meaning that neither operator needs to actually be pre-computed.

\input{algorithms/randomized-kernel}

(Note that the select access above of $K_{\j,\j}$ may be performed via a single matvec $K_{\j,\j} = \Id_{\j,:} K \Id_{:,\j} $ at each iteration.)
Overall, Algorithm~\ref{alg:randomized-cholesky} performs a subset of the computations in Algorithm~\ref{alg:exact-cholesky} at the cost of including a randomized score approximation step. 

\subsection{Matrix-free inverse Laplacian rank reduction \label{s:alg-randomized-laplacian}}

In Algorithm~\ref{alg:randomized-laplacian}, we present a method for matrix-free column selection using nuclear maximization with respect to an input rescaled Laplacian $L$ with stationary eigenvector $h$, given pseudoinverse $K \Eq L^+$ amd exact factorization $K = C C^\t$.
While operators $K$ and $C$ are taken as inputs, only matvec access is needed for them, meaning that neither operator needs to actually be pre-computed.

\input{algorithms/randomized-laplacian}

The time complexity of Algorithm~\ref{alg:randomized-laplacian} is $O(n k^2 z + k z \psi)$, where  $\psi$ is the greater cost of performing a matvec by $K$ or $C$.
In particular, Algorithm~\ref{alg:randomized-laplacian} uses only matvec access to operators $K$ and $C$, enabling the use of iterative methods with inexact factorizations, as described in Section~\ref{s:sqrt-alg}.

\subsection{Explanation of computational complexities \label{s:complexity-details}}

Here we explain further the computational complexities stated in Table~\ref{tab:complexities}.

\paragraph{Deterministic kernel approximation} (Algorithm~\ref{alg:exact-cholesky}): \newline
Updating $U$ incurs $O(k n)$ cost in every iteration, yielding $O(k^2 n)$ expense overall.
Additionally, a single full matvec is needed in every iteration, yielding $O(k \nnz(K))$ expense.
A partial matvec is also required, which is upper bounded by this expense.

\paragraph{Matrix-free kernel approximation} (Algorithm~\ref{alg:randomized-cholesky}): \newline
The complexities from deterministic approximation still hold, except that $z$ matvecs by $K$ are performed at each iteration, yielding $O(k z \nnz(K))$ expense.
Projections of the randomized scores using past column choices yield an additional $O(n k^2 z)$ expense.
Generally speaking, randomization does not help in this particular context if $K$ has been precomputed (as a sparse or dense matrix).
It is primarily useful when $K$ possesses more structure (such as a factorization) or when a matrix-free approach is required.

\paragraph{Deterministic CUR decomposition} (Algorithm~\ref{alg:exact-cholesky} with $K = A A^\t$ and $K=A^\t A$): \newline
Updates of $U$ require $O(m k^2)$ work for row selection and $O(n k^2)$ for column selection, yielding $O((m+n) k^2)$ combined complexity.
In the worst case, $A A^\t$ might be fully dense, making the calculation of $\Diag(K^2)$ via $(A A^\t)^\otwo \One$ require $O(\nnz n)$ or $O(\nnz m)$ cost. 
(Depending on the sparsity pattern of $A$, this might be reduced.)
A single matvec by $K = A A^\t$ or $K = A^\t A$ is required at each iteration, yielding $O(k \nnz(A))$ complexity.

\paragraph{Matrix-free CUR decomposition}  (Algorithm~\ref{alg:randomized-cholesky} with $K = A A^\t$ and $K=A^\t A$): \newline
Relative to deterministic factorization, the diagonals no longer need to be computed exactly, removing the $O(\nnz (m + n))$ expense, and $z$ matvecs are performed at each iteration, incurring $O(k z \nnz(A))$ cost.
Projections of the randomized scores using past row and column choices yield an additional $O((m + n) k^2 z)$ expense.

\paragraph{Deterministic Laplacian reduction} (Algorithm~\ref{alg:exact-laplacian}):\newline
We shall assume that $K$, a dense matrix, has been precomputed in $O(\vartheta)$ time.
(This might incur $O(n^3)$ complexity in itself.)
Since $K$ is generally dense, $O(k)$ matvecs incur $O(k n^2)$ expense.
As noted in the main text, Algorithm~\ref{alg:exact-laplacian} is primarily useful for smaller system sizes, owing to its lack of scalability.

\paragraph{Matrix-free Laplacian reduction}  (Algorithm~\ref{alg:randomized-laplacian}): \newline
We shall assume a sparse Cholesky factorization $L \approx R R^\t$ with $\nnz(R) \sim \tilde{O}(\nnz(L))$ and that $R^+ L R^{+\t}$ has condition number $O(1)$ in its image.
With these caveats, Chebyshev and conjugate gradient computations require $O(1)$ iterations, and matvecs by $R$, $R^+$, and $L$ are each computable in $\tilde{O}(\nnz(L))$ time, yielding a cost of $\tilde{O} (k z \nnz(L))$ overall.
Updates of $U$ require $O(n k^2)$ work, and projections of the randomized scores for subsequent iterations require $O(n k^2 z)$ work.

\subsection{Computational details \label{a:computational-details}}

All computational studies were run on an Apple M3 Max Macbook Pro (36 GB RAM, 10 performance cores, 4 efficiency cores).
Algorithms were prototyped and verified in Python, then implemented and benchmarked in C++17.
The armadillo matrix library (version 12.6.6) was used for dense and sparse matrix linear algebra routines and data structures \cite{Sanderson2016-ud,Sanderson2019-ac}.
Graph Laplacian systems were preprocessed using the minimum degree ordering algorithm from the Boost Graph library \cite{Siek2001-ed}, after which the \texttt{rchol} library \cite{Chen2020-ya} was used to compute the necessary approximate sparse Cholesky factorizations.
The Intel oneAPI Threading Building Blocks library \cite{Reinders2007-gi} was used to implement parallelism.
NUPACK 4 \cite{Fornace2022-ar} was used to sample random DNA sequences and to derive the studied graph Laplacians.

For graph Laplacians, preconditioned conjugate gradient was used for $L^+$ matvecs, and Chebyshev expansion was used for $(R^+ L R^{+\t})^{+/2}$ matvecs.
In both cases, the stationary mode corresponding to $h$ was explicitly projected out during each iteration of the algorithm to prevent accumulation of numerical error in this mode.
For deterministic algorithms, due to numerical floating point precision, the ratio $\frac{(\tK^2(\I))_{\j,\j}}{\tK_{\j,\j}(\I)} = \frac{((K - K_{:,\I} (K_{\I,\I})^\Mo K_{\I, :})^2)_{\j,\j}}{(K - K_{:,\I} (K_{\I,\I})^\Mo K_{\I, :})_{\j,\j}}$ was observed to be error-prone in cases where both the numerator and denominator were very close to 0.
(Essentially, we believe this to be a case of precision loss in the computation of $(a-b)/(c-d)$ when $a \approx b$ and $c \approx d$.)
It was assumed that the correct exact ratio for such cases would be close to 0.
Thus, in our deterministic algorithms, entries $\j$ for which $\tK_{\j,\j}(\I) < 10^{\Minus 8} K_{\j,\j}$ were not considered for selection.
(This issue did not appear in matrix-free algorithms, so no modifications were implemented therein.)

CUR decompositions were computed independently for row and column selection as described in the text.
Let $S_R$ and $U_R$ denote the computed matrices of Algorithm~\ref{alg:exact-cholesky} or Algorithm~\ref{alg:randomized-cholesky} for row selection ($K = A A^\t$), and let $S_C$ and $U_C$ be those matrices computed for column selection ($K = A^\t A$).
Finally, let $\I$ and $\J$ denote the rows and columns computed, respectively.
Then the Frobenius error may be efficiently computed as
\begin{eqn}
    \norm{A - A_{:, \J} (A_{:, \J}^+ A A_{\I, :}^+) A_{\I, :}}_F^2 = \One^\t A^\otwo \One - \One^\t (S_R^\t A_{:,\J} U_C)^\otwo \One,
\end{eqn}
a formula we used in our benchmark evaluation. 
(A corresponding formula can be worked out involving $U_R$ and $S_C$.
Note that the usage of $U$ here is in the context of upper Cholesky factorization, i.e. reflecting the matrices in the mentioned algorithms, rather than the U in the CUR name.)
One may also work out that the approximate factorization yielded is equal to $A_{:,\J} U_{C} U_C^\t A^\t_{:,\J} S_R U_R^\t A_{\I,:}$ or $A^{+\t} S_C S_C^\t A^+ S_R S_R^\t A^{+\t}$, among other possible expressions.

\begin{table}
    \centering
    \begin{tabular}{crc}
        \toprule        
        System & Sequence & Length (nt) \\
        \midrule
        1 & \texttt{TAACTAGGCTTCAAGTCGTC} & 20 \\
        2 & \texttt{TAACTAGGCTTCAAGTCGTCGACGCGA} & 27 \\
        3 & \texttt{ATACTTTGTCGGCGTAGTTGTCGACTGCTCTTCTTTTCCC} & 40 \\
        \bottomrule
    \end{tabular}
    \caption{
        DNA sequences for investigated graph Laplacians.
    }
    \label{tab:laplacian-sequences}
\end{table}

For our studies of graph Laplacians derived from DNA secondary structure ensembles, respective sequences are given in Table~\ref{tab:laplacian-sequences}; see Table~\ref{tab:laplacian-inputs} for more details and computed quantities.
For these systems, subspace iteration was used to calculate $\kappa_0$ (the condition number of $L$) by separately converging the top twelve eigenvalues of $L$ and $L^+$ until each of those eigenvalues changed by less than $10^{\Minus 8}$ between iterations.
Similarly, subspace iteration was used to calculate $\kappa$ (the condition number of $R^{\pinv} L R^{\pinv \t}$) by separately converging the top twelve eigenvalues of $R^{\pinv} L R^{\pinv \t}$ and $R L^+ R^\t$ until each of those eigenvalues changed by less than $10^{\Minus 8}$ between iterations.

\section{Objective formulation details \label{a:objective-formulation}}

In this section we provide proofs and explanations of the objectives formulated in our work (Section~\ref{s:objectives}).

\subsection{Variational eigenvalue principle \label{a:eig-proof}}

Here we give a short (and elementary) proof showing that unconstrained Frobenius error minimization (or maximization of $\Lk_K$) yields the same result as a truncated eigendecomposition.

\begin{proposition}
Suppose that $K$ is a SPSD $n \times n$ matrix of rank at least $k$, possessing eigenvalues $\lambda_1 \geq \cdots \geq \lambda_n$.
Then:
\begin{eqn}
	\max_{F \in \mathbb{R}^{n,k}} \Tr[F^\t K^2 F (F^\t K F)^\Mo] = \sum_{i=1}^k \lambda_i
    \req{ky-fan-1}
\end{eqn}
\label{prop:ky-fan}
\end{proposition}
\begin{proof}
We will invoke the reduced $QR$ decomposition of $K^{1/2} F$, which exists without restriction \cite{Golub2013-an}.
Then $R$ is $(k, k)$, invertible, and upper triangular, and $Q$ is $(n, k)$ with $Q^\t Q = \Id$.
Therefore:
\begin{eqn}
	\Tr[F^\t K^2 F (F^\t K F)^\Mo] &= \Tr[R^\t Q^\t K Q R (R^\t Q^\t Q R)^\Mo]
	= \Tr[Q^\t K Q]
\end{eqn}
We may thus translate the maximization in \eq{ky-fan-1} to one over orthonormal matrices $Q$:
\begin{eqn}
    \max_{F \in \mathbb{R}^{n,k}} \Tr[F^\t K^2 F (F^\t K F)^\Mo] = \max_{Q \in \mathbb{R}^{n,k}, Q^\t Q = \Id} \Tr[Q^\t K^2 Q (Q^\t K Q)^\Mo],
\end{eqn}
and the desired statement is then yielded by the Ky Fan theorem (e.g. \cite{Zhang2011-tv}, Theorem 8.17).
\end{proof}

\subsection{Trace augmentation identities \label{a:trace-split-column}}

Here, we provide proofs of how our objectives for kernel approximation $\Lk_K$ and inverse Laplacian reduction $\Lap_L^h$ may be efficiently maximized by choosing columns one by one.

\begin{proposition}
Suppose that $K$ is a SPSD matrix, and suppose that $K_{\I, \I}$ and $K_{\Ij, \Ij}$ are non-singular submatrices.
\begin{eqn}
    \Lk_K(\Ij) &= \Tr[(K^2)_{\Ij, \Ij}(K_{\Ij, \Ij})^\Mo] 
    = \Lk_K(\I) + \frac{(\tK^2(\I))_{\j,\j}}{(\tK(\I))_{\j,\j}}
\end{eqn}
\label{prop:split-trace}
\end{proposition}
\begin{proof}
    Note that $\tK_{\j,\j}(\I) = K_{\j,\j} - K_{\j,\I} (K_{\I,\I})^\Mo K_{\I,\j}$ by definition.
    From standard block matrix inverse identities \cite{Bernstein2009-bv}, we have that:
    \begin{eqn}
        K_{\Ij, \Ij})^\Mo = \begin{pmatrix}
        (K_{\I,\I})^\Mo + (K_{\I,\I})^\Mo K_{\I, \j} K_{\j,\I} (K_{\I,\I})^\Mo / \tK_{\j,\j}(\I) & \Minus (K_{\I,\I})^\Mo K_{\I, \j} / \tK_{\j,\j}(\I) \\
        \Minus K_{\j,\I} (K_{\I,\I})^\Mo / \tK_{\j,\j}(\I) & 1 / \tK_{\j,\j}(\I)
        \end{pmatrix}
    \end{eqn}
    Thus (exploiting $\Tr[A B] = \One^\t (A \odot B^\t) \One$, $\odot$ denoting Hadamard multiplication):
    \begin{eqn}
        \Tr[K_{:, \Ij}&(K_{\Ij, \Ij})^\Mo] \\
        &= \Tr[(K^2)_{\I, \I}(K_{\I, \I})^\Mo] \\ &+ \Tr \lrs{(K^2)_{\I \cup \j,\I \cup \j} \begin{pmatrix} 
            (K_{\I,\I})^\Mo K_{\I, \j} K_{\j,\I} (K_{\I,\I})^\Mo & \Minus (K_{\I,\I})^\Mo K_{\I, \j} \\
            \Minus K_{\j,\I} (K_{\I,\I})^\Mo & 1
        \end{pmatrix}} / \tK_{\j,\j}(\I) \\
        &= \Tr[(K^2)_{\I, \I}(K_{\I, \I})^\Mo] \\ &+ \frac{K_{\j,:} K_{:,\j} - 2 K_{\j,:} K_{:,\I} (K_{\I,\I})^\Mo K_{\I, \j} + K_{\j,\I} (K_{\I,\I})^\Mo K_{\I,:} K_{:,\I} (K_{\I,\I})^\Mo K_{\I,\j}   }{\tK_{\j,\j}(\I)} \\
        &= \Tr[(K^2)_{\I, \I}(K_{\I, \I})^\Mo] + \frac{(\tK^2(\I))_{\j,\j}}{\tK_{\j,\j}(\I)}
    \end{eqn}
    and the rest follows from the definition of $\Lk_K$.
\end{proof}


In the following proposition, we show that a similar augmentation identity holds for rescaled Laplacians, given that the augmented set is non-empty.

\begin{proposition}
Let $L$ be a rescaled Laplacian as defined in Definition~\ref{def:laplacian} with stationary eigenvector $h$, and let $\I$ be a non-empty set of indices.
Defining the following:
\begin{eqn}
    \hat{K}(\I) \Eq \tilde{K}(\I) + \frac{(h - K_{:,\I} (K_{\I, \I})^\Mo h_{\I})(h^\t - h_{\I}^\t (K_{\I, \I})^\Mo K_{\I,:})}{h_{\I}^\t (K_{\I, \I})^\Mo h_{\I}},
\end{eqn}
then for $\j \notin \I$:
\begin{eqn}
    \Lap_L^h(\Ij) = \Lap_L^h(\I) + \frac{(\hat{K}^2(\I))_{\j,\j}} {(\hat{K}(\I))_{\j,\j}}
\end{eqn}
\end{proposition}
\begin{proof}
From \eq{laplacian-complement-limit} and \eq{laplacian-objective},
\begin{eqn}
    \Lap_L^h(\Ij) - \Lap_L^h(\I) = \lim_{\Shift \rightarrow 0} \Lk_{K + \Shift^\Mo h h^\t}(\Ij) - \Lk_{K + \Shift^\Mo h h^\t}(\I)
\end{eqn}
Let $M \Eq (L+\Shift h h^\t)^\Mo = K + \Shift^\Mo h h^\t$, $\Shift > 0$. 
Then, letting $\tilde M(\I) \Eq M - M_{:, \I} (M_{\I,\I})^\Mo M_{\I, :}$, Proposition~\ref{prop:split-trace} yields that:
\begin{eqn}
    \Tr[M_{:, \Ij}(M_{\Ij, \Ij})^\Mo] - \Tr[(M^2)_{\I, \I}(M_{\I, \I})^\Mo] = \frac{(\tilde M^2(\I))_{\j,\j}}{\tilde M_{\j,\j}(\I)}
\end{eqn}
Next, upon use of the Sherman-Morrison formula and simplification:
\begin{eqn}
    \tilde M(\I) &= M - M_{:,\I} (M_{\I,\I})^\Mo M_{\I,:} \\
    &= K + \Shift^\Mo h h^\t - (K_{:,\I} + \Shift^\Mo h h_{\I}^\t) (K_{\I,\I}+\Shift^\Mo h_{\I} h_{\I}^\t)^\Mo (K_{\I,:} + \Shift^\Mo h_{\I} h^\t) \\
    &= K + \Shift^\Mo h h^\t - (K_{:,\I} + \Shift^\Mo h h_{\I}^\t) ((K_{\I,\I})^\Mo + \frac{(K_{\I,\I})^\Mo h_{\I} h_{\I}^\t (K_{\I,\I})^\Mo}{\Shift + h_{\I} (K_{\I,\I})^\Mo h_{\I}}) (K_{\I,:} + \Shift^\Mo h_{\I} h^\t) \\
    &= K - K_{:,\I} (K_{\I,\I})^\Mo K_{\I,:} + \frac{(h - K_{:,\I} (K_{\I,\I})^\Mo h_{\I})(h^\t - h_{\I} (K_{\I,\I})^\Mo K_{\I,:})}{\Shift + h_{\I}^\t (K_{\I,\I})^\Mo h_{\I}}
\end{eqn}
Taking the limit of $\Shift \rightarrow 0$ yields that $\tilde M(\I) \rightarrow \hat K(\I)$ , giving the desired result.
\end{proof}

\subsection{CUR decomposition \label{a:cur-objective}}

Here we give full proofs of error relations for CUR decomposition, which are elementary.

\TraceCUR*
\begin{proof}
    \begin{eqn}
        \epsilon_A(\J) &= \norm{A_{:, \J} (A^\t_{:,\J} A_{:,\J})^\Mo A_{:, \J}^\t A - A}_F^2 \\
        &= \Tr[A^\t A] - 2\Tr[A^\t A_{:, \J} (A^\t_{:, \J} A_{:, \J})^\Mo A_{:, \J}^\t A] \\
        & \quad \quad \quad + \  \Tr[A^\t A_{:, \J} (A^\t_{:, \J} A_{:, \J})^\Mo A^\t_{:, \J} A_{:, \J} (A^\t_{:, \J} A_{:, \J})^\Mo A_{:, \J}^\t A] \\
        &= \Tr[A^\t A] - \Tr[A^\t A_{:, \J} (A^\t_{:, \J} A_{:, \J})^\Mo A_{:, \J}^\t A] \\
        &= \Tr[K] - \Tr[K_{:, \J} (K_{\J, \J})^\Mo K_{\J, :}]
    \end{eqn}
    where we have expanded the definition of the Frobenius norm and substituted the definition of $K$. The last expression is precisely $\mathcal{E}_K (\J)$, following \eq{error-def}.
\end{proof}

\IndependentCUR*
\begin{proof}
    Compute: 
    \begin{eqn}
        \norm{A - C C^+ A R^+ R}_F
        &\leq \norm{A-  C  C^+ A}_F + \norm{C C^+ A - C  C^+ A R^+ R}_F \\
        &\leq \norm{A-  C  C^+ A}_F + \norm{ C  C^+}_2 \norm{A - A R^+ R}_F \\ 
        &\leq \norm{A-  C  C^+ A}_F + \norm{A - A R^+ R}_F \\
        &=  \epsilon_{A}^{1/2} (\mathcal{J}) +  \epsilon^{1/2}_{A^\t} (\mathcal{I})
    \end{eqn}
    using the triangle inequality, the fact that the spectral norm of a projector matrix $C  C^+$ is 1, and \eq{cur_err}. Then Lemma \ref{lem:cur-reduction} completes the proof.
\end{proof}

\subsection{Graph Laplacian reduction \label{s:markov}}

In this section we will elaborate further upon the random walk interpretation of the rescaled Laplacian as used in Section~\ref{s:laplacian-analysis}.
(We shall pursue this connection further in future work, so we give only the necessary basics here.)
Suppose we are given rescaled Laplacian $L$ and elementwise positive eigenvector $h$ satisfying $L h$ and $\norm{h}_2 = 1$.
Furthermore, suppose that the graph corresponding to $L$'s adjacency matrix is connected.
Given these assumptions, $L$ has a single zero eigenvector ($h$) and is symmetric positive semidefinite.

We shall consider the rate matrix (e.g. \cite{Aldous1995-bg}) as formed via $R \Eq \Minus \Diag(h^{\Mo}) L \Diag(h)$.
Let $\pi \Eq h^\otwo$, such that $\One^\t \pi = \One^\t h^\otwo = \norm{h}_2^2 = 1$.
Then the rate matrix $R$ satisfies:
\begin{eqn}
    R_{i,i} &< 0 & \forall i \\
    R_{i,j} &\geq 0 & \forall i \neq j \\
    \pi_i R_{i,j} &= \pi_j R_{j,i} & \forall i, j \\
    \pi^\t R &= \Zero^\t \\
    R \One &= \Zero
\end{eqn}
Then $R$ may be interpreted as the generator of a reversible, recurrent, time-homogenous continuous time Markov chain (CTMC) \cite{Aldous1995-bg}, and $\pi$ as its associated stationary distribution. Conversely, any such $R$ yields a suitable rescaled Laplacian $L$.

\subsubsection{Interpretation of inverses of Laplacian submatrices \label{a:markov-interpretation}}

Given the continuous time Markov chain interpretation above, the probability of a system starting in state $i$ being in state $j$ after time $t \geq 0$ is then:
\begin{eqn}
    P_{i,j}(t) = (\exp(R t))_{i,j}
\end{eqn}
Consideration of the series expression for the matrix exponential yields that
\begin{eqn}
    \exp(\Minus L_{\I, \I} t) = \Diag(h) \exp(R_{\I, \I} t) \Diag(1/h)
\end{eqn}
and therefore
\begin{eqn}
    (L_{\I,\I})^\Mo = \int_0^\infty \exp(\Minus L_{\I, \I} t) \odif t = \Minus \Diag(h) (R_{\I,\I})^\Mo \Diag(1/h) = \int_0^\infty \Diag(h) \exp(R_{\I, \I} t) \Diag(1/h)  \odif t
\end{eqn}
so that $\frac{h_j}{h_{i}} ((L_{\I,\I})^\Mo)_{i,j}$ is the expected time that a system starting in state $i$ will spend in state $j$ before leaving subset of states $\I$.
Thus, $\Tr[(L_{\Comp \I, \Comp \I})^\Mo]$ measures a weighted average of the time a system starting in a given state spends in that same state before reaching any of the states in $\I$. This provides an interpretation for \eq{laplacian-min} that we will develop more in future work. 


\subsubsection{Schur complement formulation for rescaled Laplacians \label{a:laplacian-schur}}

In the following lemma, we show how a Schur complement form can be still be posed to minimize $\Tr[(L_{\Comp \I, \Comp \I})^\Mo]$ even for a rescaled Laplacian $L$ possessing a single stationary mode.

\LaplacianComplement*
\begin{proof}
    We shall show the desired result by simplifying the following regularized inverse and taking the limit afterwards:
    \begin{eqn}
        \Tr[(L_{\Comp \I, \Comp \I})^\Mo] = \lim_{\Shift \rightarrow 0} \Tr[ ( (L + \Shift h h^\t)_{\Comp \I, \Comp \I} )^\Mo]
        \req{laplacian-complement-limit}
    \end{eqn}
    Noting that $(L + \Shift h h^\t)^\Mo = K + \Shift^\Mo h h^\t$ for finite $\Shift$, we may use the block matrix Schur complement as follows:
    \begin{eqn}
        \Tr[((L + \Shift h h^\t)_{\Comp \I, \Comp \I})^\Mo] &= \Tr[(L + \Shift h h^\t)^\Mo] - \Tr[((L + \Shift h h^\t)^2)_{\I,\I} ((L + \Shift h h^\t)_{\I,\I})^\Mo] \\
        &= \Tr[K] + \Shift^\Mo - \Tr[((L^2)_{\I,\I} + \Shift^2 h_{\I} h_{\I}^\t) (L_{\I,\I} + \Shift h_{\I} h_{\I}^\t)^\Mo]
    \end{eqn}
    The rightmost trace may be reduced via the Sherman-Morrison formula to find (after simplification) that:
    \begin{eqn}
        \Tr[((L + \Shift h h^\t)_{\Comp \I, \Comp \I})^\Mo]
        = \Tr[K] - \Tr[(K^2)_{\I, \I} (K_{\I,\I})^\Mo] + \frac{1+h_{\I}^\t (K_{\I,\I})^\Mo (K^2)_{\I,\I} (K_{\I,\I})^\Mo h_{\I}}{\Shift + h_{\I}^\t (K_{\I,\I})^\Mo h_{\I}}
    \end{eqn}
    and therefore taking the limit $\Shift \rightarrow 0$ in \eq{laplacian-complement-limit} trivially yields \eq{laplacian-complement}.
\end{proof}


\subsubsection{Matrix-free factorization \label{a:laplacian-factor}}

Finally, we provide a more in-depth explanation of our diagonal estimator used for inverse Laplacian reduction in noninitial iterations \eq{laplacian-complicated-estimator}.
For non-empty subset $\I$, recall our definitions \eq{laplacian-complicated-defs} that:
\begin{eqn}
    \tau(\I) \Eq h_{\I}^\t (K_{\I,\I})^\Mo h_{\I} \qquad
    y(\I) \Eq h - K_{:,\I} (K_{\I,\I})^\Mo h_{\I}
\end{eqn}
For the denominator, omitting $(\I)$ where implied and recalling that $\tilde C \tilde C^\t = \tilde K$ \eq{Ctilde},
\begin{eqn}
    \Diag(\hat K) &= \Diag\lrp{\tilde K + y y^\t / \tau} \\ 
    &= \Diag(\tilde K) + y^\otwo / \tau \\
    &= \Ex_{x \sim \mathcal{N}(\Zero,\Id)}[(\tilde C x)^\otwo \One] + y^\otwo / \tau
\end{eqn}
Next, for the numerator, working backwards and using $\tilde{K} h = \Zero$ and $\h^\t h = 1$:
\begin{eqn}
    \Ex_{x \sim \mathcal{N}(\Zero,\Id)}&[[(\tilde K + \tau^\Mo y (y-h)^\t) x]^\otwo] + y^\otwo / \tau^2 \\
    &= \Diag((\tilde K + \tau^\Mo y (y-h)^\t) (\tilde K + \tau^\Mo (y-h) y^\t) + y y^\t / \tau^2) \\
    &= \Diag(\tK \tK + y y^\t \tK / \tau + \tK y y^\t / \tau + y y^\t y y^\t /\tau^2) \\
    &= \Diag((\tK + y y^\t / \tau)^2) = \Diag(\hat K^2)
\end{eqn}
such that:
\begin{eqn}
    \frac{(\hat K^2)_{\j,\j}}{\hat K_{\j,\j}} 
    = \frac{\Ex_{x \sim \mathcal{N}(\Zero,\Id)}[[(\tilde K + \tau^\Mo y (y-h)^\t) x]^\otwo] + y^\otwo / \tau^2}
    {\Ex_{x \sim \mathcal{N}(\Zero,\Id)}[[\tilde C x]^\otwo] + y^\otwo / \tau}
    \req{laplacian-vector-estimator}
\end{eqn}
Note that the non-expectation right-hand terms in \eq{laplacian-vector-estimator} are trivially element-wise non-negative.

\section{Properties of the $k$-DPP expectation \label{a:dpp-proofs}}

In this section, we will prove elementary facts about $\dpp_k(K)$ \eq{dpp-def}, i.e., the expectation value $\Ex_{\I \sim k\text{-DPP}(K)}[\Lk_K(\I)] = \Ex_{\I \sim k\text{-DPP}(K)}[\Tr[(K^2)_{\I,\I} (K_{\I,\I})^\Mo]]$ yielded by $k$-DPP sampling.
First, in Appendix~\ref{a:dpp-monotonic}, we show that $\dpp_k(K)$ is increasing with $k$, a simple consequence of Newton's inequalities.
In Appendix~\ref{a:dpp-concave}, we give our main contribution, a novel proof that $\dpp_k(K)$ is concave in $k$ (in a finite-difference sense), a fact which is unproven in the literature (to our knowledge).
In Appendix~\ref{a:dpp-schur-convex}, we prove that $\dpp_k(K)$ is Schur-convex in the eigenvalues of $K$.
This statement is not novel, but we believe our proof is simpler than some previous ones \cite{Schur1923-di,Guruswami2012-rr}.
Finally, in Appendix~\ref{a:dpp-subadditive}, we prove the subadditivity of $\dpp_k(K)$ with respect to $K$, which follows straightforwardly from literature results.

To define our terms, for any given vector $x$ of $n$ real numbers, let $e_k(x)$ be the elementary symmetric polynomial of order $k$ in $x$:
\begin{eqn}
    e_k(x) \Eq \sum_{1 \leq i_1 < \cdots i_k \leq n} \prod_{j=1}^k x_{i_j}
\end{eqn}
and $E_k(x)$ be the associated normalized symmetric polynomial, $E_k(x) \Eq e_k(x) / \Binomial{n}{k}$.
Let also $e_{k}(x) = E_{k}(x) = 0$ for all $k < 0$ or $k > n$, and recall that for any $n \times n$ SPD matrix $A$, the $k$-DPP expectation is derived to be (Proposition~\ref{prop:dpp-sym}, e.g. \cite{Guruswami2012-rr}):
\begin{eqn}
    \dpp_k(A) \Eq e_1(x) - (k+1) \frac{e_{k+1}(x)}{e_k(x)}
\end{eqn}
where $x$ is the vector of the $n$ (positive) eigenvalues of $A$.
It is also helpful to have the following identities involving ratios of the elementary symmetric polynomials.
For any $x$ of $n$ positive real numbers, and any integer $0 \leq k \leq n$,  it follows from the definitions of $E$ and $e$ that
\begin{eqn}
	(k+1) \frac{e_{k+1}(x)}{e_k(x)} &= (n-k) \frac{E_{k+1}(x)}{E_k(x)} \\
	\frac{e_{k+1}(x)}{e_k(x)} &= \frac{e_{n-k-1}(x^\Mo)}{e_{n-k}(x^\Mo)} \\
	\frac{E_{k+1}(x)}{E_k(x)} &= \frac{E_{n-k-1}(x^\Mo)}{E_{n-k}(x^\Mo)}
	\req{sym-ratios}
\end{eqn}
where $x^\Mo$ is the vector formed from the reciprocals of each term in $x$ (elementwise).
For reference, we will also commonly use Newton's inequalities in the following sections:
\begin{theorem}[Newton's inequalities \cite{Maclaurin1729-mh}]
Let $x$ be a vector of $n > 0$ positive real numbers, and let $k$ be an integer such that $0 \leq k \leq n$.
Then:
\begin{eqn}
    E_{k}^2(x) \geq E_{k-1}(x) E_{k+1}(x)
\end{eqn}
and equality only holds if all of the elements of $x$ are identical.
In terms of the non-normalized elementary symmetric polynomials,
\begin{eqn}
    (k) (n-k) e_{k}^2(x) \geq (k+1) (n-k+1) e_{k-1}(x) e_{k+1}(x)
\end{eqn}
and so $e_{k}^2(x) > e_{k-1}(x) e_{k+1}(x)$, strictly.
\label{th:newton}
\end{theorem}

Finally, we will make use of the following elegant lemma from \cite{Rosset1989-xv}. See that reference for the proof, which, to summarize, relies on Rolle's theorem applied to differentiation of the characteristic polynomial of $\Diag(x)$.
\begin{lemma}[Elementary symmetric polynomial reduction \cite{Rosset1989-xv}]
Let $x$ be a non-empty vector of $n$ non-negative numbers.
Then there exists a vector $y$ of $n-1$ non-negative numbers such that $E_k(x) = E_k(y)$ for $0 \leq k < n$.
\label{lem:newton-rolle}
\end{lemma}


\subsection{Monotonicity with respect to $k$ \label{a:dpp-monotonic}}

For completeness, we provide a short proof that the $k$-DPP expectation monotonically increases (strictly) with $k$.
\begin{theorem}
Let $x$ be a vector of $n > 0$ positive real numbers, and let $k$ be a integer satisfying $0 \leq k < n$.
Then $\dpp_k(\Diag(x)) < \dpp_{k+1}(\Diag(x))$.
\label{th:dpp-monotonic}
\end{theorem}
\begin{proof}
Simply expand the definitions of $\dpp_k$ and apply Newton's inequalities (e.g. \cite{Hardy1952-el}).
\begin{eqn}
    \dpp_{k+1}(\Diag(x)) - \dpp_k(\Diag(x)) &= (k+1) \frac{e_{k+1}(x)}{e_k(x)} - (k+2) \frac{e_{k+2}(x)}{e_{k+1}(x)} \\
    &= (n-k) \frac{E_{k+1}(x)}{E_k(x)} - (n-k-1) \frac{E_{k+2}(x)}{E_{k+1}(x)} \\
    &= (n-k) \lrp{\frac{E_{k+1}(x)}{E_k(x)} - \frac{E_{k+2}(x)}{E_{k+1}(x)} } + \frac{E_{k+2}(x)}{E_{k+1}(x)}
\end{eqn}
The parenthesized term is non-negative by Newton's inequalities and zero only if (1) the all elements of $x$ are identical and 
(2) $k<n-1$ (if $k=n-1$, $E_{k+1}(x) / E_k(x) - E_{k+2}(x) / E_{k+1}(x) = E_{k+1}(x) / E_k(x) > 0$).
The remainder term on the right is trivially non-negative.
\end{proof}
\subsection{Concavity with respect to $k$ \label{a:dpp-concave}}

Less trivially, we will also prove that $k$-DPP is concave in $k$, in a finite difference sense.
We provide a complete proof of this identity here, which is novel to the best of our knowledge.
Our proof relies on an induction in increasing numbers of arguments $n$, with a base case proved for $n=3$ and edge cases handled for the extreme limits of low and high values of $k$.
We will begin first with the following lemma, which occurs as a base case for the inductive proof of our more general lemma.

\begin{lemma}
Let $x$ be a vector of exactly 3 positive real numbers.
Then:
\begin{eqn}
    \frac{E_{0}(x)}{E_1(x)}-2 \frac{2 E_1(x)}{E_{2}(x)}+\frac{3 E_{2}(x)}{E_{3}(x)} \geq 0
\end{eqn}
\label{lem:concave-edge}
\end{lemma}
\begin{proof}
Let $[x_1, x_2, x_3]$ be the elements of $x$.
Then we will simply apply the definition of $E$ and rearrange terms:
\begin{eqn}
	&\frac{E_{0}(x)}{E_1(x)}-2 \frac{2 E_1(x)}{E_{2}(x)}+\frac{3 E_{2}(x)}{E_{3}(x)} \\
    =& \frac{E_2(x) E_3(x) -4 E_3(x) E_1(x)^2+3 E_2(x)^2 E_1(x)}{E_1(x) E_2(x) E_3(x)} \\
	=& \frac{x_2^2 x_1^3+x_3^2 x_1^3-2 x_2 x_3 x_1^3+x_2^3 x_1^2+x_3^3 x_1^2-2 x_2 x_3^3 x_1-2 x_2^3 x_3 x_1+x_2^2 x_3^3+x_2^3 x_3^2}{E_1(x) E_2(x) E_3(x)} \\
	=& \frac{x_1^3 (x_2-x_3)^2 + x_2^3 (x_1-x_3)^2 + x_3^3 (x_1-x_2)^2 }{E_1(x) E_2(x) E_3(x)}
\end{eqn}
the right-hand side of which is clearly observed to be non-negative.
\end{proof}


We proceed with the following lemma, which addresses the most general restricted case of our target theorem:
\begin{lemma}
Let $x$ be a vector of $n \geq 3$ positive real numbers.
Let $l$ be an integer in $[1:n-2]$.
Then:
\begin{eqn}
    \frac{l E_{l-1}(x)}{E_l(x)}-2\frac{(l+1) E_l(x)}{E_{l+1}(x)}+\frac{(l+2) E_{l+1}(x)}{E_{l+2}(x)} \geq 0
\end{eqn}
\label{lem:sym-l}
\end{lemma}
\begin{proof}
We will proceed by induction on increasing $n$.
Thus, suppose we have proven our lemma for all vectors containing $n-1$ positive real numbers.
Now, suppose first that $l < n-2$.
Then by Lemma~\ref{lem:newton-rolle}:
\begin{eqn}
    \frac{l E_{l-1}(x)}{E_l(x)}-2\frac{(l+1) E_l(x)}{E_{l+1}(x)}+\frac{(l+2) E_{l+1}(x)}{E_{l+2}(x)} = \frac{l E_{l-1}(y)}{E_l(y)}-2\frac{(l+1) E_l(y)}{E_{l+1}(y)}+\frac{(l+2) E_{l+1}(y)}{E_{l+2}(y)}
\end{eqn}
for some $y$ of $n-1$ positive real numbers.
Thus this quantity is non-negative by the assumption of our induction.
Second, suppose that $l = n-2$.
Then we must show that:
\begin{eqn}
    \frac{(n-2) E_{n-3}(x)}{E_{n-2}(x)}-2\frac{(n-1) E_{n-2}(x)}{E_{n-1}(x)}+\frac{n E_{n-1}(x)}{E_{n}(x)} \geq 0
\end{eqn}
Let $r \Eq x^\Mo$; trivially, $r$ is a vector of $n$ positive numbers. Then by \eq{sym-ratios}:
\begin{eqn}
    \frac{(n-2) E_{n-3}(x)}{E_{n-2}(x)}-2\frac{(n-1) E_{n-2}(x)}{E_{n-1}(x)}+\frac{n E_{n-1}(x)}{E_{n}(x)} &= \frac{(n-2) E_{3}(r)}{E_{2}(r)}-2\frac{(n-1) E_{2}(r)}{E_{1}(r)}+\frac{n E_{1}(r)}{E_{0}(r)} \\
    &= \frac{e_1(r)}{e_0(r)}-2 \frac{2 e_2(r)}{e_1(r)}+\frac{3 e_3(r)}{e_2(r)}
\end{eqn}
Suppose that $n > 3$, and, without loss of generality, let $v$ denote the first $n-1$ elements of $r$ and $s$ denote the $n$th element of $r$.
Then:
\begin{eqn}
    \frac{e_1(r)}{e_0(r)}-2 \frac{2 e_2(r)}{e_1(r)}+\frac{3 e_3(r)}{e_2(r)} &= e_1(v) + s - 4 \frac{e_2(v) + s e_1(v)}{e_1(v) + s}+3\frac{e_2(v) + s e_2(v)}{e_3(v) + s e_1(v)} \\
    &= \frac{s (e_1(v) (s-e_1(v))+2 e_2(v))^2}{e_1(v) e_1(r) e_2(r)} + \frac{e_2(v)}{e_2(r)} \lrp{\frac{e_1(v)}{e_0(v)}-2 \frac{2 e_2(v)}{e_1(v)}+\frac{3 e_3(v)}{e_2(v)}}
\end{eqn}
The first term on the right-hand side is non-negative upon inspection, while the second term is non-negative via the assumption of our induction.
Finally, we can prove our base case for $n=3$, in which case our theorem governs only $l=1$.
The inequality is directly yielded by Lemma~\ref{lem:concave-edge}.
Thus we have proven our theorem for all $n \geq 3$.
\end{proof}

\begin{theorem}
Suppose $x$ is a vector of $n$ positive real numbers, and suppose $2 \leq k \leq n$.
Then:
\begin{eqn}
	\frac{(n-k) E_{k+1}(x)}{E_k(x)} -2 \frac{(n-k+1) E_k(x)}{E_{k-1}(x)} + \frac{(n-k+2) E_{k-1}(x)}{E_{k-2}(x)} \geq 0
\end{eqn}
\label{th:sym-convex}
\end{theorem}
\begin{proof} 
Let us first handle the case of $k=n$.
Then the inequality reduces to:
\begin{eqn}
	2 \frac{ E_{n-1}(x)}{E_{n-2}(x)} - 2 \frac{E_n(x)}{E_{n-1}(x)} \geq 0
\end{eqn}
which is true as a special case of Newton's inequalities.
Second, consider $2 \leq k \leq n-1$.
Define $l \Eq n-k$, such that $1 \leq l \leq n-2$.
Then from the reciprocal relation \eq{sym-ratios}:
\begin{eqn}
    &\frac{(n-k) E_{k+1}(x)}{E_k(x)} -2 \frac{(n-k+1) E_k(x)}{E_{k-1}(x)} + \frac{(n-k+2) E_{k-1}(x)}{E_{k-2}(x)} \\
    = & \frac{l E_{k+1}(x)}{E_k(x)} -2 \frac{(l+1) E_k(x)}{E_{k-1}(x)} + \frac{(l+2) E_{k-1}(x)}{E_{k-2}(x)} \\
    = & \frac{l E_{l-1}(x^\Mo)}{E_l(x^\Mo)} -2 \frac{(l+1) E_l(x^\Mo)}{E_{l+1}(x^\Mo)} + \frac{(l+2) E_{l+1}(x^\Mo)}{E_{l+2}(x^\Mo)}
    \req{sym-map}
\end{eqn}
where $x^\Mo$ contains the elementwise reciprocals of $x$.
$x^\Mo$ is a vector of $n$ positive numbers, and hence the right-hand side of \eq{sym-map} is non-negative by Lemma~\ref{lem:sym-l}.

\end{proof}

\begin{corollary}
Let $K$ be SPSD with eigenvalues $\lambda$. 
Then $\dpp_{j+1}(K) - \dpp_{j}(K) \leq \dpp_{j}(K) - \dpp_{j-1}(K)$ for all $0 < j < n$. 
\label{cor:dpp-concave}
\end{corollary}
\begin{proof} 
The proof is trivial upon substitution of the definitions of $\dpp_j$ and $E$.
Let $k=j+1$, such that $2 \leq k \leq n$.
\begin{eqn}
    2\dpp_j(K) - \dpp_{j+1}(K) - \dpp_{j-1}(K) 
    &= 2\dpp_j(\Diag(\lambda)) - \dpp_{j+1}(\Diag(\lambda)) - \dpp_{j-1}(\Diag(\lambda)) \\ 
    &= \frac{(k+1) e_{k+1}(\lambda)}{e_{k}(\lambda)} - 2 \frac{k e_{k}(\lambda)}{e_{k-1}(\lambda)} + \frac{(k-1) e_{k-1}(\lambda)}{e_{k-2}(\lambda)} \\
     &= \frac{(n-k) E_{k+1}(\lambda)}{E_k(\lambda)} -2 \frac{(n-k+1) E_k(\lambda)}{E_{k-1}(\lambda)} + \frac{(n-k+2) E_{k-1}(\lambda)}{E_{k-2}(\lambda)}
\end{eqn}
which is non-negative by Theorem~\ref{th:sym-convex}.
\end{proof}

\subsection{Schur-convexity \label{a:dpp-schur-convex}}

Generally, if one averages together two eigenvalues of a given input matrix $K$, the DPP expectation $\dpp_k(K)$ can only decrease, i.e., worsening its low-rank approximation quality.
Inspecting the definition of $\dpp_k(K)$ \eq{dpp-def}, we can simply show this by proving that the ratio $(k+1) e_{k+1}(x) / e_k(x)$ is Schur-concave.
This was proven by Schur in \cite{Schur1923-di}.
We shall give our own self-contained proof, which is rather simple.


\begin{theorem}[Schur-convexity of $k$-DPP]
    If $\lambda$ is a vector of $n$ non-negative numbers, then $\dpp_k(\Diag(\lambda))$ is Schur-convex in $\lambda$.
    Equivalently, consider two symmetric positive semidefinite $n \times n$ matrices $A$ and $B$ with eigenvalues $\lambda_A$ and $\lambda_B$, respectively.
    Then if $\lambda_A \prec \lambda_B$, then $\dpp_k(A) \leq \dpp_k(B)$ for $k = 0, \dots, n$.
    \label{th:schur-convexity}
\end{theorem}
\begin{proof}
The proof is trivial for $k=1$, so we shall consider only $k \geq 2$.
A sufficient condition for Schur-concavity is the Schur-Ostrowski condition (e.g. \cite{Schur1923-di,Marshall1979-lp}), which states that a function $f$ is Schur-concave if it is symmetric (invariant to reordering of arguments) and if:
\begin{eqn}
    (x_{i}-x_{j})\left({\frac {\partial f}{\partial x_{i}}}-{\frac {\partial f}{\partial x_{j}}}\right)\leq 0
\end{eqn}
for all $x$ containing non-negative elements.
Clearly $\dpp_k(\Diag(\lambda))$ is symmetric from its definition.
\newcommand{\xab}{x \cup \lrb{a,b}}
Consider now a vector of non-negative elements $x$ and non-negative numbers $a$ and $b$.
It is immediate in general that $\pdv{e_{k}(x \cup \lrb{a})}{a} = e_{k-1}(x)$,
so taking the derivative of $e_{k+1}(\xab) / e_k(\xab)$ with respect to $a$ yields that
\begin{eqn}
    \pdv{\frac{e_{k+1}(\xab)}{e_k(\xab)}}{a} = \frac{e_k(\xab) e_k(x \cup \lrb{b}) - e_{k+1}(\xab) e_{k-1}(x \cup \lrb{b})}{e_k(\xab)^2}
\end{eqn}
so that
\begin{eqn}
    \pdv{\frac{e_{k+1}(\xab)}{e_k(\xab)}}{b} & - \pdv{\frac{e_{k+1}(\xab)}{e_k(\xab)}}{a} \\ &= \tfrac{
        e_k(\xab) e_k(x \cup \lrb{a}) - e_{k+1}(\xab) e_{k-1}(x \cup \lrb{a})
        e_k(\xab) e_k(x \cup \lrb{b}) - e_{k+1}(\xab) e_{k-1}(x \cup \lrb{b})
    }{e_k(\xab)^2} \\
    &= (a-b) \frac{
        e_k(\xab) e_{k-1}(x) - (a-b) e_{k+1}(\xab) e_{k-2}(x)
    }{e_k(\xab)^2}
    \req{schur-long}
\end{eqn}
where the final line follows from isolation of the symmetric polynomial terms dependent on $a$ and $b$.
Considering the numerator of \eq{schur-long}, we see that the Schur-Ostrowski condition amounts to proving that:
\begin{eqn}
    e_{k}(\xab) e_{k-1}(x)  \geq e_{k+1}(\xab) e_{k-2}(x)
    \req{schur-numerator}
\end{eqn}
for all sets $x$ and all positive numbers $a,b$ and $k \geq 2$.
We shall use the basic identity:
\begin{eqn}
    e_k(\xab) = e_k(x) + a e_{k-1}(x) + b e_{k-1}(x) + a b e_{k-2}(x)
\end{eqn}
for any integer $k \geq 2$.
Substitution of this identity into \eq{schur-numerator} yields that:
\begin{eqn}
    e_{k}(\xab) e_{k-1}(x) - e_{k+1}(\xab) e_{k-2}(x) \\ = (a+b) (e_{k-1}(x)^2 - e_k(x)  e_{k-2}(x)) + (e_{k-1}(x)  e_k(x) - e_{k+1}(x)  e_{k-2}(x))
\end{eqn}
Let us show that each of these two terms are non-negative.
The first term is non-negative from Newton's inequalities (Theorem~\ref{th:newton}).
The second term is non-negative as well:
\begin{eqn}
    e_{k-1}(x)  e_k(x) - e_{k+1}(x)  e_{k-2}(x) &= e_{k-1}(x)  e_k(x) - e_{k+1}(x) e_{k-1}(x) e_{k-2}(x) / e_{k-1}(x) \\
    &\geq e_{k-1}(x)  e_k(x) - e_k^2(x) e_{k-2}(x) / e_{k-1}(x) \\
    &= \frac{e_k(x)}{e_{k-1}(x)} (e_{k-1}(x)^2 - e_k(x) e_{k-2}(x)) \\
    &\geq 0
\end{eqn}
by, again, two applications of Newton's inequalities.
Thus the Schur-Ostrowski criterion is fulfilled, and $\dpp_k(\Diag(\lambda))$ is Schur-convex.
\end{proof}

\subsection{Subadditivity \label{a:dpp-subadditive}}

For our linear programming bounds to relate DPP sampling to nuclear maximization, we must prove the subadditivity of $k$-DPP.
Thankfully, this follows as a straightforward corollary of the following result from the literature:

\begin{lemma}[Superadditivity of ratios of symmetric polynomials \cite{Marcus1957-cl}]
    Let $x$ and $y$ each be vectors of $n$ non-negative real numbers.
    Then, for $k = 0, \dots, n$, supposing that there are at least $k$ nonzero elements in each of $x$ and $y$:
    \begin{eqn}
        \frac{e_{k+1}(x+y)}{e_{k}(x+y)} \geq \frac{e_{k+1}(x)}{e_{k}(x)} + \frac{e_{k+1}(y)}{e_{k}(y)}
    \end{eqn}
    \label{lem:sym-superadditive}
\end{lemma}
See Theorem 1 of \cite{Marcus1957-cl} or Theorem 16 of \cite{Beckenbach2012-lx} for full proof.
(The standard proof relies on induction over increasing $k$, with $k=1$ provable as the base case.)
Using this result and standard majorization inequalities for the sum of Hermitian matrices, we can derive the following relation:
\begin{corollary}[Subadditivity of $k$-DPP]
    Let $A$ and $B$ each be symmetric positive semidefinite matrices of size $n \times n$. Then for $k = 0, \dots, n$:
    \begin{eqn}
        \dpp_k(A) + \dpp_k(B) \geq \dpp_k(A + B)
    \end{eqn}
    \label{cor:dpp-subadditive}
\end{corollary}
\begin{proof}
    Let $\lambda_A$, $\lambda_B$, and $\lambda_{A+B}$ denote the eigenvalues of $A$, $B$, and $A+B$ respectively, each sorted in descending order.
    Then $\lambda_{A+B} \prec \lambda_A + \lambda_B$ (e.g. \cite{Marshall1979-lp}, Theorem~G.1), and so:
    \begin{eqn}
        \dpp_k(A+B) \leq \dpp_k(\Diag(\lambda_A + \lambda_B)) \leq \dpp_k(\Diag(\lambda_A)) + \dpp_k(\Diag(\lambda_B)) = \dpp_k(A) + \dpp_k(B)
    \end{eqn}
    with the first inequality is from the Schur-convexity of $\dpp_k$ (Theorem~\ref{th:schur-convexity})
    and the second from Lemma~\ref{lem:sym-superadditive}.
\end{proof}

\section{Proofs of theoretical error bounds \label{a:error-bounds}}

\subsection{Chebyshev approximation \label{s:chebyshev-analysis}}

In this section, we consider the Chebyshev interpolation of the inverse square root function, showing that the number of matvecs required for our matrix-free approach to inverse Laplacian reduction requires only $O(\sqrt{\kappa} \log(\kappa))$ matvecs for fixed relative error given preconditioned condition number $\kappa$.
This characterizes the computational cost of our diagonal estimation used for inverse Laplacian reduction (Section~\ref{s:sqrt-alg}). 

\begin{lemma}
    Consider Chebyshev interpolation of the function $f(x) = \lrp{\frac{a+b+(b-a)x}{2}}^\Mh$ on the interval $[ \Minus 1, 1]$, and let $f_n(x)$ be its Chebyshev interpolant of order $n$.
    Then there exists $n^* = n^* (\ve, \kappa)$ satisfying
    \begin{eqn}
        n^* \lesssim \frac{1}{2} \sqrt{\kappa } \log ( \kappa / \ve )
    \end{eqn}
    in the $\ve \rightarrow 0^+$ and $\kappa \rightarrow \infty$ limit, such that 
    \begin{eqn}
        \req{chebrel}
        \frac{ \max_{x \in [-1,1]} \  \Abs{f(x) - f_{n}(x)} } { \min_{x \in [-1,1]}  \ \Abs{f(x)} } \leq \ve
    \end{eqn}
    for any $n \geq n^*$.
    \label{lem:chebyshev}
\end{lemma}
\begin{proof}
    We shall follow the presentation of \cite{Trefethen2020-cm} in the preliminaries.
    Note that $f(\Mo) = a^{\Mh}$ and $f(+1) = b^{\Mh}$, and therefore $\min_{x \in [-1,1]} \ \Abs{f(x)} = b^\Mh$.
    Next, $f(x)$ is analytically continuable to the Bernstein ellipse $E_\rho$ where 
    \begin{eqn}
        0 < \rho < \lrp{\frac{(a+b)^2}{(b-a)^2}-1}^{1/2}-\frac{a+b}{a-b}
    \end{eqn}
    On any such ellipse, 
    \begin{eqn}
        |f(x)| \leq M \Eq 2 \lrp{\frac{\rho}{a (\rho +1)^2-b (\rho -1)^2}}^{1/2}
    \end{eqn}
    Then by the guarantees of Chebyshev projection (\cite{Trefethen2020-cm}, Theorem~8.2), the Chebyshev interpolant $f_n$, $n \geq 1$, of $f$ satisfies
    \begin{eqn}
        \Abs{f(x) - f_n(x)} \leq \frac{4 M \rho^{\Minus s}}{\rho-1 } 
    \end{eqn}
    Next, consider the choice of $\rho = (1 - (2n)^\Mo) \frac{2 \sqrt{a b}+a+b}{b-a}$, which gives the tightest bound asymptotically for $n \rightarrow \infty$. 
    Upon substitution of this choice, one has that:
    \begin{eqn}
        \frac{|f(x) - f_n(x)|}{b^{-1/2}} \leq -\frac{\left(\sqrt{\kappa }-1\right) 2^{n+\frac{7}{2}} n \left(\frac{\left(\sqrt{\kappa }+1\right)^2 (2 n-1)}{(\kappa -1) n}\right)^{-n}}{\left(\sqrt{\kappa }-4 n+1\right) \sqrt{\frac{\kappa +\sqrt{\kappa } (2-8 n)+1}{\kappa  (1-2 n) n}}}
        \req{chebyshev-exact-bound}
    \end{eqn}
    where $\kappa = b/a$.
    In the $\ve \rightarrow 0$ limit, we must take $n \rightarrow \infty$, so we can approximate the relative error bound in Chebyshev interpolation as
    \begin{eqn}
        \frac{|f(x) - f_n(x)|}{b^\Mh}
        \sim  \tilde{\ve} \Eq
        \sqrt{2} \left(\sqrt{\kappa }-1\right) \sqrt[4]{\kappa } \sqrt{n} \left(\sqrt{\kappa }+1\right)^{-2 n} (\kappa -1)^n
    \end{eqn}
    Let $n^*$ be the minimum value of $n$ such that maximum relative error $\ve$ is assured.
    Then by solving the right-hand side above, we see that can achieve \eq{chebrel} by taking 
    \begin{eqn}
        n^* \lesssim -\frac{W_\Mo\left(\frac{\ve ^2 \left(\log (\kappa -1)-2 \log \left(\sqrt{\kappa }+1\right)\right)}{\left(\sqrt{\kappa }-1\right)^2 \sqrt{\kappa }}\right)}{4 \log \left(\sqrt{\kappa }+1\right)-2 \log (\kappa -1)}
    \end{eqn}
    using Lambert's W function \cite{Corless1996-pe}.
    Asymptotic expansion yields the desired result.
\end{proof}

In practice, for a given $\kappa$ and error requirement $\ve$, we calculated $n^*$ via the above formula and performed a bisection search of \eq{chebyshev-exact-bound} afterwards to fine-tune $n$ to be as low as possible.
The search yielded only minor changes to the analytical formula for $n^*$ in our examples.

Let $B = R^{+} L R^{+\t}$ for short, and let $v$ be an arbitrary nonzero vector. We will explain now how to perform the approximate matvec $B^{+/2} v$. Given an expansion of our approximating polynomial $f_n = \sum_{k=0}^{n} c_k T_k $ in terms of Chebyshev polynomials $T_k$, we approximate $B^{+/2} v \approx f_n (B) v = \sum_{k=0}^{n} c_k T_k(B) v$. The matvecs $T_k (B) v$, $k=0, \ldots , n$  can all be constructed using the three-term recurrence for the Chebyshev polynomials, requiring $n$ total matvecs by $B$.

The relative error $\Vert B^{+/2} v - f_n (B) v \Vert / \Vert  B^{+/2} v \Vert$ of the matvec in the Euclidean norm is bounded by 
$\ve$, provided $n$ is taken to be at least as large as $n^* (\ve, \kappa)$ as in the statement of Lemma~\ref{lem:chebyshev}. In turn the relative error in the infinity norm is bounded by $\ve \sqrt n$.

\subsection{Randomized diagonal estimation \label{a:concentration-proofs}}

In this section we provide proofs of the concentration achieved by randomized diagonal estimation within our matrix-free approach (Section~\ref{s:concentration}).

\DiagonalApprox*
\begin{proof}
    For each $i=1,\ldots ,n$, the $i$-th entry $\frac{1}{z} [C Z Z^\t C^\t]_{i,i}$ of the vector $\frac{1}{z} 
     \Diag(C Z Z^\t C^\t)$ can be viewed precisely as the $z$-shot Hutchinson trace estimator for the trace $X_{i,i} = \Tr[ Y^\t e_i e_i^\t Y]$. Then Lemma 2.1 of \cite{Meyer2021-pn}, which is a concentration bound for the Hutchinson trace estimator, together with the union bound over $i=1,\ldots ,n$, guarantees that there exists a universal $c$ such that for $\epsilon > 0$,  $\delta \in (0,1/2]$, and $z \geq c \log(n/\delta) / \epsilon^2$, the following inequality holds:
    \begin{eqn}
    \left\vert \,  \frac{1}{z}  [C Z Z^\t C^\t]_{i,i} - X_{i,i} \, \right\vert \leq \epsilon \, \Vert Y^\t e_i e_i^\t Y \Vert_F
    \req{hutchbound}
    \end{eqn}
    for all $i=1,\ldots, n$. Note that the Frobenius norm can be explicitly computed as
    \begin{eqn}
    \Vert Y^\t e_i e_i^\t Y \Vert_F^2 = \Tr[ Y^\t e_i e_i^\t Y Y^\t e_i e_i^\t Y ] = (e_i^\t Y Y^\t e_i) (e_i^\t Y Y^\t e_i) = (X_{i,i})^2.
    \end{eqn}
    Then it follows from \eq{hutchbound} that 
    \begin{eqn}
    \left\vert \, \frac{1}{z}  [C Z Z^\t C^\t]_{i,i} - X_{i,i} \, \right\vert \leq \epsilon \, X_{i,i},
    \end{eqn}
    from which the lemma follows.
\end{proof}

\StochasticRatio*
\begin{proof}
    For any given index $\j$, denote the left-hand maximand as $\tilde{g}_{\j}$ and the right-hand maximand as $g_{\j} \Eq \max_{\j} \frac{({K}^2)_{\j,\j}}{ K_{\j,\j}}$.
    Let $g_\mathrm{opt}$ be $\max_{\j} g_{\j}$ and $\tilde{g}_\mathrm{opt} = \max_{\j;g_{\j} = g_\mathrm{opt}} \tilde{g}_{\j}$.
    Then failure may occur only if, for any $\j$, $\tilde{g}_{\j} \geq \tilde{g}_\mathrm{opt}$ and $g_{\j} < \Orand^\Mo g_\mathrm{opt}$.
    Suppose we apply Theorem~\ref{th:trace-concentration} with $\epsilon \Eq \frac{\Rand-2 \sqrt{\Rand+1}+2}{\Rand}$.
    Taking a union bound over the numerator and denominator of the compared ratios (4 terms), then with at least probability $1-4\delta$, $\tilde{g}_\mathrm{opt} \geq \tfrac{1-\epsilon}{1+\epsilon} g_\mathrm{opt}$ and $\tilde{g}_{\j} \leq \tfrac{1+\epsilon}{1-\epsilon} g_{\j}$ for all $\j$.
    Equivalently, $\tilde{g}_{\j} \leq \tilde{g}_\mathrm{opt}$ for all $\j$ such that $g_{\j}/g_\mathrm{opt} < \tfrac{(1+\epsilon)^2}{(1-\epsilon)^2} = \Orand^\Mo$.
    Finally, $\epsilon \geq (3 - 2 \sqrt{2})\Rand$ for $\Rand \in (0, 1]$, so substituting the bound $\epsilon^{\Minus 2} \leq (17+2 \sqrt{2}) / \Rand^2$ into the $z$ requirement from Theorem~\ref{th:trace-concentration} yields the desired result.
\end{proof}

\subsection{Linear programming relations \label{a:lp-relations}}


As motivation for our general linear programming bound, we will recapitulate the original linear program posed in \cite{Nemhauser1978-ci} to bound the performance of greedy optimization of submodular functions (cf. Definition~\ref{def:submodular} in the main text):

\begin{restatable}[Submodular inequalities \cite{Nemhauser1978-ci}]{lemma}{SubmodularInequalities}
    For a submodular function $\subm$ and $t = 1, \dots, k$, consider the greedily chosen subset $\Gg_t$, its yielded gain $g^{(0)}_t$, such that the final objective value is $\sum_{t=1}^k g^{(0)}_t$, and the optimal subset $\Opt_s$ of size $s \leq k$ (Definition~\ref{def:greedy}).
    Then:
    \begin{eqn}
        \subm(\Opt_s) \leq 
        \sum_{i=1}^{t-1} g^{(0)}_t + s g^{(0)}_t \qquad t=1, \dots, k
        \req{submodular-example}
    \end{eqn}
    \label{lem:submodular-program}
\end{restatable}
\begin{proof}
    One may obtain \eq{submodular-example} by the following logic, in short.
    Consider each set $\Gg_{t-1}$ of $t-1$ greedily selected elements versus the optimal set $\mathcal{O}_s$.
    Then there are at most $s$ elements in the optimal set which are not in $\Gg_{t-1}$, trivially.
    Considering adding each of those $s$ elements into $\Gg_{t-1}$. Then:
    \begin{eqn}
        \subm(\Gg_{t-1} \cup \Opt_s) &\leq \subm(\Gg_{t-1}) + \sum_{\j \in \mathcal{O_s}} \subm(\Gg_{t-1} \cup \lrb{\j}) - \subm(\Gg_{t-1}) \\
        &\leq \subm(\Gg_{t-1}) + s \max_{\j \notin \Gg_{t-1}} \subm(\Gg_{t-1} \cup \lrb{\j}) - \subm(\Gg_{t-1}) \\
        &= \subm(\Gg_{t-1}) + s (\subm(\Gg_{t}) - \subm(\Gg_{t-1})) \\
        &= \sum_{i=1}^{t-1} g^{(0)}_t + s g^{(0)}_t
    \end{eqn}
    where the first line is from submodularity \eq{submodular-def} and the last line is from the definition of $g^{(0)}$.
\end{proof}


Inspired by the form of Lemma~\ref{lem:submodular-program} and the associated bounds of \cite{Nemhauser1978-ci}, we next consider solution of the following generic linear program:

\GeneralSubmodularityConstraint*
\begin{proof}
We shall establish this inequality by considering the minimization of $\One^\t \Y$ over all vectors $\Y = (y_1 ,\ldots , y_k )^\t$ which satisfy \eq{inequalities}.
Collecting the coefficients of \eq{inequalities} in matrix form, we then consider the solution of the linear program:
\begin{eqn}
    \min_\Y  \ \One^\t \Y \text{ such that } M \Y \geq \One \\
\end{eqn}
where $M$ is the lower triangular $k \times k$ matrix:
\begin{eqn}
    M \Eq \Lhs^\Mo \cdot \begin{pmatrix} 
        f_1  \\
        1 & f_2  \\
        \vdots & \ddots & \ddots  \\
        1 &  \cdots & 1  & f_k
    \end{pmatrix}
    \req{constraint-matrix}
\end{eqn}
This problem is reducible in its dual form (e.g. \cite{Matousek2007-fx}):
\begin{eqn}
    \max_x \ \One^\t x \text{ such that } M^\t x = \One, \  x \geq 0
\end{eqn}
which immediately yields the bound that $\One^\t \Y \geq (1 - \One^\t M^\Mo \One) z$.
In particular, backward substitution of the triangular matrix $M$ explicitly shows that:
\begin{eqn}
    \frac{\One^\t M^\Mo \One}{\Lhs} = 1 - \prod_{i=1}^k (1 - f_i^\Mo) > 1 - \exp \Big( \Minus \sum_{i=1}^k f_i^\Mo \Big)
\end{eqn}
where the right-hand side follows from $1-c^\Mo < e^{\Minus c^\Mo}$ for $c > 0$, yielding \eq{submodular-general-bound}.
\eq{simple-bound} proceeds by simple substitution.
\end{proof}


One advantage to the linear programming formulation is that additional constraints may be added to tighten the results.
We give proofs for two corollaries showing how this may be done, targeting scenarios in which column selection is performed with simple stopping conditions (rather than fixed $k$ in advance):

\AccumulatedConstraint*
\begin{proof}
    We again consider a linear programming formulation:
    \begin{eqn}
        \min_{\Y} b^\t \Y \text{ such that } M \Y \geq c
    \end{eqn}
    where we define the following matrices and vectors:
    \begin{equation}
        M \Eq \Lhs^\Mo \cdot \begin{pmatrix} 
            f_1 &  \\
            1 & f_2 &  \\
            \vdots & \ddots   & \ddots &  \\
            \vdots &  \ddots & \ddots  & f_{k}  \\
            1 &  \cdots & \cdots & 1  & f_{k+1}  \\
            \alpha & \cdots & \cdots & \alpha & \Mo 
        \end{pmatrix} \quad 
         \Y \Eq \begin{pmatrix} \Y_1 \\ \Y_2 \\ \vdots \\ \Y_{k+1} \end{pmatrix} \quad
        b \Eq \begin{pmatrix} 1 \\ \vdots \\ 1 \\ 0 \end{pmatrix} \quad
        c \Eq \begin{pmatrix} 1 \\ \vdots \\ \vdots \\ 1 \\ 0 \end{pmatrix}
        \req{constraint-matrix-sum}
    \end{equation}
    To make our dual system fully determined, let $\phi \in \RR$,  $\bar b$ be the result of appending $\phi$ to $b$, and let $\bar M$ be the result of appending column $\begin{pmatrix}0 & \cdots & 0& 1 \end{pmatrix}^\t$  to $M$.
    Then $\bar M$ is square and lower triangular, and our program is posable in dual form as:
    \begin{eqn}
        \max_x c^\t x \text{ such that } \bar M^\t x = \bar b, x \geq 0
    \end{eqn}
    whence we deduce that $c^\t x' = \bar b^\t \bar M^\Mo c$ for maximizing solution $x'$.
    Backward substitution yields that:
    \begin{eqn}
        (\bar M^\Mt \bar b)_i = \Lhs \cdot \begin{cases}
                f_i^\Mo(1 -  (f_{k+1}^\Mo + \alpha) \phi) \prod_{j=i+1}^k (1-f_j^\Mo) & i \leq k \\
                \phi f_{k+1}^\Mo & i = k+1 \\
                \phi & i = k+2
        \end{cases}
    \end{eqn}
    Our additional constraint will be tight if $(\bar M^{\Minus \t} \bar b)_1 = 0$.
    Solving for $\phi'$ to satisfy this condition, we find that:
    \begin{eqn}
        \phi' &= \lrp{\alpha +f_{k+1}^\Mo}^\Mo \\
        \frac{\One^\t \Y}{\Lhs} &\geq \bar b^\t \bar M^\Mo c \; \big\rvert_{\phi=\phi'} = (1 + \alpha f_{k+1})^\Mo
    \end{eqn}
    and \eq{submodular-sum-bound} is yielded via rearrangement.
\end{proof}

\InitialConstraint*
\begin{proof}
    We again consider a linear programming formulation:
    \begin{eqn}
        \min_y \ b^\t \Y \text{ such that } M {\Y} \geq c
    \end{eqn}
    where we define the following matrices and vectors:
    \begin{eqn}
        M \Eq \Lhs^\Mo \cdot \begin{pmatrix} 
            f_1 &  \\
            1 & f_2 &  \\
            \vdots & \ddots   & \ddots  \\
            \vdots &  \ddots & \ddots  & f_{k}  \\
            1 &  \cdots & \cdots & 1  & f_{k+1}  \\
            \beta & 0 & \cdots & 0 & \Mo 
        \end{pmatrix}  \quad 
         \Y \Eq \begin{pmatrix} \Y_1 \\ \Y_2 \\ \vdots \\ \Y_{k+1} \end{pmatrix} \quad
        b \Eq \begin{pmatrix} 1 \\ \vdots \\ 1 \\ 0 \end{pmatrix} \quad
        c \Eq \begin{pmatrix} 1 \\ \vdots \\ \vdots \\ 1 \\ 0 \end{pmatrix}
        \req{constraint-matrix-first}
    \end{eqn}
    Again, let $\phi \in \RR$,  $\bar b$ be the result of appending $\phi$ to $b$, and let $\bar M$ be the result of appending column $\begin{pmatrix}0 & \cdots & 0& 1 \end{pmatrix}^\t$  to $M$.
    Then our program is posable in dual form as:
    \begin{eqn}
        \max_x \ c^\t x \text{ such that } \bar M^\t x = \bar b, x \geq 0
    \end{eqn}
    whence we deduce that $c^\t x' = \bar b^\t \bar M^\Mo c$ for maximizing $x'$.
    Backward substitution yields that:
    \begin{eqn}
        (\bar M^\Mo c)_i = \Lhs \cdot \begin{cases}
            \Minus\beta \phi f_i^\Mo \delta_{i,1} + f_i^\Mo (1 - \phi f_{k+1}^\Mo) \prod_{k=i+1}^k (1 - f_k^\Mo) & i \leq k \\
            \phi f_{k+1}^\Mo & i = k+1 \\
            \phi & i = k+2
        \end{cases}
    \end{eqn}
    Analogously, our additional constraint will be tight if $(\bar M^\Mo \bar b)_1 = 0$.
    Solving for $\phi'$ to satisfy this condition, we find that:
    \begin{eqn}
        \phi' &= \lrp{ \beta \prod_{i=2}^k (1-f_i^\Mo)^\Mo +f_{k+1}^\Mo }^\Mo \\
        \frac{\One^\t \Y}{\Lhs} &\geq \bar b^\t \bar M^\Mo c \; \big\rvert_{\phi=\phi'} = 
        1-\frac{1}{\prod_{j=2}^k (1-f_j^\Mo)^\Mo+\lrp{\beta  f_{k+1}}^\Mo}
    \end{eqn}
    and \eq{submodular-first-bound} is yielded via rearrangement.
\end{proof}

Given the above generic bounds, in the following corollary, we relate our linear programs to the maximization of submodular functions, generally extending the results of \cite{Nemhauser1978-ci}.
While we do not address function that are strictly submodular (over all sets, including $\emptyset$) in the current work, this extension is immediate and may have additional use-cases:

\begin{corollary}[Approximate maximization of submodular functions]
    Given a submodular function $\subm$ (Definition~\ref{def:submodular}), consider $\Ge_t$, $\G_t$, and $\Opt_s$ as specified in Definition~\ref{def:greedy}.
    Then:
    \begin{eqn}
        1 - \frac{ \sum_{i=1}^k y_i }{\subm(\Opt_s)} \leq \prod_{i=1}^k \lrp{1 - \frac{1}{\Orand s}}^k < \exp \Big( \Minus \frac{k}{\Orand s} \Big)
    \end{eqn}
    Next suppose that we have the constraint that $\G_k \leq \alpha \,  \subm(\Ge_{k-1})$, i.e., that the gain at iteration $k$ is less than $\alpha$ times the prior accumulated gain.
    \begin{eqn}
        1 - \frac{  \sum_{i=1}^k y_i  }{\subm(\Opt_s)} \leq \frac{\Orand \alpha s}{1 + \Orand \alpha s}
    \end{eqn}
    Suppose that instead, we have the constraint that $\G_k \leq \beta \subm(\Ge_{1})$, i.e., that the gain at iteration $k$ is less than $\beta$ times the initial gain $\G_1$.
    Then:
    \begin{eqn}
        1 - \frac{  \sum_{i=1}^k y_i  }{\subm(\Opt_s)} \leq \frac{\Orand \beta s}{\frac{\beta (\Orand s)^k}{(\Orand s-1)^{k-1}}+1} < \frac{\Orand \beta s}{1 + \Orand \beta s e^{(k-1)/(\Orand s)}}
    \end{eqn}
\end{corollary}
\begin{proof}
    These results are immediate upon substitution of  $\Lhs = \subm(\Opt_s)$ and $f_1 = \cdots = f_k = \Orand s$ into \eq{submodular-general-bound}, \eq{submodular-sum-bound}, and \eq{submodular-first-bound} of Theorem~\ref{th:submodular}, Corollary~\ref{cor:accumulated-constraint}, and Corollary~\ref{cor:initial-constraint} (respectively).
\end{proof}

\subsection{Error bounds for kernel approximation \label{a:kernel-approximation-error}}

In this section, we give proofs for our stated error bounds relating nuclear maximization to DPP sampling and spectral bounds.
First, we give an expanded set of inequalities bounding nuclear maximization to DPP sampling:

\GreedyDPP*
\begin{proof}
    By substitution of Theorem~\ref{th:submodular-dpp} into Theorem~\ref{th:submodular}:
    \begin{eqn}
        1 - \frac{\Lk_K(\Ge_k)}{\dpp_s(K)} &\leq \lrp{1 - \frac{1}{\Orand s}}^k < e^{\Minus k / (\Orand s)}
    \end{eqn}
    If it is guaranteed that $\G_k \leq \Lk_K(\Ge_k)$, i.e. as a stopping condition for the algorithm, then by Corollary~\ref{cor:accumulated-constraint}:
    \begin{eqn}
        1 - \frac{\Lk_K(\Ge_k)}{\dpp_s(K)} &\leq \frac{\Orand \alpha s}{1 + \Orand \alpha s} 
    \end{eqn}
    On the other hand, if it is guaranteed that $\G_k \leq \beta \G_1$, then from Corollary~\ref{cor:initial-constraint}:
    \begin{eqn}
        1 - \frac{\Lk_K(\Ge_k)}{\dpp_s(K)} &\leq \frac{\Orand \beta s}{1 + \Orand \beta s e^{(k-1)/(\Orand s)}} 
    \end{eqn}
\end{proof}


Next, we provide a proof for the multiplicative error bound given in Section~\ref{s:kernel-error}:

\NuclearRelativeError*
\begin{proof}
    We start from the two inequalities
    \begin{eqn}
        \dpp_s(K) &\geq \Tp{r}{K} - \frac{r}{s-r+1}(\Tr[K] - \Tp{r}{K})
    \end{eqn}
    from Lemma~\ref{lem:dpp-error} (using the Schur-convexity of $\dpp_s$, Theorem~\ref{th:schur-convexity}) and
    \begin{eqn}
        \Lk(\Ge_k) &\geq \lrp{1 - e^{\frac{\Minus k}{\Orand s}}} \dpp_s(K)
    \end{eqn}
    from Corollary~\ref{cor:dpp-greedy-bound}.
    Combining lines, we find that:
    \begin{eqn}
        1 - \frac{\Lk_K(\Ge_k)}{\Tp{r}{K}} \leq \tilde \Rel \Eq e^{\frac{\Minus k}{\Orand s}} + \lrp{1 - e^{\frac{\Minus k}{\Orand s}}} \frac{r}{s-r+1} \nu
        \label{eq:866}
    \end{eqn}
    where $s$ is the number of $s$-DPP reference columns, and $r \leq s \leq k$.
    Our first task is to find the value of $s$ which yields the tightest bound on $\tilde \Rel$.
    Assuming small enough $\tilde \Rel$, $s \gg r$, so:
    \begin{eqn}
        \tilde \Rel \approx e^{\frac{\Minus k}{\Orand s}} + \lrp{1 - e^{\frac{\Minus k}{\Orand s}}} \frac{r \nu}{s}
        \req{approx-eps}
    \end{eqn}
    Minimizing \eq{approx-eps} with respect to $s$ and making use of Lambert's $W$ function \cite{Corless1996-pe}, we find that:
    \begin{eqn}
        \hat{s} &= \frac{\nu  k r}{k + \Orand r \nu (1 - W_{0}(e^{\frac{k}{\Orand r \nu}+1}))} 
        \approx \frac{k \nu  r}{k + \Orand \nu  r \left(1-\left(\log \left(e^{\frac{k}{\Orand \nu  r}+1}\right)-\log \log \left(e^{\frac{k}{\Orand \nu  r}+1}\right)\right)\right)} \\
        &= \frac{k}{\Orand \log \left(\frac{k}{\Orand r \nu}+1\right)} 
        \approx \frac{k}{\Orand \log \left(\frac{k}{\Orand r \nu}\right)}
    \end{eqn}
    where in the first approximation we have used the asymptotic expansion of $W_0(z)$ for large $z$, and we have again considered the large $k$ limit in the second approximation.
    Directly plugging in $\hat{s}$ for $s$ in \eq{866}, one now finds that
    \begin{eqn}
        1 - \frac{\Lk_K(\Ge_k)}{\Tp{r}{K}} \leq \Rel \lesssim \nu  r \left(\frac{1+\Rand}{k} + \frac{1-\frac{\Orand \nu  r}{k}}{\frac{k}{\Orand \log \left(\frac{k}{\Orand \nu  r}\right)}+1-r}\right) 
        \sim \frac{\nu  r}{\frac{k}{\Orand \log \left(\frac{k}{\Orand \nu  r}\right)}-r+1}
    \end{eqn}
    with the second approximation following from assumption of large $k$.
    Inverting this relation to find $k^*$ in terms of $\Rel$, one finds that:
    \begin{eqn}
        k^* &\lesssim \Minus \Orand \left(\frac{\nu  r}{\Rel }+r-1\right) W_\Mo\left(-\frac{\nu  r}{\frac{\nu  r}{\Rel }+r-1}\right) \\
        &\sim \Orand \left(r \nu / \Rel+ r - 1\right) \left(\log \left(\frac{r-1}{r \nu}+\frac{1}{\Rel }\right)+\log \log \left(\frac{r-1}{r \nu}+\frac{1}{\Rel }\right)\right) \\
        &\sim \Orand \left(r \nu/\Rel +r-1\right) \log \left(\frac{r-1}{r \nu}+\frac{1}{\Rel }\right)
    \end{eqn}
    where we have used the asymptotic expansion of $W_\Mo(z)$ for small negative $z$ in the second line and dropped the $\log \log$ dependence in the final line.
\end{proof}

By rearrangement of Theorem~\ref{th:greedy-general-bound}, we can form an additive $r,\ve$ bound for nuclear maximization:

\NuclearAdditiveError*
\begin{proof}
    We may trivially between the relative error bound $\Rel$ and the additive error bound $\ve$ by multiplying both sides of the relative error bound by $\nu^\Mo$ and rearranging:
    \begin{eqn}
        1 - \frac{\Lk_K(\Ge_k)}{\Tp{r}{K}} &\leq \omega \\
        \frac{\Tp{r}{K} - \Lk_K(\Ge_k)}{\Tr[K] - \Tp{r}{K}} &\leq \nu^\Mo \omega \\
        \Tp{r}{K} - \Lk_K(\Ge_k)  &\leq \nu^\Mo \omega (\Tr[K] - \Tp{r}{K}) \\
        \Tr[K] - \Lk_K(\Ge_k)  &\leq (1 + \nu^\Mo \omega) (\Tr[K] - \Tp{r}{K})
    \end{eqn}
    That is, if an algorithm yields an $r,\ve$ bounded solution, then that solution is within relative error $\ve \nu$ of the partial trace $\Tp{r}{K}$.
    Substituting $\omega = \ve \nu$ into \eq{general-relative-bound} immediately yields \eq{general-re-bound}.
\end{proof}

\subsection{Error bounds for graph Laplacians \label{a:laplacian-error}}

Here we provide proofs of our error bounds for inverse Laplacian reduction using nuclear maximization (Section~\ref{s:laplacian-analysis}).
To begin, we provide a proof that our objective $\Lap_L^h$ is submodular on non-empty sets:

\LaplacianSubmodularity*
\begin{proof}
We will proceed in proof by studying the random walk associated to graph Laplacian $L$.
Note that, for any non-empty set of indices $\I$:
\begin{eqn}
    (L_{\bar \I, \bar \I})^\Mo = \int_0^\infty \exp(\Minus L_{\bar \I, \bar \I} t) \odif t
\end{eqn}
where $\bar \I = \{1 \dots n\} \setminus \I$.
The term $\exp(\Minus L_{\bar \I, \bar \I} t)_{\alpha, \beta}$ corresponds to the probability of a continuous-time random walk starting at state $\alpha$ being at state $\beta$ at time $t$, where the original random walk has been modified to be killed on all states in set $\I$.
Likewise, $(L_{\bar \I, \bar \I})^\Mo_{\alpha,\beta}$ corresponds to the expected time the original random walk started from $\alpha$ spends in $\beta$ before reaching any state in $\I$ (see Appendix~\ref{a:markov-interpretation} for more details).
(Note that this expectation converges to a finite value, given non-empty $\I$, by our assumption that the original random walk is recurrent.)

Let us consider all paths $\mathcal{P}$ from an arbitrary starting point $\alpha$ which end in $A$, and partition these paths according to the sequence of the order of the first visits to subsets $A$, $B$, and $C$.
Thus the set $\mathcal{P}$ may be partitioned into five disjoint subsets:
\begin{eqn}
    \mathcal{P} = P_{\alpha \To A}
    \cup P_{\alpha \To B \To A}
    \cup P_{\alpha \To C \To A}
    \cup P_{\alpha \To B \To C \To A}
    \cup P_{\alpha \To C \To B \To A}
\end{eqn}
where:
\begin{itemize}
	\item $P_{\alpha \To A}$ comprises the set of all paths from $\alpha$ to $A$ which do not traverse any state in $B$ or $C$.
	\item $P_{\alpha \To B \To A}$ comprises the set of all paths from $\alpha$ to $A$ which traverse $B$ but not $C$.
	\item $P_{\alpha \To C \To A}$ comprises the set of all paths from $\alpha$ to $A$ which traverse $C$ but not $B$.
	\item $P_{\alpha \To B \To C \To A}$ comprises the set of all paths from $\alpha$ to $A$ which traverse both $B$ and $C$, and in which the first visit to $B$ occurs before the first visit to $C$.
	\item $P_{\alpha \To C \To B \To A}$ comprises the set of all paths from $\alpha$ to $A$ which traverse both $B$ and $C$, and in which the first visit to $C$ occurs before the first visit to $B$.
\end{itemize}

Then using this decomposition, we may express the expected time $((L_{\bar A, \bar A})^\Mo)_{\alpha,\beta}$ as follows, (where $p(x)$ is the probability of path $x$):
\begin{eqn}
	\tfrac{h_\beta}{h_\alpha} ((L_{\bar A, \bar A})^\Mo)_{\alpha,\beta} \hspace{0.1in} = \hspace{0.1in} \begin{aligned}
	&\sum_{x \in P_{\alpha \To A}} p(x) T^\beta_{\alpha \To A}(x) \\
	+&\sum_{x \in P_{\alpha \To B \To A}} \hspace{-0.1in} p(x) \big (T^\beta_{\alpha \To B}(x) + T^\beta_{B \To A}(x) \big) \\
	+&\sum_{x \in P_{\alpha \To C \To A}} \hspace{-0.1in} p(x) \big (T^\beta_{\alpha \To C}(x) + T^\beta_{C \To A}(x) \big) \\
	+&\sum_{x \in P_{\alpha \To B \To C \To A}} \hspace{-0.2in} p(x) \big (T^\beta_{\alpha \To B}(x) + T^\beta_{B \To C}(x) + T^\beta_{C \To A}(x) \big) \\
	+&\sum_{x \in P_{\alpha \To C \To B \To A}} \hspace{-0.2in} p(x) \big (T^\beta_{\alpha \To C}(x) + T^\beta_{C \To B}(x) + T^\beta_{B \To A}(x) \big)
	\end{aligned} \hspace{0.1in} \eqqcolon \hspace{0.1in} \begin{aligned}
		 & c_a \\ 
		 + &c_b + c_c \\ 
		 + &c_d + c_e \\ 
		 + &c_f + c_g + c_h \\
		 + &c_i + c_j + c_k
	\end{aligned}
    \req{path-decomp}
\end{eqn}
where, for brevity, we have defined scalar constants $c_{:}$ as the summations represented by visual analogy.
That is, 
$c_a = \sum_{x \in P_{\alpha \To A}} p(x) T^\beta_{0 \To A}(x)$, 
and, for example,
$c_g = \sum_{x \in P_{\alpha \To B \To C \To A}} p(x) T^\beta_{B \To C}(x)$.
We may enumerate the terms in the submodularity relation \eq{submodular-trace} as follows:
\begin{enumerate}
    \item $\tfrac{h_\beta}{h_\alpha} ((L_{\Comp{A},\Comp{A}})^\Mo)_{\alpha,\beta}$ is the expected time the random walk from $\alpha$ spends in $\beta$ before reaching $A$.
    \item $\tfrac{h_\beta}{h_\alpha} ((L_{\Comp{A \cup B},\Comp{A \cup B}})^\Mo)_{\alpha,\beta}$ is the expected time the random walk from $\alpha$ spends in $\beta$ before reaching $A$ or $B$.
    \item $\tfrac{h_\beta}{h_\alpha} ((L_{\Comp{A \cup C},\Comp{A \cup C}})^\Mo)_{\alpha,\beta}$ is the expected time the random walk from $\alpha$ spends in $\beta$ before reaching $A$ or $C$.
    \item $\tfrac{h_\beta}{h_\alpha} ((L_{\Comp{A \cup B \cup C},\Comp{A \cup B \cup C}})^\Mo)_{\alpha,\beta}$ is the expected time the random walk from $\alpha$ spends in $\beta$ before reaching $A$, $B$, or $C$.
\end{enumerate}
While $((L_{\bar A, \bar A})^\Mo)_{\alpha,\beta}$ represents the most inclusive summation here, the subsequent quantities (e.g., $((L_{\Comp{A \cup B},\Comp{A \cup B}})^\Mo)_{\alpha,\beta}$) incorporate only subsets of the $c$ terms.
Using the above numbering, it may be helpful to depict the path segments which contribute to the summation in \eq{path-decomp} as follows, in visually analogous order:
\begin{eqn}
	& \xrightarrow[a]{1,2,3,4} A & (\text{from paths } P_{\alpha \To A}) \\
	& \xrightarrow[b]{1,2,3,4} B \xrightarrow[c]{1,3} A & (\text{from paths } P_{\alpha \To B \To A})  \\
	& \xrightarrow[d]{1,2,3,4} C \xrightarrow[e]{1,2} A & (\text{from paths } P_{\alpha \To C \To A}) \\
	& \xrightarrow[f]{1,2,3,4} B \xrightarrow[g]{1,3} C \xrightarrow[h]{1} A & (\text{from paths } P_{\alpha \To B \To C \To A}) \\
	& \xrightarrow[i]{1,2,3,4} C \xrightarrow[j]{1,2} B \xrightarrow[k]{1} A & (\text{from paths } P_{\alpha \To C \To B \To A}) \\
\end{eqn}
where the top indices reflect the constraint index for those paths which incorporate summations of $c$ subscripted by the bottom indices.
Considering all restrictions, the total integrals are representable as the summations:
\begin{eqn}
    \tfrac{h_\beta}{h_\alpha} ((L_{\Comp{A},\Comp{A}})^\Mo)_{\alpha,\beta} &= c_a + c_b + c_c + c_d + c_e + c_f + c_g + c_h + c_i + c_j + c_k \\
    \tfrac{h_\beta}{h_\alpha} ((L_{\Comp{A \cup B},\Comp{A \cup B}})^\Mo)_{\alpha,\beta} &= c_a + c_b + c_d + c_e + c_f + c_i + c_j \\
    \tfrac{h_\beta}{h_\alpha} ((L_{\Comp{A \cup C},\Comp{A \cup C}})^\Mo)_{\alpha,\beta} &= c_a + c_b + c_c + c_d + c_f + c_g + c_i \\
    \tfrac{h_\beta}{h_\alpha} ((L_{\Comp{A \cup B \cup C},\Comp{A \cup B \cup C}})^\Mo)_{\alpha,\beta} &= c_a + c_b + c_d + c_f + c_i
    \req{path-sums}
\end{eqn}
from which we quickly obtain that:
\begin{eqn}
    ((L_{\Comp{A},\Comp{A}})^\Mo)_{\alpha,\beta}
    + ((L_{\Comp{A \cup B},\Comp{A \cup B}})^\Mo)_{\alpha,\beta}
    + ((L_{\Comp{A \cup C},\Comp{A \cup C}})^\Mo)_{\alpha,\beta}
    + ((L_{\Comp{A \cup B \cup C},\Comp{A \cup B \cup C}})^\Mo)_{\alpha,\beta} = \tfrac{h_\alpha}{h_\beta} (c_h + c_k) \geq 0
\end{eqn}
where the last inequality holds since $c_h$ and $c_k$ are non-negative.
It is trivial now by summation over $\alpha=\beta=1, \dots, n$ that 
\begin{eqn}
    \TrL{\Comp{A}} - \TrL{\Comp{A\cup B}} - \TrL{\Comp{A \cup C}}  + \TrL{\Comp{A\cup B\cup C}} \geq 0
\end{eqn}
and thus our considered objective $\Lap_L^h$ is submodular on non-empty sets.
\end{proof}


To establish an error bound for nuclear maximization, we must first bound the value of $\Lap_L^h$ after a single iteration, which is an edge case due to the stationary mode in $L$:

\LaplacianFirstBound*
\begin{proof}
    From \eq{first-column-laplacian} we have that:
        $\Tr[(L_{\overline{ \{i\} }, \overline{ \{i\} }})^\Mo] = \Tr[L^+] + \frac{(L^+)_{i,i}}{h_i^2}$,
    so that, using $\norm{h}_2 = 1$:
    \begin{eqn}
        \min_i \left\{ \Tr[(L_{\overline{ \{i\} }, \overline{ \{i\} }})^\Mo] \right\} = \Tr[L^+] + \min_i  \left\{ \frac{(L^+)_{i,i}}{h^2_i} \right\} \leq \Tr[L^+] + \frac{\sum_i (L^+)_{i,i}}{\sum_i h^2_i} = 2 \Tr[L^+]
    \end{eqn}
    and, by the same general argument:
    \begin{eqn}
        \Tr[(L_{\overline{ \{\j\} }, \overline{ \{ \j \} }})^\Mo] = \Tr[L^+] + \frac{(L^+)_{\j,\j}}{h^2_\j} \leq \Tr[L^+] + \Orand \min_i \frac{(L^+)_{i,i}}{h^2_i} \leq (2 + \Rand) \Tr[L^+]
    \end{eqn}
    Then  \eq{one-iter-3} follows by the definition of $\Lap_L^h$. 
\end{proof}


Following from Lemma~\ref{lem:laplacian-initial-bound} and the submodularity of $\Lap_L^h$ on non-empty sets (Theorem~\ref{th:laplacian-is-submodular}), we can bound the result of nuclear maximization relative to the optimal column selection subset as follows:

\LaplacianTraceBound*
\begin{proof}
    We shall generally assume the first element $\Ge_1$ fixed and apply submodularity relations to the rest.
    For convenience, define the vector of gains as follows:
    \begin{eqn}
        g_i = \begin{cases}
            \Lap_L^h(\Ge_1) & i = 1 \\
            \Lap_L^h(\Ge_i) - \Lap_L^h(\Ge_{i-1}) & i > 1
        \end{cases}
    \end{eqn}
    Then, given that $\Lap_L^h$ is submodular over non-empty subsets, we obtain the following inequalities:
\begin{eqn}
    \Lap_L^h(\Opt_s) &\leq g_1 + \Orand s g_2 \\
    \Lap_L^h(\Opt_s) &\leq g_1 + g_2 + \Orand s g_3 \\
    \dots & \\
    \Lap_L^h(\Opt_s) &\leq \sum_{i=1}^{k-1} g_i + \Orand s g_k
    \req{laplacian-submodularity-1}
\end{eqn}
    Compared to \eq{submodular-example}, this system is missing the initial inequality $\Lap_L^h(\Opt_s) \leq \Orand s g_1$ due to the handling of the stationary state.
    Nevertheless, subtracting $g_1$ from both sides of \eq{laplacian-submodularity-1} yields an analogous linear program in $k-1$ equations (intuitively, with $\Lhs=\Lap_L^h(\Opt_s) - g_1$), yielding the bound:
    \begin{eqn}
        \Lap_L^h(\Opt_s) - g_1 - \sum_{i=2}^{k} g_i = \Lap_L^h(\Opt_s) - \Lap_L^h(\Ge_k) \leq \lrp{1-\tfrac{1}{\Orand s}}^{k-1} (\Lap_L^h(\Opt_s) - g_1)
    \end{eqn}
    The result \eq{laplacian-trace-bound} follows from 
    Lemma~\ref{lem:laplacian-initial-bound}, which shows that $g_1 \geq \Minus \Orand \Tr[L^+]$,
    \eq{laplacian-objective}, which shows that $\Lap_L^h(\Opt_s) \leq \Tr[L^+]$,
    and finally the inequality $1-c^\Mo < e^{\Minus c^\Mo}$ for all $c > 0$.
\end{proof}


\end{document}

%% file: math.tex
\usepackage[parfill]{parskip}
\usepackage[utf8]{inputenc}    
\usepackage{amsmath,amssymb,amsfonts}
\usepackage{graphicx}
\usepackage{hyperref}
\usepackage{etoolbox}
\usepackage{titling}
\usepackage[normalem]{ulem}
\usepackage{algpseudocodex}
\usepackage{algorithm}
\usepackage{mathtools}
\usepackage{scalerel}
\usepackage{derivative}
\usepackage{abstract}
\usepackage{braket}

\newcommand{\Nystrom}{Nystr\"om}

\newcommand{\pinv}{+}

\DeclareMathOperator{\nnz}{nnz}
\DeclareMathOperator{\Diag}{diag}

\DeclareMathOperator{\I}{\mathcal{I}}
\DeclareMathOperator{\J}{\mathcal{J}}

\renewcommand{\j}{\ell}

\renewcommand{\t}{{\protect \scaleobj{0.8}{\top}}}
\newcommand{\Id}{\mathbb{I}}

\newcommand{\Ex}{\mathbb{E}}

\newcommand{\mvec}{\mathrm{vec}}

\newcommand{\Tr}{\mathrm{Tr}}

\newcommand{\lrp}[1]{{\left( #1 \right) }}
\newcommand{\lrs}[1]{{\left[ #1 \right] }}
\newcommand{\lrb}[1]{{\left \{ #1 \right \} }}

\newcommand{\One}{\mathbf{1}}
\newcommand{\Zero}{\mathbf{0}}
\newcommand{\RR}{\mathbb{R}}
\newcommand{\Eq}{\coloneqq}

\newcommand{\ve}{\varepsilon}

\newenvironment{eqn}{\begin{equation}\begin{aligned}}{\end{aligned}\end{equation}\ignorespacesafterend}

\DeclareMathOperator*{\argmax}{arg\,max}

\DeclareMathOperator*{\logdet}{log\,det}

\DeclarePairedDelimiterX{\norm}[1]{\lVert}{\rVert}{#1}
\DeclarePairedDelimiter\Abs{\lvert}{\rvert}%

\newcommand{\Minus}{-}
\newcommand{\Mo}{{\Minus 1}}
\newcommand{\Mt}{{\Minus \t}}
\newcommand{\Mh}{{\Minus 1/2}}

\newcommand{\otwo}{{\odot 2}}
\newcommand{\SetDef}[2]{\lrb{#1 \mid #2}}

\newcommand{\Binomial}[2]{{#1 \choose #2}}
\newcommand{\req}[1]{\label{eq:#1}}
\newcommand{\eq}[1]{(\ref{eq:#1})}

\usepackage{amsthm}

\usepackage{framed}
\usepackage{thmtools}
\renewenvironment{leftbar}{%
  \MakeFramed {\advance\hsize-\width \FrameRestore}}%
 {\endMakeFramed}

\declaretheoremstyle[headfont=\bfseries,%
spaceabove=2pt,spacebelow=2pt,%
notefont=\bfseries,%
notebraces={}{},%
headpunct=,%
bodyfont=,%
headformat=\color{black}\NAME~\NUMBER:~{\NOTE}\hfill\smallskip\linebreak,%
preheadhook=\begin{leftbar},%
postfoothook=\end{leftbar},%
]{customDefinition}
\declaretheorem[style=customDefinition]{definition}
\declaretheorem[style=customDefinition]{theorem}
\declaretheorem[style=customDefinition, numberwithin=section]{lemma}
\declaretheorem[style=customDefinition, numberwithin=section]{proposition}
\declaretheorem[style=customDefinition, numberwithin=theorem]{corollary}

\makeatletter
\def\th@plain{%
  \thm@notefont{}
  \itshape 
}
\def\th@definition{%
  \thm@notefont{}
  \normalfont 
}
\makeatother

\makeatletter
\g@addto@macro\normalsize{%
  \setlength\abovedisplayskip{7pt}
  \setlength\belowdisplayskip{7pt}
  \setlength\abovedisplayshortskip{2pt}
  \setlength\belowdisplayshortskip{2pt}
}
\makeatother

%% file: common.tex
\usepackage[parfill]{parskip}
\usepackage[utf8]{inputenc}    
\usepackage{amsmath,amssymb,amsfonts,thmtools}
\usepackage{graphicx}
\usepackage{hyperref}
\usepackage{etoolbox}
\usepackage{titling}
\usepackage[normalem]{ulem}
\usepackage{algpseudocodex}
\usepackage{algorithm}
\usepackage{blindtext}
\usepackage{mathtools}
\usepackage{todonotes}
\usepackage{scalerel}
\usepackage{derivative}
\usepackage{pdfpages}
\usepackage{setspace}
\usepackage{wrapfig}
\usepackage{microtype}
\usepackage{booktabs}
\usepackage[margin=1cm]{caption}

\usepackage{enumitem}
\setlist{nosep}

\usepackage{abstract}

\hypersetup{
    colorlinks,
    citecolor=blue,
    filecolor=red,
    linkcolor=black,
    urlcolor=blue
}

%% file: symbols.tex
\newcommand{\tK}{\tilde{K}}
\newcommand{\Opt}{\mathcal{O}}
\newcommand{\Rand}{\zeta}
\newcommand{\Orand}{(1+\zeta)}
\newcommand{\Y}{y}
\newcommand{\Gg}{\mathcal{G}^{(0)}}
\newcommand{\G}{g^{(\Rand)}}
\newcommand{\Ge}{\mathcal{G}^{(\Rand)}}
\newcommand{\Rel}{\omega}

\newcommand{\Tp}[2]{\Tr^{(#1)} \left[#2 \right]}
\newcommand{\Comp}[1]{\overline{#1}}
\newcommand{\To}{\shortrightarrow}
\newcommand{\h}{h}
\newcommand{\dpp}{\mathcal{D}}
\newcommand{\Lk}{\mathcal{L}}
\newcommand{\Lap}{\mathcal{L}}

\newcommand{\Lhs}{\mathcal{R}}
\newcommand{\subm}{\xi}
\newcommand{\Shift}{\gamma}

\newcommand{\TrL}[1]{\Tr [(L_{#1, #1})^\Mo ]}
\newcommand{\Ij}{\I \cup \{\j\}}

%% file: algorithms/exact-kernel.tex
\begin{algorithm}[H] \begin{algorithmic}
    \setstretch{1.2}
    \Procedure{EfficientNuclearMaximization}{$K,k$}
        \State $d \gets \Diag(K)$ \Comment{Denominator term}
        \State $w \gets \Diag(K^2)$ \Comment{Numerator term, equivalent to $w \gets (K \odot K) \One$}
        \State $U \gets$ matrix of zeros of size $(k, n)$ \Comment{$n$ is the dimension of input $K$}
        \State $S \gets$ empty matrix of size $(n,k)$
        \State $\I \gets \emptyset$
        \For{$t = 1, \dots, k$}
            \State $\j \gets \argmax_{\j \notin \I} w_\j / d_\j$ \Comment{Choose maximum gain}
            \State $\I \gets \I \cup \j$
            \State $U_{t,\j} \gets 1$ \Comment{Perform normal Cholesky updates}
            \State $U_{:,\j} \gets U_{:,\j} / \sqrt{d_\j}$
            \State $S_{:,t} \gets \Minus K_{:, \I} U_{:t, \j}$
            \State $U_{:t, \Comp \I} \gets U_{:t, \Comp \I} + U_{:t, \j} S_{\Comp \I, t}^\t$
            \State $d \gets d - S_{:, t}^\otwo$ \Comment{Update denominator}
            \State $w \gets w + S_{:,t} \odot (2 S_{:, :t} S_{:, :t}^\t S_{:, t} - 2 K S_{:, t} - S_{:, t} S_{:, t}^\t S_{:, t})$ \Comment{Update numerator}
        \EndFor
    \EndProcedure
    
    \caption{Efficient nuclear maximization using upper Cholesky factorization 
        \label{alg:exact-cholesky}
    }
\end{algorithmic} \end{algorithm}

%% file: algorithms/exact-laplacian.tex
\begin{algorithm}[H] \begin{algorithmic}
    \setstretch{1.2}
    \Function{ExactLaplacianScores}{$t, w, d, g, c, y$}
        \If{$t = 1$}
            \Return $-d / h^\otwo$ 
        \Else
            \; \Return $\dfrac{w + 2 g^\Mo \cdot y \odot c + g^{\Minus 2} \cdot y^\t y \cdot y^\otwo} 
            {d + g^\Mo \cdot y^\otwo}$
        \EndIf
    \EndFunction
    \\
    \Procedure{LaplacianNuclearMaximization}{$K,h,k$}
        \State $d \gets \Diag(K)$ \Comment{Unadjusted denominator term}
        \State $w \gets \Diag(K^2)$ \Comment{Unadjusted numerator term}
        \State $U \gets$ zero matrix of size $(k, n)$ \Comment{$n$ is the dimension of input $K$}
        \State $S \gets$ empty matrix of size $(n,k)$
        \State $y \gets h$  \Comment{Projected stationary state}
        \State $c \gets$ zero vector of length $n$
        \State $g \gets 0$
        \State $\I \gets \emptyset$
        \For{$t = 1, \dots, k$}
            \State $v \gets \Call{ExactLaplacianScores}{t, w, d, g, c, y}$
            \State $\j \gets \argmax_{\j \notin \I} v_{\j}$
            \State $\I \gets \I \cup \j$
            \State $U_{t, \j} \gets 1$ \Comment{Normal Cholesky updates}
            \State $U_{:, \j} \gets U_{:, \j} / \sqrt{d_{\j}}$ 
            \State $S_{:,t} \gets \Minus K_{:, \I} U_{:t, \j}$ \Comment{Single reduced matvec by $K$}
            \State $d \gets d - S_{:,t}^\otwo$
            \State $U_{:t, \Comp \I} \gets U_{:t, \Comp \I} + U_{:t, \j} S_{\Comp \I,t}^\t$
            \State $\tau \gets U_{:t, \j}^\t h_{\I}$ \Comment{Handle stationary state}
            \State $g \gets g + \tau^2$
            \State $y \gets y + \tau S_{:,t}$
            \State $f \gets K S_{:, t} - S_{:, :t-1} S_{:, :t-1}^\t S_{:, t}$ 
            \State $c \gets c + \tau f - S_{:,t}^\t y S_{:,t}$
            \State $w \gets w + S_{:,t}^\t S_{:,t} S_{:,t}^\otwo - 2 f \odot S_{:,t}$
        \EndFor
    \EndProcedure
    \caption{Deterministic nuclear maximization for inverse graph Laplacians.
        \label{alg:exact-laplacian}
    }
\end{algorithmic} \end{algorithm}

%% file: algorithms/randomized-kernel.tex
\begin{algorithm}[H] \begin{algorithmic}
    \setstretch{1.2}
    \Function{RandomizedScores}{$\I,t,z,K,C,U,S$}
        \State $Z^{(1)}, Z^{(2)} \gets$ random i.i.d. Gaussian matrices of with $z$ columns
        \State $D \gets C Z^{(1)}$
        \State $N \gets K Z^{(2)}$
        \State \Return $\dfrac{(N + S_{:,:t-1} U_{:t-1, \I}^\t N_{\I,:})^\otwo \One}  
                        {(D + S_{:,:t-1} U_{:t-1, \I}^\t D_{\I,:})^\otwo \One}$
    \EndFunction
    \\
    \Procedure{MatrixFreeNuclearMaximization}{$K,C,k,z$}
        \State $U \gets$ matrix of zeros of size $(k,n)$ \Comment{$n$ is the dimension of input $K$}
        \State $S \gets$ empty matrix of size $(n,k)$
        \State $\I \gets \emptyset$
        \For{$t = 1, \dots, k$}
            \State $y \gets \Call{RandomizedScores}{\I,t,z,K,C,U,S}$
            \State $\j \gets \argmax_{\j \notin \I} y_{\j}$
            \State $\I \gets \I \cup \j$
            \State $U_{t,\j} \gets 1$ \Comment{Perform normal Cholesky factorization updates}
            \State $U_{:,\j} \gets U_{:,\j} / (K_{\j,\j} - S_{\j,:t-1} S_{\j,:t-1}^\t)^{1/2}$ \Comment{Accomplishable via single matvec}
            \State $S_{:,t} \gets \Minus K_{:,\I} U_{:t, \j}$  \Comment{Accomplishable via single matvec}
            \State $U_{:t,\Comp \I} \gets U_{:t,\Comp \I} + U_{:t, \j} S_{\Comp \I, t}^\t$ 
        \EndFor
    \EndProcedure
    \caption{Matrix-free nuclear maximization using upper Cholesky factorization
    \label{alg:randomized-cholesky}
    }
\end{algorithmic} \end{algorithm}

%% file: algorithms/randomized-laplacian.tex
\begin{algorithm}[H] \begin{algorithmic}
    \setstretch{1.2}
    \Function{RandomizedLaplacianScores}{$t, z, K, C, h,S, U,y, g$}
        \State $Z^{(1)}, Z^{(2)} \gets$ random Gaussian matrices of size $(n,z)$
        \State $D \gets C Z^{(1)}$
        \If{$t = 1$}
            \Return $\Minus \dfrac{D^\otwo \One}{z h^\otwo}$
        \Else
            \; \Return $\dfrac{z^\Mo ((K - S_{:,:t-1} S_{:,:t-1}^\t + g^\Mo y (y - h)^\t) Z^{(2)})^\otwo \One + g^{\Minus 2}  \odot y^\otwo} 
            {z^\Mo (D + S_{:,:t-1} U_{:t-1, \I}^\t D_{\I, :})^\otwo \One + g^\Mo \cdot y^\otwo}$
        \EndIf
    \EndFunction
    \\
    \Procedure{MatrixFreeLaplacianSubset}{$K,C,h,k,z$}
        \State $U \gets$ zero matrix of size $(k, n)$ \Comment{$n$ is the dimension of input $K$}
        \State $S \gets$ empty matrix of size $(n,k)$
        \State $y \gets h$ 
        \State $g \gets 0$
        \State $\I \gets \emptyset$
        \For{$t = 1, \dots, k$}
            \State $v \gets \Call{RandomizedLaplacianScores}{t, z, K, C, h, S, U, y, g}$
            \State $\j \gets \argmax_{\j} v_{\j}$
            \State $\I \gets \I \cup \j$
            \State $U_{t,\j} \gets 1$
            \State $U_{:,\j} \gets U_{:,\j} / \sqrt{K_{\j,\j} - S_{\j,:t-1} S_{\j,:t-1}^\t}$ \Comment{Calculate $K_{\j,\j}$ via single matvec}
            \State $S_{:,t} \gets \Minus K_{:, \I} U_{:t, \j}$ \Comment{Single reduced matvec by $K$}
            \State $U_{:\j, \Comp \I:} \gets U_{:\j, \Comp \I} + U_{:t, \j} S_{\Comp \I,t}^\t$
            \State $\tau \gets U_{:t, \j}^\t h_{\I}$
            \State $g \gets g + \tau^2$
            \State $y \gets y + \tau S_{:,t}$
        \EndFor
    \EndProcedure
    \caption{Inverse Laplacian rank reduction via matrix-free nuclear maximization
        \label{alg:randomized-laplacian}
    }
\end{algorithmic} \end{algorithm}